\documentclass[11pt]{article} 
 
\newcommand{\version}{23-7-2011}

\usepackage{amsthm,amsfonts,amsmath,amssymb}
\swapnumbers
\theoremstyle{plain}

\newtheorem{thm}{THEOREM}[section]
\newtheorem{lm}[thm]{LEMMA}

\newtheorem{prop}[thm]{PROPOSITION}

\theoremstyle{definition}

\theoremstyle{remark}

\newcommand{\upchi}{\raise1pt\hbox{$\chi$}}
\newcommand{\R}{{\mathord{\mathbb R}}}

\newcommand{\N}{{\mathord{\mathbb N}}}

\newcommand{\F}{{\mathcal{F}}}

\renewcommand{\|}{{\Vert}}
\numberwithin{equation}{section}  
\def\dd{{\rm d}}

\def\a{\alpha}
\def\b{\beta}
\def\d{\delta}
\def\e{\epsilon}

\def\g{\gamma}
\def\k{\kappa}
\def\l{\lambda}

\def\s{\sigma}

\def\L{\Lambda}

\def\FF{{\mathcal F}}
\def\II{{\mathcal I}}
\def\BB{{\mathcal B}}
\def\AA{{\mathcal A}}
\def\md{{\overline m}}
\def\RR{{\mathcal R}}

\begin{document}

\markboth{\scriptsize{CO \version}}{\scriptsize{CO \version}}

\title{\bf{ Stability of planar  fronts for a non--local phase kinetics equation with a conservation law in $D \le 3$}}
\author{\vspace{5pt} Eric A. Carlen$^{1}$ and Enza Orlandi$^{2}$
 \\
\vspace{5pt}\small{1 Department of Mathematics, Rutgers University, NJ, USA;}
\\   \vspace{5pt}\small{carlen@math.rutgers.edu}
\\
\small{2 Dipartimento di Matematica, Universit\'a degli Studi,  Roma
Tre,   P. S. Murialdo 1,
  00146 Roma, Italy;  }\\ \small {orlandi@mat.uniroma3.it}
\\}
\date{\version}
\maketitle \footnotetext    {  Work of Eric Carlen partially supported by U.S.
National Science Foundation
grant DMS 0901632   and  INDAM;     Work of Enza Orlandi partially supported by PRIN -2009-11, Universit\'a degli Studi,  Roma
Tre.      \\
\copyright\, 2010 by the authors. This paper may be reproduced, in
its entirety, for non-commercial purposes. }
 
 \begin{abstract}   
 We   consider, in a $D-$dimensional cylinder,  a non--local evolution
equation that describes the  evolution of the local magnetization in a continuum 
limit of an Ising spin system with Kawasaki dynamics and Kac potentials.
 We consider sub--critical temperatures,  for which there are two local  spatially  homogeneous  equilibria,
and  show   a local nonlinear stability result for the  minimum free
energy profiles for the magnetization at the interface between regions of these two
different local equilibrium; i.e., the planar  fronts: We show   that an initial perturbation
 of a front  that is sufficiently small in $L^2$ norm, and sufficiently
localized  
 yields a solution that relaxes to
another front, selected by a conservation law, in the $L^1$ norm at an algebraic rate
that we explicitly estimate.  We also obtain rates for the
relaxation in the
$L^2$ norm  and the rate of decrease of the excess free energy.

\end{abstract}

\bigskip
\centerline{Key words:  phase kinetics, fronts, non-linear stability, non local equation}
\centerline{Mathematics Subject Classification Numbers: 35A15, 35B40, 35K25, 35K45, 35K55, 35K65 }

\newpage

\tableofcontents

\section{Introduction and    main results    }
 
We consider   
   the nonlocal and nonlinear  evolution equation 
\begin{eqnarray}\label{1.1}
   \frac {\partial}{\partial t}  m(x,t)    &=&  \nabla\cdot\bigl(
\nabla  m(x,t)  - 
  {\beta (1-m(x,t)^2) (J\star
\nabla m )(x,t)} \bigr)  \quad x \in \R \times \L, \nonumber\\
\phantom{1} &{\phantom=}& {\phantom1}\nonumber\\
m(0,x)&=&  m_0(x), \quad m_0(x) \in [-1,1]  
\end{eqnarray}
in the $D-$dimensional cylinder  $ \R \times \L$ where   $ \L$ is a $(D-1)$-dimensional  torus of side-length   $L >0$, (equipped with the periodic Euclidean metric),  $ \beta>1 $,    $\star$ denotes convolution, and    $J$ is smooth,
spherically symmetric probability density on $\R^D$
with compact support.  In the following we set  $d= D-1$.

This equation first appeared
in the literature in a study \cite{LOP} of the dynamics of Ising systems
with a long--range interaction and so--called ``Kawasaki'' or ``exchange''
dynamics. In this physical context, $m(x,t)$ is the
magnetization density at $x$ at time $t$, viewed on the length scale of the
interaction J, and $\beta$ is the inverse temperature.
The derivation of   (\ref{1.1})  from the
underlying stochastic dynamics with $x$ taking
values in a torus $T^d$   is done in \cite{GL1}.    Equation  (\ref{1.1}) has been object of several studies that shall be quoted later.  
\vskip0.5cm

Our investigation in this paper turns on the fact that  
 the equation  (\ref{1.1}) can be written 
in a gradient flow form: Introduce the {\bf Gates-Penrose-Lebowitz free energy functional}
$\FF$ defined on all measurable functions from $\R\times \L$ by 
\begin{equation}\label{freeen}
{\cal F } (m) = \int_{ \R \times \L}  [f(m(x)) - f(m_\beta)] {\rm d} x+ {1\over 4} 
\int_{ \R \times \L}\int_{ \R \times \L}  J(x-y) [m(x)-m(y)]^2 {\rm d} x {\rm d} y 
\end{equation}  
where  $ f(m)$ is
\begin{equation*}f(m) = -{1\over 2} m^2 + {1\over \beta}\biggl[
\biggl({1+m\over 2}\biggr)\ln\biggl({1+m\over 2}\biggr) + 
\biggl({1-m\over 2}\biggr)\ln\biggl({1-m\over 2}\biggr)\biggr].
\end{equation*}
For $\beta >1$, this  potential function $f$  is a symmetric double well potential on
$[-1,1]$.  We denote the positive minimizer of $f$ on $[-1,1]$ by
$m_\beta$.
It is easy to see that $m_\beta$ is the positive solution of the equation
\begin{equation*}
m_\beta = \tanh (\beta m_\beta). 
\end{equation*}
The functional  (\ref{freeen}) is well defined on the set of measurable functions from $\R\times \L$ to $[-1,1]$, although it might be  infinity. 
The equation   (\ref{1.1}) can be written in the gradient flow form
\begin{equation}\label{gradflow}
{\partial\over \partial t}m    =  \nabla\cdot\left[ 
\sigma(m)\nabla\biggl({\delta {\cal F}\over \delta  m } \biggr)\right]
\end{equation}
where the {\it mobility} $\sigma(m)$ is given by
\begin{equation}\label{mobility}
\sigma(m) = \beta (1 -m^2). 
\end{equation}
{}From  this it follows, at least on a formal level, that $\FF$ is decreasing along the flow described by (\ref{1.1}):
The formal Frechet derivative
of the free energy  ${\displaystyle {\delta {\cal F}\over \delta  m }}$  is 
\begin{equation}\label{formd}
{\delta {\cal F}\over \delta  m } =  \frac {1}{\b}   {\rm arctanh}(m) -  J\star m, 
\end{equation} 
and thus, one formally  derives
\begin{equation}\label{dissipate}
{{\rm d}\over {\rm d} t}{\cal F}(m(t))   =  -\int_ {\R \times \L} \biggl|
  \nabla \left ({\delta {\cal F}\over \delta  m } \right) \biggr|^2
\sigma(m(t)){\rm d}x =: - {\cal I} (m(t)) \ .
\end{equation}

Based on this calculation, one might hope that $\FF$ would be a Lyapunov function governing the approach of solutions of (\ref{1.1})
to a minimizer of $\FF$.

The global minimizers of $\FF$ are of course the two constant profiles  $m(x) = m_\beta$ and  $m(x) = -m_\beta$ for all $x$ in the cylinder  $\R\times \L$.
Here we study a more interesting class of profiles $m$ under the constraint that $m(x)$ is very close to $-m_\beta$ far to the left in the cylinder, and is
very close to $m_\beta$ far to the right in the cylinder.  

More precisely, let us write $x = (x_1,x_1^\perp)$ where the first coordinate $x_1$ runs along the length of the cylinder, and $x^\perp$ along the cross section $\L$.
Consider the class ${\mathcal C}$ of measurable  functions $m$ from $\R\times\L$ to $[-1,1]$  such that for almost every $x_1^\perp\in \Lambda$,
$$\lim_{x_1\to -\infty}m(x_1,x_1^\perp) = -m_\beta\qquad{\rm and}\qquad  \lim_{x_1\to\infty}m(x_1,x_1^\perp) = m_\beta\ .$$

The minimizers of  $\FF$ over ${\mathcal C}$ can be expressed in terms of the minimizers of a simpler functional of one dimensional profiles.
More specifically, 
  in \cite{DOPT4}  is  was shown that there exists an  unique function $
\overline m_0(\cdot)$,  such
that
\begin{equation}\label{emnaught}
{{\cal F}_1}(\overline m_0) = \inf \left \{\ {{\cal F}_1}(m)\ \biggl|\ {\rm sgn}(x_1)m(x_1)\ge
0, \lim_{x_1 \to \pm \infty} {\rm sgn}(x_1)m(x_1) >0  
 \right \},\
 \end{equation}
where   $\FF_1$ is  is  the functional    
\begin{equation*}
\FF_1(m) = \int_{ \R}  [f(m(x_1)) - f(m_\beta)] {\rm d} x_1+ {1\over 4} 
\int_{ \R }\int_{ \R }  \overline J(x_1-y_1) [m(x_1)-m(y_1)]^2 {\rm d} x_1 {\rm d} y_1, 
\end{equation*}
and 
\begin{equation}\label{tildeJ}
\overline J(x_1) = \int_{\L}J(x_1,x^\perp) 
  {\rm d}x^{\perp}.
  \end{equation}
Furthermore it is shown that $\overline m_0$ is an odd, $C^\infty (\R)$,
increasing function,  and
\begin{eqnarray}\label{decay}
0 &< &m_\beta^2 - \overline m_0^2(x_1)  \le Ce^{-\alpha |x_1|}\ , \nonumber\\
0 &<& \quad {\overline m}'_0 (x_1)\ \quad\ \le Ce^{-\alpha |x_1|}\ ,\nonumber\\ 
0 &<& \quad |{\overline m}''_0 (x_1)|  \quad\     \le Ce^{-\alpha |x_1|}\ , \nonumber
\end{eqnarray}
for positive constants $C$ and $\alpha $ depending on $J$ and
$\beta$.  
The first   two of these  estimates are proved in  \cite{DOPT4} and the third one
in \cite{DGP}.  A review of  these  and     related  results  can be found in Chapter 8 of the book    \cite{Pr}. 
The subscript $0$ on the minimizer refers to the fact that the constraint
imposed in (\ref{emnaught}) breaks the translational invariance of the free
energy. For any $a$ in $\R$,   
 define
\begin{equation}\label{ema1}
\overline m_a (x)= \overline m_a (x_1, x_1^\perp)=\overline m_a(x_1) = \overline m_0(x_1-a), \quad  x \in \R \times \L.
\end{equation}
 Clearly 
\begin{equation}\label{euler}
{\cal F}(\overline m_a) = {\cal F}(\overline m_0),  \qquad     \frac{\delta \FF} {\delta m}(\overline m_a) =    \frac {1}{\b}  {\rm arctanh}(\overline m_a) -    J\star  \overline m_a =0.   
\end{equation}  
Thus the profiles $\overline m_a$ are at least critical points of the free energy $\F$  in the class $\mathcal {C}$.
Since they are built out of minimizing one dimensional profiles, it is natural to guess that they are in fact  minimizers in $\mathcal{C}$. 
This has been proved by  Alberti and Bellettini  \cite {AB}, who showed moreover  that every minimizer of $\F$ in $\mathcal {C}$ is of this form.  
The functions in this one parameter family of minimizers of the free energy
$\overline m_a$, $a\in \R$,  are the   stationary solutions of  (\ref{1.1})  whose stability is to be investigated here.

Because the free energy is 
reflection invariant, there is also another
family, obtained by reflecting the previous one. However, these two families of minimizers are well separated
in all of the metrics in which we shall work, and it suffices to consider only one
of them.

We shall be concerned here with the evolution of small perturbations of $m $ from   $\overline m_0 $, and their relaxation to $\overline m_a $ for some $a$
under the dynamics introduced above. We shall show that if the perturbation is suitably small, then this happens, and moreover, we shall find the value of $a$, and estimate the rate of convergence.

The equation (\ref{1.1})  not only has a Lyapunov function; it has a conservation
law as well: For times $t$ in any interval on which $m(x,t) - {\rm sgn } (x_1) m_\b $ is
integrable,
\begin{equation*}
\frac{{\rm d}}{{\rm d}t}\int_{\R \times \L}\bigl(m(x,t) - \overline m_b(x)\bigr){\rm d}x = 0
\end{equation*} 
for any $b$. Therefore, if one defines $a$  in terms of the
initial data $m_0$ for
(\ref{1.1}) by
\begin{equation}\label{conser}
\int_ {\R \times \L}\bigl(m(x,0) - \overline m_a(x)\bigr){\rm d}x = 0\ ,
\end{equation}
 one has for the
solution
\begin{equation*}
\int_ {\R \times \L} \bigl(m(x,t) - \overline m_a(x)\bigr){\rm d}x = 0
\end{equation*}
 for all $t$
or at least all $t$ such that $m(s, x)- {\rm sgn } (x_1) m_\b $ is integrable for all
$s\le t$.

Now, formally invoking the Lyapunov function and the conservation law,
it is easy to guess the result of solving (\ref{1.1}) for initial data $m_0$
that is a small perturbation of the front $\overline m_0$: The decrease of the
excess free energy ``should'' force the solution $m(t)$ to tend to the family
of fronts, and the conservation law ``should'' select $\overline m_a$ as the front it should
be converging to, so the result ``should'' be that
$$\lim_{t\to\infty} \int_ {\R \times \L} |m(x,t) - \overline m_a(x)| {\rm d}x  =0$$
with $a$ given in terms of the initial data $m_0$ through (\ref{conser}).

There is a fundamental obstacle in the way of this optimistic line of reasoning:  Consider a very small $\e>0$ and
consider
$$
m(x) := \overline m_0(x) + \epsilon 1_{[\e^{-3/2},2e^{-3/2}]}(x_1)\ .
$$
Then it is very easy to see that
$$\F(m) - \F(\overline m_0) = {\mathcal O}(\epsilon^{1/2}) \quad{\rm while}\quad \inf_{a\in \R} \|m - \overline m_a\|_1 =    {\mathcal O}(\epsilon^{-1/2})\ .$$
That is, perturbations of a minimizer with extremely small excess free energy can be very far from any minimizer in the $L^1$ norm:  
If the only information on the evolution that one had was that the excess free eenrgy was decreasing to zero, one could not rule out the posibility
that the $L^1$ distance to the nearest minimizer might be increasing to infinity. As we shall see, all profiles $m$ with small excess free energy
and a large $L^1$  distance to the nearest minimizer are very delocalized perturbations of minimizers, spread out on a very large scale, as in the example we have given. To rule this out, we have to assume moment conditions that prevent our perturbations from being too delocalized at the beginning, and then we must work to show that this localization does not deteriorate too rapidly.  This accounts for the moment conditions in the
theorem stated below. These moment conditions are essential for bounding the rate of convergence; a small initial excess free energy is not enough.

In the following, whenever there is no ambiguity, we  denote by
$  \|u\|_p$      the $L^p (\R \times \L)$  of a function $u$. If  $u \in W^{s,2} (\R \times \L)$, $s \in \N$, the space of functions $u \in  L^2 (\R \times \L) $ whose  distributional  derivatives of order $\le s$ are in  $ L^2 (\R \times \L)$,   we   denote by $ \| u \|_{W^{s,2}} = \sum_{ |\alpha | \le s }  \|D^\alpha u\|_2$   its norm.   
 
 We have the following main result.
\begin{thm}\label{main1}  Let $m(t)$ be the solution of equation  (\ref{1.1}) in the $D$-dimensional cylinder $\R \times \L$, $D\leq 3$, and with   initial data $m_0$  such that 
$$\int_{\R \times \L} x_1^2  (m_0(x) - \overline m_0(x))^2   {\rm d}x \le c_0\ ,$$  where  $ c_0 $ is any
positive constant. Then for any $\d > 0$ there is a strictly positive constant
$\e =\e(\d,c_0,\b,J, L)$   depending only on $\d$, $c_0$, $\b$,  $J$  and   $L$   such that  for
all initial data $m_0$ with $-1 \le m_0 \le 1$,  and with
$$  \int_{\R \times \L}(m_0(x) - \overline m_0(x))^2 {\rm d}x
 \le \e \ ,$$ the excess free energy $\FF(m(t)) - \FF(m_0)$ of the
corresponding solution
$m(t)$ of  (\ref{1.1}) satisfies
\begin{equation}\label{result111}
{\cal F}(m(t)) - {\cal F}(\overline m) \le c_2 (1 + c_1 t)^{-(9/13 -\d)} 
\end{equation}
and
\begin{equation}\label{result222}
\|m(t) - \overline m_a\|_1 \le c_2 (1 + c_1 t)^{-(5/52 -\d)}  
\end{equation}
 where $c_1$ and $c_2$ are
finite constants depending only on $\d$, $c_0$,
$J$, $\beta$ and $L$ 
 and $a$ is given by (\ref{conser}).
 \end{thm}

 \vskip0.5cm
   In  $D=1$ the same stability   problem  for the equation  (\ref{1.1})    was addressed
in the papers   \cite{CCO1} and   \cite{CCO2}.   The strategy used in these papers was    applied in      \cite{CCO3}  to  show  the  local non-linear stability of the interface solution  for the Cahn-Hilliard equation, always in    $D=1$.

    The  method  applied     in $D=1$  has been adapted   in this paper to   show  local non-linear stability of the interface solution  of  (\ref{1.1})  when  dimension $D \ge 2$.     To apply the previous strategy in    $D \ge 2$ one needs  to   control the transverse contribution  of the perturbation to the planar fronts.

This is done   by  a suitable splitting of a  function  in $\R \times  \L$  as the sum of  two  functions, one   depending   only   on  $x_1\in \R$ and  the other  with mean zero in the direction  orthogonal to $x_1$.  This  allows us to effectively decouple the problem into transverse and longitudinal parts,
and to  control  the gradient of the function in the  transverse direction  applying  the  Poincar\'e inequality.
 
The method  is  robust enough  and it should allow to deal with   nonlinear  local stability problems  for   other equations of  Cahn-Hilliard type.   

There are very few results in the literature regarding  stability of the planar   fronts  in  infinite domain   for    equations of  Cahn-Hilliard type.   
 The only paper to our knowledge dealing with the stability of the planar front  for    Cahn-Hilliard equation in   $ \R^D$,  for $D \ge 3$ is the paper by   Korvola,  Kupianen and Taskinen \cite{KKT}.    They proved that  the leading asymptotic of the solution is characterized by a length  scale proportional  to $t^\frac {D-1}3$ instead of the usual  $t^\frac D2$  typical to parabolic problems.   In contrast to the one dimensional  and to  $D-$ dimensional cylinder  setting,  considered in this paper,   they show 
that    the translation of the front tends to zero  as time tends to infinity.  This is because a localized perturbation is not able to produce a constant shift in the whole transverse space $\R^{D-1}$.           In our case, a perturbation of  an equilibrium front  need not return  asymptotically 
 to the initial front.  Indeed, there is no easy argument using only decrease of free energy to show that the perturbation does not cause the 
front to ``run away to infinity''. Our method provides a proof, with quantitate estimates, on the size of the shift that can result as the perturbation is dissipated way.   

  The restriction to $D\le 3$ is for reasons that are surely technical; the condition $D \le 3$ is used only in proving certain regularity estimates 
 that are required in our central arguments. Very likely with more labor (and more pages), these could be proved in higher dimension as well.  
    
To implement the  heuristics discussed before the theorem is not so simple as one might hope.  There are
several reasons for this.  The first has to do with the relevant norms.

To explain the physical relevance of the $L^2$ norm, we note that $a \mapsto  \|m - \overline m_a\|_2^2$ is strictly convex whenever
for some $a\in \R$, 
\begin{equation*}
\|(m- \overline m_a)' \|_2  < \|\overline m_0' \|_2\ .
\end{equation*}
Using this we show in  Theorem~\ref{5} that under suitable smallness assumptions on $m - \overline m_a$ for some $a$, there exists a 
unique $b\in \R$ so that 
$$\|m - \overline m_b\|_2 = \inf_{a\in \R}\{  \|m - \overline m_a \|_2 \}\ .$$
Further we   show, see Lemma~\ref{A2},   that the excess
free energy measures the distance to this closest front in the $L^2$ metric in the sense that   
\begin{equation}\label{zaazoo}
\FF (m ) - \FF (\overline m_b)  \simeq  C\|m-\overline m_b\|_2^2 
\end{equation}
 under suitable smallness assumptions on $m - \overline m_b$. 

We use the smoothing
properties of (\ref{1.1}) to obtain the condition on the derivative of $m$  for all
$t \ge t_0$, for some finite $t_0$  so that we can apply Lemma~\ref{A2}. 

On account of (\ref{zaazoo}), for any solution $m(t)$ of  (\ref{1.1}), 
define $a(t)$ to be that value of $a$ such that
\begin{equation}\label{madef}
\|m(t) - \overline m_{a(t)}\|_2 = \inf_{a\in \R}\{ \|m(t) - \overline m_a\|_2 \}
\end{equation}
and note that $a(t)$ is a well--defined function
as long as $\|m(t) - \overline m_{a(t)}\|_2$ stays sufficiently small 
since then the minimum is uniquely attained. (We shall do all of our analysis
in this paper for times $t$ in an interval $(t_0,T_0)$ on which
$\|m(t) - \overline m_{a(t)}\|_2$ does stay small, and then at the end we shall show that
$T_0 = \infty$.)
Hence, if one proves that the excess free energy decreases to zero, the
best one can obtain from this is that
$$\lim_{t\to\infty}\|m(t)-\overline m_{a(t)}\|_2 = 0.$$
However, this does not yield any information on $a(t)$ -- 
and it cannot by the translation 
invariance  of the free energy. 
The conservation law would give us information on $a(t)$, but to use it
we require $L^1$ control on $m(\cdot,t) - \overline m_{a(t)}(\cdot)$. Since
$$\|m(\cdot,t) - \overline m_{a(t)}(\cdot)\|_\infty \le 2$$
{\it a--priori}, $L^1$ control would give us $L^2$ control through
$$\|m(\cdot,t) - \overline m_{a(t)}(\cdot)\|_2^2 \le 2
\|m(\cdot,t) - \overline m_{a(t)}(\cdot)\|_1$$
but not {\it vice--versa}. 
In order to use the conservation law to show that 
$
 \lim_{t \to \infty} a(t) =a 
$
where $a$ is given by (\ref{conser}), we must, and shall, show that
$$\lim_{t\to\infty}\|m(\cdot,t)-\overline m_{a(t)}\|_1 = 0\ .$$

Before discussing the $L^1$ behavior of
perturbations of fronts, we make the following convention, 
to be used throughout the paper, whenever
some solution $m(x,t)$ is under discussion:
\begin{equation}\label{perturbdef}
v(x,t) = m(x,t) -\overline m_{a(t)} (x)
\end{equation}
where $a(t)$ is given in (\ref{madef}), and moreover 
\begin{equation}\label{mbarconv}
\overline m(x) \qquad{\rm denotes}\qquad \overline m_{a(t)}(x) .
\end{equation}

As explained above, we shall have to look into the details of the free energy dissipation ${\mathcal I}$, see \eqref {dissipate},  in order to understand whatever stability properties our equation may have. To begin this, we write
\begin{eqnarray}\label{8.1}
  {\cal I} (m)  &=& 
 \int _ {\R \times \L} \s (m(x)) \left 
 [\frac {\partial}{ \partial x_1 }\left(
{1\over \b}{\rm arctanh}m(x)-(J \star m)(x)
\right ) \right ]^2  { \rm  d } x \nonumber\\
   &+&
\int_ {\R \times \L} \s (m(x)) 
 \left [  \nabla^{\perp} \left(
{1\over \b}{\rm arctanh}m(x)-(J \star m)(x)
\right ) \right ]^2  { \rm  d } x,   
\end{eqnarray} 
where 
 $  \nabla^{\perp} $ is the gradient in the orthogonal direction of $x_1$.

One   result  of the paper, Theorem~\ref{51}, gives a  lower bound on the rate 
of dissipation of the excess free energy, whenever the dimension $D \le 3 $.
 For any $\e>0$,
\begin{eqnarray}\label{doodah}
  {d\over dt}\bigl[{\cal F}\big(m(t))-{\cal F}\big(\overline
m     \big)\bigr]  &=& -  
  \II \big (m( t)\big ) \nonumber\\            
    & \le& -   (1-3 \e)   \sum_{i\ge1} \int_{\R \times \L}
  \sigma(\overline m)\bigl[(\BB v (t))_{x_i}  \bigr]^2{\rm d}x, 
 \end{eqnarray} 
where $  \II \big (m( t)\big )$ is given in  (\ref{8.1})  and, recall \eqref {mbarconv},  $\bar m := \bar m_{a(t)}$.
To get this  result we need smoothness estimates to hold for the solution $m(t)$, in order to apply Sobolev inequalities.
 Namely (\ref{doodah}) holds  when the  $ ||v||_{W^{s,2}} \le  \kappa_1(\b,J,L,\e)$ and 
$||v||_2 \le  \d_1(\b, J,L, \e)$ for some strictly positive 
constants $\kappa_1(\b,J,L,\e)$ and $\d_1(\b, J,L, \e)$, where    $ ||v||_{W^{s,2}}$ is the Sobolev norm, see in the Appendix,  Lemma~\ref{A1} and $s> \frac {D} 2$. 
We have  quantitative estimates, see Theorem~\ref{smoo},  of    the derivative of all order  of $m(t)$ only when $D \le 3$.   
This is the rather technical reason, noted above,  for which  we impose   the constraint on the dimension. 
The  ${\cal B}$ in  (\ref{doodah})
denotes the second variation of the free energy $\FF$ at $\overline m$.
By our convention, $\overline m$ denotes $\overline m_{a(t)}$, and while it is
occasionally preferable to write ${\cal B}_{a(t)}$ to make this explicit,
we shall generally simply write ${\cal B}$, and leave the dependence on
$a(t)$ implicit. However, in recalling the definition, we shall be 
explicit:
\begin{equation}\label{secondvar}
\langle u,{\cal B }_a u\rangle_{L^2} = {{\rm d}^2\over {\rm d}s^2}{\cal F} 
(\overline m_a +
su)\biggr|_{s=0}. 
\end{equation}
The  properties of ${\cal B}$ that we shall  use in our analysis
  are discussed in Section 3.  
Because of the derivatives, the quadratic form on the right
in (\ref{doodah}) has no 
spectral gap. If it did, this together with (\ref{zaazoo}) would provide an
exponential rate of decrease of the excess free energy, and hence of 
$\|v(t)\|_2$.  
Since there is no spectral gap here, one needs additional monotonicity, or
at least {\it a--priori} boundedness properties to exploit (\ref{doodah}),
as explained in \cite{CCO1} and \cite{CCO2}.
In the study of parabolic equations
\begin{equation}\label{parabolic}
{\partial u\over \partial t} = 
\nabla\cdot\bigl({\bf D}(u)\nabla u\bigr)\ ,
\end{equation}
where $ {\bf D}(\cdot) $ is the diffusivity  matrix for which there is also no spectral gap,
\begin{equation}\label{down}
\frac{{\rm d}}{{\rm d}t}\int |u(x,t)|{\rm d}x \le 0
\end{equation}
which trivially provides the additional monotonicity required to 
show that
$$\sup_{t \ge 0}\|u(t)\|_1  \le  \|u(0)\|_1\ .$$
Then a standard argument with the Nash 
inequality allows one to conclude that  $\|u(t)\|_2$  decreases to zero at an 
algebraic rate, at least when the diffusivity ${\bf D}(\cdot) $  in  (\ref{parabolic}) is
bounded from below.

This route is closed to us since the analog of (\ref{down}) does not
hold  for $v(t)$ when $m(t)$ is a solution to (\ref{1.1}),  see \cite{CCO1} and   \cite{CCO2} for further details.
Moreover, there are other
problematic non-dissipative features,  the
maximum principle fails to hold for (\ref{1.1}) and 
the 
free energy is {\it not Frechet differentiable} on the natural set of functions
that is invariant under  the evolution prescribed by (\ref{1.1}).
  Namely, recall   (\ref{formd}), for any $m$ with $-1\le  m \le 1$, $ J\star m$ is bounded, but 
${\rm arctanh}(m) = \pm\infty$ on $\{x \ |\ m(x) = \pm 1\}$.   This
means that some care must be taken with the use of the key dissipativity property
(\ref{dissipate}) whose formal derivation depends on this Frechet differentiability.
Even worse, however, is that the mobility (\ref{mobility}) vanishes where
$m = \pm 1$, and with it the local contribution to the dissipation in (\ref{dissipate}).

One  way to   obtain bounds on the decay  in the $L^1$  norm 
is to apply  a  strong formulation of  the 
``uncertainty principle'';  as done in  \cite{CCO2}.  

 We  illustrate this in the case  of  the heat equation in   the  $D-$ dimensional cylinder $ \R \times \L$.  
Recall that we denote  $x=(x_1, x_1^{\perp})$, $x_1 \in \R$.  
Consider a solution $u$ of the heat equation
\begin {equation*} {\partial \over \partial t}u(x,t)  =  \Delta u(x,t), \quad x \in \R \times \L  \end{equation*}
with integrable initial data $u_0$, 
and suppose that
\begin{equation*}
 \int_{\R \times \L}  u_0(x){\rm d}x = 0.\ 
 \end{equation*}
 We show in Theorem 2.1 of  \cite{CCO2}
that under the constraint
$$\int \psi(x_1){\rm d}x_1 = 0$$
one has
\begin{equation}\label{uncertainty2}
\biggl(\int x_1^2|\psi(x_1)|^2{\rm d}x_1\biggr)
\biggl(\int|\psi'(x_1)|^2{\rm d}x_1 \biggr) \ge  {9\over 4} \biggl(\int|\psi(x_1)|^2{\rm
d}x_1\biggr)^2. 
\end{equation}
 Define 
\begin{equation}\label{ex.0}
f(t)= \int_{\R \times \L} |u(x,t)|^2{\rm d}x\qquad{\rm and}
\qquad \phi(t) = \int_{\R \times \L}  x_1^2|u(x,t)|^2{\rm d}x +1.
\end{equation}
One then computes that
\begin{equation}\label{ex.10}
 \frac{{\rm d}}{{\rm d}t}f(t) = -2\int|\nabla u(x,t)|^2{\rm d}x. 
 \end{equation}
\begin{eqnarray}\label{heatex}
 \frac{{\rm d}}{{\rm d}t}\phi(t)  &=&  -4 \int x_1      e_1 \cdot u(x,t)(\nabla u)(x,t) {\rm d}x  - 2 \int x_1^2      |(\nabla u)(x,t)|^2 {\rm d}x \nonumber\\
  &\leq&  2f(t).  
  \end{eqnarray}
  We would like, as it will be clear in the following, to write equations   (\ref{ex.10}) 
  and   (\ref{heatex})
   in a  closed   form, i.e.  to write the right hand side of  (\ref{ex.10}) in term of $f(t)$ and $\phi(t)$. 
  Denote 
 $$ J(u(t))= \int|\nabla u(x,t)|^2{\rm d}x. $$
 Since the dependence on $t$ does not play any role in the following calculations,  we do not write it  explicitly. Set 
 $$u(x) = v_1(x_1) + w(x)  $$ where 
 $$ v_1(x_1):= \frac 1 {L^d} \int_{\L} u (x_1, x^{\perp}) d  x^{\perp}, \quad   x_1 \in \R. $$ 
By construction  
\begin{equation}\label{ex1} 
\int_{\L} w (x_1, x^{\perp}) d  x^{\perp}=0,   \quad 
\forall x_1 \in \R.  
\end{equation}
 We have 
\begin{eqnarray*}
 J(u) &=&  \int|\nabla u(x)|^2{\rm d}x = \int|u_{x_1}(x)|^2{\rm d}x  +  \int|\nabla^{\perp} u(x)|^2{\rm d}x  \nonumber\\
 &=&  \int |v'_1(x_1)  + w_{x_1} (x)|^2{\rm d}x  +  \int|\nabla^{\perp} w(x)|^2{\rm d}x, 
   \end{eqnarray*}
    where  $v'_1 $ is  the  spatial derivative of $v_1$.
 Notice that, because of  (\ref{ex1}),
 $$ \int_{\R \times \L} v'_1(x_1)  w_{x_1} (x){\rm d}x =0.$$
 So we  get 
\begin{equation}\label {ex.2a}
  J(u) =  \int_{\R \times \L}  |v'_1(x_1)|^2 {\rm d}x    +  \int_{\R \times \L}   |w_{x_1} (x)|^2{\rm d}x    +  \int|\nabla^{\perp} w(x)|^2{\rm d}x.    
  \end{equation}
Again by    (\ref{ex1}) 
we can   bound from below the last term of    (\ref{ex.2a})   using   the
 Poincar\'e inequality in $\L$.  It states    that 
 for $g \in L^2(\L)$ and $\nabla g \in L^2 (\L)$ 
$$c(d)\|g- \overline g\|^2_2 \le L^2 \|\nabla g\|_2^2 $$ 
where $\overline g= \frac 1 {|\L|} \int_\L g(x) dx $ and $c(d)$  depends on dimensions.

 Applying this for each $x_1$ in   (\ref{ex.2a}), and using the fact that because of   (\ref{ex1}),   $ \overline w(x_1)=0$,  we obtain
\begin{eqnarray}\label{M.u1} 
     \|  \nabla^{\perp}  w\|^2_2  &=&   \int _\R   { \rm  d } x_1   \left ( \sum_{i\ge2}   \int_{ \L}  { \rm  d } x^\perp |  ( w(x_1, x^\perp))_{x_i} |^2 \right )\nonumber\\
      &\geq& 
\frac {c(d)}   {L^{2}} \int_\R  \int_{ \L} |   w(x_1, x^\perp) |^2  dx.   
\end{eqnarray}
To  lower bound the first  term of    (\ref{ex.2a}) we apply the uncertainty principle  (\ref{uncertainty2}) as following: 
\begin{eqnarray}\label{M.u} 
  \int_{\R \times \L}  |v'_1(x_1)|^2 {\rm d}x   &\ge&  \frac 9 4  L^d \frac {  \left (  \int_{\R}  |v_1(x_1)|^2  \right )^2 }   { \biggl(\int x_1^2|v_1(x_1)|^2{\rm d}x_1\biggr)} \nonumber\\
&=& \frac 9 4   \frac {  \left (  \int_{\R \times \L}  |v_1(x_1)|^2  \right )^2 }   { \biggl( \int_{\R \times \L} x_1^2|v_1(x_1)|^2{\rm d}x_1\biggr)}.     
   \end{eqnarray}
Therefore, taking into account     (\ref{M.u1})  and (\ref{M.u})  we  lower bound  \eqref {ex.2a} as following:  
\begin{eqnarray} \label {ex.5}
 J(u) &\ge&  
  \frac 9 4   \frac {  \left (  \int_{\R \times \L}  |v_1(x_1)|^2  \right )^2 }   { \biggl( \int_{\R \times \L} x_1^2|v_1(x_1)|^2{\rm d}x_1\biggr)}    +  \frac {c(d)}  {L^{2 }}  \int_{\R \times \L}  |w(x)|^2{\rm d}x   \nonumber\\
   &\ge&
    \frac 9 4   \frac 1    { \biggl( \int_{\R \times \L} x_1^2|v_1(x_1)|^2{\rm d}x_1\biggr) + 1}   \left \{    \left (  \int_{\R \times \L}  |v_1(x_1)|^2  \right )^2     + \frac 49    \frac {c(d)}  {L^{2 }}  \int_{\R \times \L}  |w(x)|^2{\rm d}x   \right \}.  
   \end{eqnarray}
 Notice that we  obtained  this estimate by dropping the contribution of the second term in \eqref {ex.2a}, i.e.   $\int_{\R \times \L}   |w_{x_1} (x)|^2{\rm d}x$.  When dealing with   equation \eqref {1.1} we     shall need to keep this term to control the  non linearity. 
  
      By orthogonality
$$ \|u\|^2_{L^2(\R \times \L)}= \|v_1\|^2_{L^2(\R \times \L)}+ \|w\|^2_{L^2(\R \times \L)}. $$
Therefore 
\begin{eqnarray*} 
 \|v\|^4_{L^2(\R \times \L)}&=& \left (  \|v_1\|^2_{L^2(\R \times \L)}+ \|w\|^2_{L^2(\R \times \L)} \right )^2 \nonumber\\
  &=&  \|v_1\|^4_{L^2(\R \times \L)} +  \left ( \|w\|^2_{L^2(\R \times \L)} +2  \|v_1\|^2_{L^2(\R \times \L)}\right )  \|w\|^2_{L^2(\R \times \L)}.
   \end{eqnarray*}  
Suppose that 
$$  \left ( \|w\|^2_{L^2(\R \times \L)} +2  \|v_1\|^2_{L^2(\R \times \L)}\right )  \le \frac 49    \frac {c(d)}  {L^2 }. $$ 
Therefore
\begin{equation}\label{ex.6}
\|v\|^4_{L^2(\R \times \L)} \le  \|v_1\|^4_{L^2(\R \times \L)} + \frac 49    \frac {c(d)}  {L^2 }  \|w\|^2_{L^2(\R \times \L)} .  
\end{equation}
Further
\begin{eqnarray}\label{ex.8} 
   \int_{\R \times \L} x_1^2|v_1(x_1)|^2{\rm d}x   &\le&
 \int_{\R \times \L} x_1^2  \left [ |v_1(x_1)|^2  +  |w(x)|^2 \right ] {\rm d}x   \nonumber\\
& =&  \int_{\R \times \L} x_1^2  \left [ |v_1(x_1)+ w(x)|^2\right ] {\rm d}x 
=  \int_{\R \times \L} x_1^2   |u(x)|^2 {\rm d}x.  
\end{eqnarray}
Taking into account  (\ref{ex.0}),  \eqref{ex.5},    (\ref{ex.6}),  (\ref{ex.8}) we have
$$ J(u(t)) \ge    \frac 9 4  \frac { f^2(t)} { \phi (t)}.  $$ 
Therefore, from  (\ref{ex.10}) we get 
 $$ \frac{{\rm d}}{{\rm d}t}f(t) = -   \frac 9 2   \frac { f^2(t)} { \phi (t)}.  $$ 
Recalling  (\ref{heatex}) we get  the system of differential inequalities,
\begin{eqnarray}\label{system} 
 \frac{{\rm d}}{{\rm d}t}f(t) 
&\leq& - A{f(t)^2\over \phi(t)}\nonumber\\
 \frac{{\rm d}}{{\rm d}t}\phi(t) &\leq& Bf(t)
\end{eqnarray}
with $A = 9/2$ and $B = 2$. Theorem 5.1 of \cite{CCO1} says that for 
any  non negative solution of  
(\ref{system}),  the following holds

\begin{eqnarray}\label{may1} 
f(t) &\leq& f(0)^{1-q}\phi(0)^q 
\left({\phi(0)\over f(0)} +  (A+B) t  \right)^{-q} \nonumber\\
\phi(t) &\leq& f(0)^{1-q}\phi(0)^q 
\left({\phi(0)\over f(0)} +  (A+B) t  \right)^{1-q} 
\end{eqnarray}
where
 $$q = {A\over A+B}\ .$$  
In the case at hand, this is
$$q = {9\over 13}\ .$$
Since this value exceeds $1/2$, we get $L^1$ decay in the following way:
We prove in     Section 6, Lemma~\ref{5.2a}    that
for any function $f$, so that  $\|(1 + x_1^2)^{1/2}f\|_{L^2 (\R \times \L)}$ is finite   and for   any $0 <\d < 1$  we have
$$ \|f\|_1 \le C(\d, L)\|(1 + x_1^2)^{1/2}f\|_2^{(1+\d)/2}\|f\|_2^{(1-\d)/2}$$
where
$C(\d,L)$ is  given explicitly in the lemma.
Since $9/13 > 1/2$ for $\d$ sufficiently small, we have that
$\|(1 + x_1^2)^{1/2}u(t)\|_2^{(1+\d)/2}$, see  (\ref{may1}),  increases more slowly
than $\|u(t)\|_2^{(1-\d)/2}$ increases, and so $\|u(t)\|_1$
decreases to zero. In fact, the rate one gets for $\|u(t)\|_1$ is arbitrarily close to
$t^{-5/13}$, for $\d$ sufficiently small.    Actually, one can do better for the heat equation.  One can obtain,  as      in \cite{CCO1},   $ \frac{{\rm d}}{{\rm d}t}\phi(t) \le  \frac 32 f(t)$ (we have 2 in   (\ref{heatex})).  Then  $B= \frac 32$  and   the  rate one gets for $\|u(t)\|_1$ is arbitrarily close to
$t^{-\frac 1 4}$, for $\d$ sufficiently small. 

The previous argument presented for the heat equation can be implemented   for equation (\ref{1.1}). 
We define 
\begin{equation*}
f(t) = \FF(\overline m+v(t)) - \FF(\overline m) \qquad{\rm and}
\qquad \phi(t) = L^d+ \int \sigma(\overline m)x_1^2|{\cal B}v(x,t)|^2{\rm d}x\ ,
\end{equation*}
where $v$ as in \eqref {perturbdef} and $ \BB$ as in \eqref {secondvar}.
We could  estimated the time derivatives of these quantities   obtaining  bounds of the form given in (\ref{system}),
but with inexplicit constants $A$ and $B$.

Now, the rate of decay that one gets by this method depends very much on the
ratio of the constants $A$ and $B$ in (\ref{system}). To get $L^1$ decay, we need this ratio
to be fairly close to the ratio $9/13$ obtained for the heat equation.

We do this by exploiting the following alternatives:
for any $\e_1>0$,  at any time $t$,  one has either
\begin{equation}\label{mucho}
{\cal I}(m(t)) \le \e_1 \bigl[{\cal F}(m(t)) - {\cal F}(\overline m)\bigr]\ ,
\end{equation}
or 
\begin{equation}\label{meno}
{\cal I}(m(t)) \ge \e_1 \bigl[{\cal F}(m(t)) - {\cal F}(\overline m)\bigr]\ ,
\end{equation}
where ${\cal I}$ is the dissipation functional (\ref{dissipate}). 

We prove in Section 4 that
for any $\e > 0$, there are strictly positive constants 
$\d_0(\b,J,L,\e)$, $\kappa_0(\b,J,L,\e)$ and $\e_1(\b,J,L,\e)$
   depending only on $\b$, $J$, $L$  and $\e$,
such that for all $t$ for which   (\ref{mucho}) is satisfied together with 
$||v(t)||_{W^{s,2}} <\kappa_0$, $ s > \frac D 2$, $||v (t)||_2 <\d_0$ and $|a(t)| \le 1$,
it is the case that        
$$\frac{{\rm d}}{{\rm d}t}\bigl[{\cal F}\big(m(t))-{\cal F}\big(\overline m  \big)\bigr]
\le -9(1-\e)(1-\sigma(m_\b))^2    
{\bigl[{\cal F}\big(m( t))-{\cal F}\big(\overline m  \big)\bigr]^2 \over
\phi ( t )}\ . $$
We then show in Section 5 that under the same assumptions of Section 4, 
it is the case that        
$$\frac{{\rm d}}{{\rm d}t}\phi(t) \le 
(1+\e)4(1-\sigma(m_\b))^2
\bigl[{\cal F}(\overline m+v) - {\cal F}(\overline m)\bigr]\ .$$
Notice the condition that $|a(t)| \le 1$, to which we shall return. Thus,
when (\ref{mucho}) holds, we have
\begin{eqnarray}\label{system16} 
 \frac{{\rm d}}{{\rm d}t}f(t) 
&\leq& - \tilde A{f(t)^2\over \phi(t)}\nonumber\\
 \frac{{\rm d}}{{\rm d}t}\phi(t) &\leq& 
\tilde Bf(t)\end{eqnarray}
with the difference between $\tilde A/\tilde B$ and $9/13$ arbitrarily small
for $\e$ small enough for all times $t$ such that $\|v(t)\|_2$,
  $||v(t)||_{W^{s,2}} $,  $ s > \frac D 2$, are sufficiently small and $|a(t)|\le 1$.

On the other hand, when (\ref{mucho}) is violated and (\ref{meno}) holds, the dissipation
is large, and this works in our favor. In Section 6, we exploit this 
alternative to prove Theorem~\ref{main1}. The proof is still somewhat intricate,
and it would have been simplified had we been able to show the existence
of a time $t_\star$ such that (\ref{mucho})   holds for all $t  \ge t_\star$. If this
were the case, the constants
$\tilde A$ and $\tilde B$ above would govern the decay, and we would
obtain a bound on the excess fee energy of the form
$$\bigl[{\cal F}(m(t)) - {\cal F}(\overline m)\bigr] \le
C(1 + A(1-\sigma(m_\b))^2 t)^{-q}$$
where $A$ does not depend on $\b$. Since $(1-\sigma(m_\b))^2$
vanishes as the critical temperature is approached,  this would indicate how the rate
of relaxation slows in this limit. In any case, our results do
show that it is possible to estimate the {\it exponent} in the rate of
relaxation independently of $\b$.

To explain why (\ref{mucho}) enables us to obtain what are essentially
heat equation constant in (\ref{system16}), one has to view it as a smoothness
condition. Indeed,  it follows from Theorem~\ref{51}  that
$$ (1-\e) \sigma(m_\b)\|   |\nabla ({\cal B}v )| \|_2^2 \le  {\cal I}(\overline m + v)$$
for any $\e$, under appropriate 
 conditions on $v$.
Hence, by Lemma~\ref{A2}  which compares $\|{\cal B}v\|_2^2$
and the excess free energy of $\overline m + v$, when (\ref{mucho}) holds,
\begin{equation*}
\|  \nabla \left ({\cal B}v \right) \|_2^2 << \|{\cal B}v\|_2^2\ .
\end{equation*}
Next, the action of ${\cal B}$ on functions $w$ that satisfy
$$\|\nabla w \|_2 << \|w\|_2$$
is particularly simple: As shown in  Lemma~\ref{TT1} in the appendix,
$${\cal B }w \approx \tilde\alpha w$$
where $\tilde\alpha = 1/\sigma(m_\b) -1$.  Once one may replace ${\cal B}$ with
multiplication by $\alpha$, the linearized version of (\ref{1.1}) does
become essentially the heat equation. This discussion is heuristic, but
in no way misleading, and hopefully motivates the technical preliminaries in 
Section 2.

Before turning  to   Section
2 we state  the notations which  will be used  in  the following sections.
For   a   function $f(\cdot, t)$,  of one spatial  variable and time  $t$, we denote by $f' (\cdot, t)$  and $f''(\cdot, t )$ the first    and   the second derivative with respect to   the spatial variable.  We will denote by $ \dot  f  (\cdot, \cdot) $ the   time   derivative.    Further, we will denote by $ C= C( \beta, J,d) $ a positive constant  which depends only on  these quantities  and which might change from one occurrence to one other.  

\section{Smoothing estimates and differentiability  of the free energy }

In this section we state some technical results upon which our analysis
in the following sections depends. 

It is not hard to show that classical solutions of our equation exist and are unique  by adapting to our non-local equation standard fixed-point
arguments for semi-linear equations, and it is easy to see from a maximum principle argument that these classical solutions satisfy
$|m(x,t)| \leq 1$ for all $x$ and $t$ since the non-local drift term vanishes wherever $m(x,t) = 1$. 

The integration by parts leading from (\ref{gradflow})  to (\ref{dissipate}) is problematic if $m(x,t)$ is not bounded away from $\pm1$. 
Therefore, fix $0 < \lambda <1$, which we shall take increasing to $1$ shortly, and consider the function
$t\mapsto \F(\lambda m  (t))$ where $m$ solves our equation.  Now there is no trouble integrating by parts for any $t>0$, as long as $\nabla m$ is square integrable, and hence
we find 
\begin{eqnarray}
\frac{{\rm d}}{{\rm d}t} \F(\lambda m(t)) &=& \lambda \int_{\R\times\Lambda}  [\frac{1}{\beta} {\rm arctanh}(\lambda m) - \lambda J\star m] \nabla \cdot[ \nabla m(t) - \sigma(m)\nabla J\star m(t)]\dd x\nonumber\\
&=& -\lambda \int_{\R\times \Lambda}  \left[ \frac{\l}{\beta} \frac{|\nabla m(t)|^2}{1 - \lambda^2 m^2(t)}  -  \lambda \frac{1-m^2(t)}{1 - \lambda^2 m^2(t)}\nabla m \cdot \nabla J\star m(t) +\lambda \sigma(m) |\nabla J\star m(t)|^2\right]\dd x  \nonumber
\\
&=&  -\lambda^2 \int_{\R\times \Lambda}  \left[ \frac{1}{\beta} \frac{|\nabla m(t)|^2}{1 - \lambda^2 m^2(t)}  -   \frac{1-m^2(t)}{1 - \lambda^2 m^2(t)}\nabla m \cdot \nabla J\star m(t) + \sigma(m) |\nabla J\star m(t)|^2\right]\dd x. 
\nonumber
\end{eqnarray}

Therefore,
\begin{multline}
 \F(\lambda m(0)) \geq  \F(\lambda m(t)) +\\
\lambda^2  \int_0^t  \int_{\R\times \Lambda}  \left[ \frac{1}{\beta} \frac{|\nabla m(t)|^2}{1 - \lambda^2 m^2(t)}  -   \frac{1-m^2(t)}{1 - \lambda^2 m^2(t)}\nabla m \cdot \nabla J\star m(t) + \sigma(m) |\nabla J\star m(t)|^2\right]\dd x  \ .\nonumber
\end{multline}

Now a simple argument with Fatou's Lemma shows that
$$
{\mathcal I}(m(t)) \leq \liminf_{\lambda \to 1}
\int_{\R\times \Lambda}  
\left[ \frac{1}{\beta} \frac{|\nabla m(t)|^2}{1 - \lambda^2 m^2(t)}  
-2 \lambda \frac{1-m^2(t)}{1 - \lambda^2 m^2(t)}\nabla m \cdot \nabla J\star m(t) 
+\lambda \sigma(m) |\nabla J\star m(t)|^2\right] \dd x \ ,  $$
and that
$$ \F(m(t))  \leq \liminf_{\lambda \to 1}\F(\lambda m(t)) \ .$$
On the other hand, for $m(0)$ strictly and uniformly  bounded away from $\pm1$, it is easy to show that
$$ \F(m(0))  =  \liminf_{\lambda \to 1}\F(\lambda m(0)) \ .$$

Thus, one obtains for every classical solution such that $|\nabla m(x,t)|^2$ is integrable in $x$ for each $t$, and which has
initial data  strictly and uniformly  bounded away from $\pm1$ that
\begin{equation}\label{intfor}
\F(m(t)) + \int_0^t  {\mathcal I}(m(s))  \dd s \leq \F(m(0))\ .
\end{equation}
In Section 7 we shall prove that classical solutions whose initial data is a small $L^2$ perturbation of an instanton do have square integrable gradients
for all positive times.   Note that even if the gradient of $m(t)$  is square integrable for all $t$, it is still possible for
${\mathcal I}(m(t))$ to be infinite, though by (\ref{intfor}), the set of such $t$ must be a null set. 

In what follows, we shall need to integrate by parts frequently, and it will be convenient to use the energy dissipation relation with equality
in place of  the inequality in (\ref{intfor}).  
Our main goal in this section is to explain results showing that if the initial data is a small perturbation of an instanton in the $L^2$ norm, then, after waiting a short time $t_0$ later,
the solution is regular and is a small perturbation of the instanton in Sobolev norms that guarantee that it is also a small perturbation in the $L^\infty$ norm. Thus if we wait a short time, starting from initial data that is a sufficiently small $L^2$  perturbation of an instanton, the solution
$m(t)$ will be strictly bounded away from $\pm1$, at the very least on an open interval to the right if $t_0$, and then we no longer need
the $\lambda$-regularization used above: We can simply integrate by parts to rigorously obtain the {\em identity} (\ref{dissipate}). 
Moreover, we shall know that $\F(m(t_0))$ is still bounded by $\F(m(0)$. We shall also show that the moment bounds on the initial data in our main theorem are effectively  propagated forward to time $t_0$. 

Therefore, the results in this section permit us throughout the rest of the paper to restrict our attention to the behavior of the 
free energy functional $\FF$ on the
set of profiles
\begin{equation}\label{bigmdef}
{\cal M} = 
\{ m: \ |\ \|m\|_\infty  <  1\quad{\rm and}\quad \|m-\bar m\|_2 <\infty  \},
\end{equation}
equipped with the the metric $d(m_1,m_2) = \|m_1-m_2\|_\infty + \|m_1-m_2\|_2$.
Note that this is not a subset of $L^\infty\cap L^2$ since profiles in ${\mathcal M}$ are never
square integrable.   However, ${\cal M}$  is open in this metric topology, and 
 $\FF$ is  Frechet differentiable 
on  ${\cal M}$, and for any differentiable curve $m(t)$  in ${\mathcal M}$, one has
$$\frac{\dd }{\dd t}\FF(m(t)) = \int_\R  \frac{\d  \FF}  {\d  m}  \frac{\partial}{\partial t}m(t)\dd x \ ,$$
where
\begin{equation*}
 \frac{\d  \FF}  {\d  m} = {1\over\beta}{\rm arctanh}(m) -
J\star (m). 
\end{equation*}
The convolution term satisfies $\|J \star  m \|_\infty \le 1$, but
on any set where $m= \pm 1$  we have
$${\rm arctanh} m  = \pm\infty\ .$$

The following theorem summarizes our discussion so far in this section, except for the fact that it remains to be proved that,
as claimed, the conditions on the initial data do indeed ensure square integrability of the gradient.  This is done in Section 7.

\begin{thm}\label{43}  There is a $\delta>0$ such that for all initial data $m_0$ in ${\cal M}$
with $\FF(m_0) < \infty$ and $\|m_0 - \bar m_a\|_2 < \delta$, some $a\in \R$, the corresponding 
solution $m(x,t)$
of  equation  (\ref{1.1}) satisfies 

\begin{equation*}
{\cal F}(m(0)) \geq  {\cal F}(m(t)) + \int_0^t{\cal I}(m(s)){\rm d}s 
 \quad {\rm for\ all}\quad   t > 0 
 \end{equation*}
where    
\begin{equation*}
 {\cal I}(m)  = \int_{\R \times \L} \sigma(m)\biggl( \nabla {\delta \FF\over \delta m}
\biggr)^2 {\rm d}x,  
\end{equation*}
is the quantity defined in  (\ref{8.1}).
 In particular, ${\cal F}(m(t))$ is monotonically decreasing.
\end{thm} 

The proofs of many the results
established in this paper depend  on certain smoothing properties of the evolution
(\ref{1.1}). The required {\it a--priori} smoothing estimates are summarized 
in the following theorem which holds only  when   $D\le 3$. The constraint on the dimensions depends on  the application  of Sobolev estimates  to control  
the $L_\infty$ norm of the gradient of $v$, see Lemma~\ref{lm4}.

\begin{thm}\label{smoo}   Let  $D \le 3$ and  $m(t)$ be
any solution of   (\ref{1.1}). Let $\epsilon>0$, $t_0>0$ and $k\in \N$ be given. Then there is a $\delta= \delta (\e, t_0,k) >0$  and $T_0$ 
such that provided
$$\|m(t) - \bar m_{a(t)}\|_2^2 \le  \delta  \qquad \hbox{\rm for  all} \qquad t \le    T_0\ ,$$
then
$$\sum_{j=0}^k \|(-\Delta)^{j/2}[m(t) - \bar m_{a(t)}]\|_2^2 \le \epsilon   \qquad \hbox{\rm  for all}\quad
t_0 \le t \le  T_0$$
  and also such that
$$\| x_1(m(t) - \bar m_{a(0)})\|_2^2 \leq   2 \| x_1(m(0) - \bar m_{a(0)})\|_2^2\ .$$ 
\end{thm}

The proof is based on several intermediate results and it is given in Section 7.

  \bigskip
  Next we  prove that $ a(\cdot)$ is differentiable and estimate $ \dot a(\cdot)$.   

\begin{thm}\label{5}Let  $m$ be a solution of (\ref{1.1}). Then
there is a $\delta_0 > 0$  so  that whenever
\begin{equation}\label{amin}
\inf_{a\in\R}\{\|m(t) - \bar m_a\|_2\} < \delta_0 
\end{equation}
there is a unique value $a(t)$ at which the infimum in (\ref{amin}) is
attained. Moreover, for any $\kappa>0$, there is a $\delta_1(\kappa,\b,J)$
such that whenever $\| v(t)\|_{W^{s,2}} \le \kappa$ for $s > \frac {D} 2$  and $\|v(t)\|_2 \le 
\delta_1$, $a(t)$ is differentiable and
\begin{equation*}
|\dot a(t)| \le D(\kappa,\b,J)\|v(t)\|_2
\end{equation*}
where $D(\kappa,\b,J)$ is a constant depending only on $\kappa$, $\beta$ 
and $J$.
\end{thm}
\bigskip
\noindent{\bf Proof:} Let $a(t)$ be any minimizer in (\ref{amin}). Clearly there
is at least one  and what we must show is that there is exactly one.
Define $d(b) = \|m(t) - \bar m_b\|_{L^2(\R \times \L)}^2$. 
We have
\begin{equation}\label{GG1}
  d'(b) = 2\int_{\R \times\L} [m (t,x) - \bar m_b(x_1)] \bar m_b' (x_1)  dx, 
  \end{equation}
and 
$$ 2\int_{\R \times\L}  \bar m_b(x_1)  \bar m_b' (x_1)  dx  =0.  $$ 
Further deriving   (\ref{GG1}) and recalling that  $m(t)= \bar m_{a(t)} + v(t) $  we have 
\begin{eqnarray}
  d''(b)& =& -2\int_{\R \times \L}  m(x,t)\bar m''_b(x_1){\rm d}x = 
-2\int_{\R \times \L}  [ \bar m_{a(t)} (x_1) + v(t, x) ]  \bar m''_b(x_1){\rm d}x  \nonumber\\
 & =&  
2\int_{\R \times \L}  \bar m'_{a(t)} (x_1)\bar m'_b(x_1){\rm d}x - 2\int_{\R \times \L}  v(x,t)\bar m''_b(x_1){\rm d}x.\nonumber
\end{eqnarray}
 Hence,
$$d''(b) \ge  2 \int_{\R \times \L} \bar m'_{a(t)} (x_1)\bar m'_b(x_1){\rm d}x - 2
\delta_0\|\bar m''_b \|_{L^2(\R \times \L)}\ .$$
But by continuity, $$\int_{\R \times \L} \bar m'_{a(t)}(x_1)\bar m'_b(x_1){\rm d}x > \frac 12 \|\bar
m'_b\|_{L^2(\R \times \L)}^2 $$ on some interval $(a(t)-c,a(t)+c)$ for some $c$ depending only on
$\beta$ and $J$. Therefore, choose 
\begin{equation*}
\delta_0 \le {\|\bar m'_b\|_{L^2(\R \times \L)}^2\over 4\|\bar m''_b\|_{L^2(\R \times \L)}}
\end{equation*}
and it follows that
$d''(b) > 0$ on $(a(t)-c,a(t)+c)$, and hence there is exactly one critical point of $d(b)$ on 
$(a(t)-c,a(t)+c)$. However, if $b$ is any value with
$$\|m(t) - \bar m_b\|_{L^2(\R \times \L)}  = \|m(t) - \bar m_{a(t)}\|_{L^2(\R \times \L)}$$
it follows that
$$\|\bar m_b - \bar m_{a(t)}\|_{L^2(\R \times \L)} \le 2\|m(t) - \bar m_{a(t)}\|_{L^2(\R \times \L)} \le 2\delta_0\
.$$ But there is a constant $C$ depending only on $\b$ and $J$ so that
$$\|\bar m_b - \bar m_a\|_{L^2(\R)} \ge {(b-a)^2\over C + (b-a)^2}$$
and thus, 
$$ L^d {(b-a)^2\over C + (b-a)^2} \le  \|\bar m_b - \bar m_a\|_{L^2(\R \times \L)}  \le 2\delta_0\ .$$
Decreasing $\delta_0$ if necessary, one can ensure that $|b-a|<c$. Hence 
any putative second minimum must occur within $(a(t)-c,a(t)+c)$ where there
is only the single critical point $a(t)$. Hence there is no other minimum.
This proves that $a(t)$ is a well--defined function under the condition 
(\ref{amin}).  
To show that $a(t)$ is continuously differentiable, we use the Implicit Function  Theorem. Define 
$$f(a,t) :=  \int_{ \R\times \L }(m(t, x) - \bar m_{a}(x_1))\bar m'_{a}(x_1){\rm d}x\ .$$
This is a $C^1$ function on $\R^2$, and in fact even $C^2$. By  what we have proved above, for each $t$,
there is exactly on $a(t)$ so that $f(a(t),t)= 0$, and at no such point does the gradient of $f$
vanish, since the $a$-component of this gradient is non-zero. Hence, by the Implicit Function  Theorem, the curve $t\mapsto (a(t),t)$ is 
continuously  differentiable. Moreover, since this curve is the graph of the function $t\mapsto a(t)$, we have that $a(t)$
is continuously differentiable, as claimed.

  We now  bound $|\dot a(t)|$. 
  Differentiating $f(a(t),t) = 0$ in $t$, one obtains    
\begin{equation*}
\dot a(t)\bigl(\|\bar m'_a\|_{L^2(\R \times \L)}^2 - \langle v,\bar m''_a\rangle_{L^2(\R \times \L)} \bigr)
= -\int_{\R \times \L} {\partial m\over \partial t} \bar m'_{a(t)}\ .  
\end{equation*}
Taking into account   (\ref{gradflow}) and integrating by part   we have 
\begin{eqnarray}
 -\int_{\R \times \L} {\partial m\over \partial t} \bar m'_{a(t)}&=&
 \int_{\R \times \L}   \sigma(m)\nabla  (\frac {\delta \FF}    { \delta  m }  )   e_1 \bar m''_{a(t)}  
 \nonumber\\
  &=&
 \int_{\R \times \L}     \frac {\delta \FF}    { \delta  m }   \nabla \cdot [\sigma(m)  e_1 \bar m''_{a(t)}] \nonumber\\
  &=&   \int_{\R \times \L}     \frac {\delta \FF}    { \delta  m } (\bar m + v)     [\sigma(m)    \bar m''_{a(t)}]_{x_1}.
  \nonumber  
\end{eqnarray}
 Assume that
 \begin{equation*}
  \|v\|_{L^2(\R \times \L)}  \le   \frac 12  \frac { \|\bar m'_a\|_{L^2(\R \times \L)}^2}  { \| \bar m''_a\|_{L^2(\R \times \L)} }. 
 \end{equation*}
  We thus obtain  
 \begin{eqnarray}|\dot a(t)| &  \le& \frac { 2 } { \|\bar m'_a\|_{L^2(\R \times \L)}^2}   \left | \int_{\R \times \L}     \frac {\delta \FF}    { \delta  m }  (\bar m + v)    [\sigma(m)    \bar m''_{a(t)}]_{x_1}  \right | \nonumber\\
  & \le&
   C(\b,J) \frac { 2 } { \|\bar m'_a\|_{L^2(\R \times \L)}^2} \| [\sigma(m)    \bar m''_{a(t)}]_{x_1}  \|_{L^2(\R \times \L)}  \|v\|_{L^2(\R \times \L)}.   
   \nonumber
   \end{eqnarray}
 Taking into account     (\ref{euler}),   we bounded the $L^2$ norm of
$\delta\FF/\delta m(\bar m + v)$    by a constant times the 
$L^2$ norm of $v$ whenever $\| v(t)\|_{W^{s,2}} \le \kappa$  with $s > \frac D 2$ and  $\kappa $ sufficiently small to guarantee that
$\|v\|_\infty \le (1-m_\b^2)/2$. 
 \qed

  \section{Bound   on  the dissipation rate of the free energy in terms of the dissipation rate for the linearized evolution}

In this section we establish
a  bound on the rate ${\cal I}(m(t))$ defined in (\ref{8.1}) at which the excess free energy 
${\cal F}(m(t)) - {\cal F}(\md)$ is dissipated  in terms of  the dissipation rate for the  {\it  linearized evoluiton}, see Theorem~\ref{51}.

To state the main result of this section we need  the following definitions. 
 Denote by $ \BB_a$, $a \in \R$ the family of  linear operators in $ L^2 ( \R \times \L)$,
 \begin{equation}\label{G.2}
  \BB_a v =   \left(
 \frac { v } {\b (1-\overline m_a^2)} -  J \star v     
\right ), 
\end{equation} 
 where   $\overline m_a $, $ a \in \R$ is   the planar front  defined in (\ref{ema1}). 
Denote,   by an abuse of notation,     $ \overline m '_a(x_1) = \frac {\partial} {\partial x_1} \bar m_a(x_1, x_1^\perp)$. 
 It is immediate to verify that 
 \begin{equation}\label{eig1}
   (\BB_a \overline m '_a) (x)  =0, \qquad x \in \R \times \L,  
\end{equation}
      $ \overline m '_a$  is  therefore the eigenfunction corresponding to  the zero  eigenvalue.
 Further  $ \BB_a$ is a  selfadojnt  operator  in $ L^2 ( \R \times \L)$ and   Weil's     theorem,  by the same argument used in \cite{DOPT3}  for the $d=1$  case,   assures the existence of a gap in the  spectrum: 
    For $ v \in  L^2 ( \R \times \L)$ 
 \begin{equation}\label{G.3}
  \int_{\R \times \L} v(x) \overline m_a'(x_1)   { \rm  d } x =0,$$
$$ \langle \BB_a v, v \rangle_{L^2( \R \times \L)} \ge \g(L) \|v\|^2_{L^2( \R \times \L)}, 
\end{equation}
where $\g(L) >0$.    A quantitative argument   given in  \cite{BDDP} proves that  $\g(L)= \frac { a(\beta,J)} {L^2}$.  
To  our aims the merely existence of a gap is enough. 
Denote by $ \AA_a$,  $a \in \R$,  the family of  linear operator in $ L^2 ( \R )$  
 \begin{equation}\label{G.2a} \AA _av :=   \left(
 \frac { v } {\b (1-\overline m_a^2)} -  \overline J \star v     
\right ),   
\end{equation} 
where $\overline J$ is defined in (\ref{tildeJ}). 
 In  \cite{DOPT3}  was shown that   $\AA_a$ has a gap $\g_0$. 
The  eigenfunction corresponding to  the zero  eigenvalue is $ \overline m_a'(\cdot)$.
Then for all $ v \in  L^2 ( \R )$ so that
$$  \int_{\R  } v(x_1) \overline m_a'(x_1)   { \rm  d } x_1 =0$$
 \begin{equation*}
  \langle \AA_a v, v \rangle_{L^2( \R)} \ge \g_0 \|v\|^2_{L^2( \R)} . 
  \end{equation*}
 The operator $ \BB_a$ is the  multidimensional version of the   operator $ \AA_a$. Notice  that   $ \AA_a$  is defined in term of  $ \bar J$ the one dimensional projection of $J$.  Both operators  have  as   eigenfunction corresponding to  the zero  eigenvalue $\bar m'_a$. 
In the following we will 
  drop the subscript $a$ if no confusion arises.  We have:

\begin{thm}\label{51} Let  $D \le 3$,  $m(\cdot,t)$ be a solution of
(\ref{1.1}) and  $ m(\cdot,t)= \overline m_{a(t)}(\cdot) + v(\cdot, t)$ where $a(t)$ is chosen so that minimizes $\|m(t)- \overline m_a\|^2_{L^2(\R \times \L)}$.    For any
  $\e>0$ small enough, 
    there is   $\d_1=\d_1 (\e,\b,J,L) >0 $  
   so that   at all time $t$ for which  $ \|v (t)\|_ {W^{s,2} } \le \d_1  $,   where    $s> \frac {D}2$,        we have that        
  \begin{eqnarray}\label{f2}
      {d\over dt}\bigl[{\cal F}\big(m(t))-{\cal F}\big(\overline
m     \big)\bigr]  &=& -  
  \II \big (m( t)\big ) \nonumber\\             
   \cr & \le& -   (1-3 \e)   \sum_{i\ge1} \int_{\R \times \L}
  \sigma(\overline m)\bigl[(\BB v(t))_{x_i}  \bigr]^2{\rm d}x, \nonumber
\end{eqnarray}
where $  \II \big (m( t)\big )$ is given in (\ref{8.1}).
\end{thm}

The proof of Theorem~\ref{51} is based on several intermediate results.  We start proving the following estimate  for the mobility $ \s
(m)=\b(1-m^2)
$. 

\begin{lm}\label{519}   Set   $ m = \overline m_{a} + v$ where $a$ is chosen so that it  minimizes $\|m - \overline m_a\|^2_{L^2(\R \times \L)}$.  
      For any
  $\e>0$ 
 there exists $\d_1( \e) >0 $ such that    
 \begin{equation*}
  (1-\e ) \s(\overline m) \le   \s(m) \le \s(\overline m)  (1+\e )
\end{equation*}
when   $ \|v\|_ {W^{s,2} } \le \d_1  $ and    $s>\frac {D} 2$.
\end{lm}

\noindent \noindent{\bf Proof:}    Write 
 $ \sigma (m) = \sigma (\overline m)  \left [1+  \frac 1 {\sigma (\overline m) }
 \b (2\overline m+v) v\right].$  One easily obtains the pointwise bound 
$$\left| \frac {1} {\sigma (\overline m) } \b (2 \overline m+v) v\right| \le    2|v| +v^2   \le   2  C(d,s) \|v\|_ {W^{s,2}}, $$
 where in the last passage we   estimated the    $ \|v\|_\infty$ by  Lemma~\ref{A1} for   $s> \frac {D} {2}$. 
 Take $\d_1$ so that 
 $2 C(d,s) \d_1 \le \e$.  \qed

\vskip0.5cm \noindent

\begin{lm}\label{DD1}   Set  $m= \overline m+ v$, $v \in L^2 (\R \times \L)$, $\int \overline m'(x_1) v(x) {\rm d} x=0$.
 For any $\e>0$, for  $s > \frac {D} {2}$   there exists  $ d_1:= d_1(\e,L,\b, d) $,  defined in  condition (\ref{D.3}),  so that    if  $ \|v\|_ {W^{s,2}}  \le d_1$, we have 
 \begin{equation}\label{S.S1} \begin {split}
     {\cal I} (m)    
&\ge   (1-2 \e)  \int_{\R \times \L} \sigma(\overline m)\bigl[(\BB v)_{x_1}  \bigr]^2{\rm d}x  \\ 
& +  \e  \int_{\R \times \L} \sigma(\overline m)\bigl[(\BB v)_{x_1}  \bigr]^2{\rm d}x -
 \frac{1}{ \e}  \int_{\R \times \L} \sigma(\overline m) [U (v)]^2  {\rm d}x  \\ 
 &+    (1-2 \e)   \sum_{i\ge2} \int_{\R \times \L}\sigma(\overline m)\bigl[(\BB v)_{x_i}  \bigr]^2{\rm d}x   
 \end {split} \end{equation}
  where $U (v)$ is defined in (\ref{3.100}). 
  \end{lm}
  
  \noindent 
\noindent{\bf Proof:}  
Since  $\overline m$ depends only on $x_1$  and  by  assumption $  m(x)= \overline m(x_1) + v (x) $ we   decompose  $ {\cal I} (m)$, see  (\ref{8.1}),  as 
 \begin{equation*}
   {\cal I} ( \overline m   + v )=    {\cal I}_1( \overline m  + v  ) +  {\cal I}_2( \overline m  + v  )
\end{equation*}
where
 \begin{equation*}
     {\cal I}_1( \overline m  + v  )=   
  \int \s (m) \left 
 [ \left(
 \frac {m_{x_1 }} {\b (1-m^2)} -  J \star m_{x_1 }   
\right ) \right ]^2  { \rm  d } x ,  
\end{equation*}
and 
 \begin{equation*}
  {\cal I}_2 ( \overline m  + v  )= \sum_{i\ge 2} 
\int \s (m) 
 \left [   \left(
 \frac { v_{x_i } }{\b (1-m^2)} -  J \star v_{x_i } 
\right ) \right ]^2  { \rm  d } x. 
\end{equation*} 
 We have
 \begin{equation*}
  \left(
 \frac {m_{x_1 }} {\b (1-m^2)} -  J \star m_{x_1 }   
\right ) = \left(
 \frac {m_{x_1 }} {\b (1-m^2)}  - 
 \frac {m_{x_1 }} {\b (1-\overline m^2)} \right )  +  \BB  m_{x_1 },  
 \end{equation*}
  where $\BB$ is the linear  operator defined in (\ref{G.2}). 
Denote  
 \begin{equation}\label{L.2}
 \tilde U (v):=  \frac {1}{\b}   \biggl({1\over {1-m^2}} - {1\over {1-\overline m^2}}\biggr) =  \frac {1}{\b}  
 {2\overline m\over (1-\overline m^2)^2}v + {1+3\overline m^2 + 2\overline m v\over (1-\overline
m^2)^2(1-m^2)}v^2. 
\end{equation} 
Since   $\BB  m_{x_1 } = \BB  \overline m'+  \BB  v_{x_1 }= \BB  v_{x_1 }$, by (\ref{eig1}),  we have 
 \begin{equation} \begin {split} \label{L.3}
  &   \frac {1}{\b}    \frac {  m_{x_1 }} {  1-m^2  } -  J \star m_{x_1 }  =  \BB  v_{x_1 }  +  \tilde U (v) ( v_{x_1 } + \overline m' )    \\
   & =  \BB  v_{x_1 }  +  {1\over
\beta}{2\overline m\overline m'\over (1-\overline m^2)^2}v  + U(v)  
 \end {split}   
\end{equation}
 where 
 \begin{equation}\label{3.100}
 U(v) =   \frac {1}{\b}  {2\overline m\over (1-\overline m^2)^2}v v_{x_1 }  +\frac 1
\b {1+3\overline m^2 + 2\overline m v\over (1-\overline m^2)^2(1-m^2)}v^2 ( v_{x_1 } + \overline m' ).   \end{equation}
But since
$${2\bar  m \bar m' \over (1-\overline m^2)^2} = {{\rm d}\over {\rm d}x_1} \left ( {1\over
 1-\overline m^2 } \right )$$  \eqref{L.3}  is the same as
 \begin{equation*} 
 \frac {1}{\b}    \frac {  m_{x_1 }} {  1-m^2  } -  J \star m_{x_1 }  = (\BB v)_{x_1} + U(v)\ .
 \end{equation*}
 Applying  Lemma~\ref{519} we have that
 \begin{equation}\label{3.A20}
   {\cal I}_1(m  +v ) \ge  (1-\e)  \int_{\R \times \L} \sigma(\overline m)\bigl[(\BB v)_{x_1} + U(v)
\bigr]^2{\rm d}x  
\end{equation}
provided $d_1 $ is less than the $ \d_1$ of  Lemma~\ref{519}.
 
We apply to  (\ref{3.A20}),   inequality (\ref{3.B1}) stated in the appendix   with $\l=1-\e$, where $   \e >0  $ is small and  arbitrarily chosen.    We obtain 
 \begin{equation*}  \begin {split}
  &   \int_{\R \times \L} \sigma(\overline m)\bigl[(\BB v)_{x_1} + U(v)
\bigr]^2{\rm d}x  \cr & \ge  (1-2 \e)  \int_{\R \times \L} \sigma(\overline m)\bigl[(\BB v)_{x_1}  \bigr]^2{\rm d}x   \\
     &  +  \e  \int_{\R \times \L} \sigma(\overline m)\bigl[(\BB v)_{x_1}  \bigr]^2{\rm d}x -
 \frac {1}{\e}  \int_{\R \times \L} \sigma(\overline m) [U (v)]^2  {\rm d}x . 
\end{split}\end{equation*}
We proceed  similarly for   
  $   {\cal I}_2 ( \overline m + v   )$, taking in account that $ \overline m$ depends only on $x_1$. 
  We have    
\begin{equation*}
 \frac {1} {\b}  \frac {  v_{x_i }} {  1-m^2  } -  J \star v_{x_i }  =    \BB v _{x_i} + \tilde U(v)   v_{x_i }, \quad i \ge 2 
 \end{equation*}
 where $ \tilde U (v)$ is given in (\ref{L.2}). 
  Then 
 \begin{eqnarray}\label{3.A22}    {\cal I}_2 ( \overline m  + v  )
&  \ge&    (1-2 \e)   \sum_{i\ge2} \int_{\R \times \L}
  \sigma(\overline m)\bigl[ \BB v_{x_i}  \bigr]^2{\rm d}x  \nonumber\\ 
      &  + &\e   \sum_{i\ge2}  \left (  \int_{\R \times \L} \sigma(\overline m)\bigl[ \BB v _{x_i}  \bigr]^2{\rm d}x -
 \frac {1}{ \e }  \int_{\R \times \L} \sigma(\overline m) [\tilde U (v)   v_{x_i} ]^2 {\rm d}x  \right) \ . 
\end{eqnarray}
  Next we show that the last line of (\ref{3.A22})   is positive when $ \|v\|^2_ {W^{s,2}}$ is  small enough. 
  By  periodicity   $ \int_\L  v_{x_j} (x)  {\rm d} x^\perp =0 $ for all $j \ge 2$. 
 This   implies that
 $\int_{\R \times \L} v_{x_j} (x)  \overline m' (x_1){\rm d} x=0$ for $j\ge 2$, therefore  by   (\ref{G.3})
 \begin{equation}\label{D.2}
  \sum_{i\ge2}  \int_{\R \times \L} \sigma(\overline m)\bigl[\BB v_{x_i}  \bigr]^2{\rm d}x  \ge \g(L)^2   \|\nabla^ {\perp}v \|^2_2 .
\end{equation}
We have
 \begin{equation*}
      \sum_{i\ge2}  \int_{\R \times \L} \sigma(\overline m) [\tilde U (v(x))   v_{x_i} ]^2 {\rm d}x   \le  c(\b)  \sup_{x \in \R \times \L} |\tilde U (v)|^2 \|\nabla^{\perp}v \|^2_2  \le  
C(d,\beta) \|v\|^2_ {W^{s,2}} \|\nabla^{\perp}v \|^2_2 , 
\end{equation*}
by  Lemma  \ref{A1}  and   $s > \frac {D} {2}$. 
 Then for any given $\e>0$  take $  d_1:= d_1(L,\b, \e) $ so that  for $  \|v\|^2_ {W^{s,2}}  \le d_1$, see (\ref{3.A22}) and (\ref{D.2}),  
 \begin{equation}\label{D.3}
  \e \g^2(L)  \ge    \frac {1}{\e}  
C(d,\b)  d_1.
\end{equation}
We then obtain (\ref{S.S1})
\qed 
\vskip0.5cm 
 \noindent 
We would like to show that  the quantity on the right hand side of   (\ref{S.S1})  is strictly positive.      There  is  no hope  to show that the second line on   the right hand side of   (\ref{S.S1})  is positive. We cannot expect to control the nonlinear contribution of  the dissipation of the free energy  by only the derivative in the $x_1$ direction.  We need to  take into account  also the gradient in the orthogonal direction of $x_1$.  To this aim we
 denote, see (\ref{S.S1}), 
 \begin{equation}\label{O.12}   G_\e (v):=   \e  \left (     \sum_{i\ge1} \int_{\R \times \L}
  \sigma(\overline m)\bigl[(\BB v)_{x_i}  \bigr]^2{\rm d}x \right )-
 \frac {1}{\e}  \int_{\R \times \L} \sigma(\overline m) [U (v)]^2  {\rm d}x.   
 \end{equation}
  In the next proposition we   show that  $ G_\e (v)$ is positive under smoothing assumptions on $v$.

 \begin{prop}\label{512}   Let   $s > \frac {D} {2}$,   $v  \in W^{s+1,2}(\R \times \L)$,    $\int_{\R \times \L}  v(x) \overline m'(x_1){\rm d}x=0$.  For any $\e>0$  
 \begin{equation}\label{3.B5}
     G_\e (v) \ge 0,    
    \end{equation}
provided  
 \begin{equation}\label{mag1}
    \| v\|_{W^{s+1,2}}^2 \le \e_0  
\end{equation}
 for   $\e_0= \e^{2+ r}$ for     $r= r(L)>0$. 
    \end{prop}
     
The proof  depends on several  intermediate  results, and it is given  at the  end of the section.
\vskip0.5cm 
 \noindent
In Lemma~\ref{DD1} we took  advantage by decomposing $ m = \overline m +v$, $v \in L^2 (\R \times \L)$.   In the following it is  helpful to split 
   $v \in L^2(\R \times \L ) $  in the manner: 
 \begin{equation}\label{P1}
  v (x):= v_1(x_1)  + w(x),    
\end{equation}
where
 \begin{equation*}
  v_1(x_1):= \frac {1 }{L^d} \int_{\L}  v (x_1, x^{\perp}) d  x^{\perp},  \quad   x_1 \in \R.  
\end{equation*} 
By construction,  
 \begin{equation}\label{NN5}
   \int_{\L} w (x_1, x^{\perp}) d  x^{\perp}=0   \quad 
\forall x_1 \in \R.  
\end{equation} 
Further if    $  \int_{\R \times \L} v (x) \overline m' (x_1) d  x =0$ then
 \begin{equation*}
   \int_{\R\times \L } w (x) \overline m' (x_1) d  x = \int_{\R } d x_1 \overline m' (x_1) \int_\L w (x)    d  x^{\perp}_1 = 0,  
 \end{equation*}
 \begin{equation}\label{NN2}  
 \int_{\R} v_1 (x_1) \overline m' (x_1) d  x =0. 
\end{equation}
 Using decomposition (\ref{P1}) we get the following  useful result.

\begin{lm}\label{R.1}     Let $v \in L^2(\R \times \L )$,  $v_{x_1}  \in
L^2( \R \times \L )$, $v=v_1+w$ as       in (\ref{P1}).   We have
 \begin{eqnarray*}
  && \int_{\R \times \L} \sigma(\overline m)\bigl[(\BB v)_{x_1}  \bigr]^2{\rm d}x =
 \int_{\R \times \L} \sigma(\overline m)\bigl[(\AA v_1)_{x_1}  \bigr]^2{\rm d}x  \nonumber\\
  & +&
 \int_{\R \times \L} \sigma(\overline m )\bigl[(\BB w)_{x_1}  \bigr]^2{\rm d}x . 
 \end{eqnarray*}
\end{lm}

\noindent \noindent{\bf Proof:}    Take $v$ as in (\ref{P1})
  \begin{eqnarray}\label{G.1}  
  && \int_{\R \times \L}    \sigma(\overline m)\bigl[(\BB v)_{x_1}  \bigr]^2{\rm d}x \cr &=&  \int \s (\overline m(x_1)) \left 
 [\frac {\partial}{ \partial x_1 }\left(
 \frac  {v_1} {\b (1-\overline m^2)} -  \overline J \star v_1    
\right )  +   \frac {\partial}{ \partial x_1 }\left (
 \frac  {w} { \b (1-\overline m_0^2)} -    J \star w   
\right )  \right ]^2  { \rm  d } x    \nonumber\\
 &=&
 \int \s (\overline m (x_1)) \left 
 [ \frac {\partial}{ \partial x_1 }\left(
 \frac  {v_1} {\b (1-\overline m^2)} -  \overline J \star v_1    
\right )     \right ]^2  { \rm  d } x    \nonumber\\ 
& +&
 \int \s (\overline m_0(x_1)) \left 
 [\frac {\partial}{ \partial x_1 }\left(
     \frac  {w} { \b (1-\overline m^2)} -    J \star w   
\right )  \right ]^2  { \rm  d } x     \nonumber\\
 & +&
2  \int \s (\overline m (x_1)) 
 \left [  \frac {\partial}{ \partial x_1 }\left(
 \frac  {v_1} {\b (1-\overline m^2)} -  \overline J \star v_1    \right ) 
  \right ]   \left  [ \frac {\partial}{ \partial x_1 } \left (
 \frac  {w} { \b (1-\overline m ^2)} -    J \star w   
\right )  \right ]   { \rm  d } x.    \nonumber\\
\end{eqnarray}
Integrating  per part with respect to $x_1$   the last term in (\ref{G.1}) we have
\begin{eqnarray*} &&  2  \int \s (\overline m(x_1)) 
 \left [  \frac {\partial}{ \partial x_1 }\left(
 \frac  {v_1} {\b (1-\overline m ^2)} -  \overline J \star v_1    \right ) 
  \right ]   \left  [ \frac {\partial}{ \partial x_1 } \left (
 \frac  {w} { \b (1-\overline m ^2)} -    J \star w   
\right )  \right ]   { \rm  d } x    \nonumber\\ & =&
- 2  \int   \frac {\partial }{ \partial {x_1}  }\left [  \s (\overline m (x_1)) 
  \frac {\partial }{ \partial {x_1}  }\left(
 \frac  {v_1} {\b (1-\overline m^2)} -  \overline J \star v_1    \right ) \right ] 
    \left  [  
 \frac  {w} { \b (1-\overline m ^2)} -    J \star w   
  \right ]   { \rm  d } x_1  { \rm  d } x_1^\perp  \nonumber\\ & =&
  - 2  \int _\R   \frac {\partial }{ \partial {x_1}  }\left [  \s (\overline m (x_1)) 
  \frac {\partial }{ \partial {x_1}  }\left(
 \frac  {v_1} {\b (1-\overline m^2)} -  \overline J \star v_1    \right ) \right ] 
  \left ( \int_\L \left  [  
 \frac  {w} { \b (1-\overline m ^2)} -    J \star w   
  \right ]     { \rm  d } x_1^\perp \right )  { \rm  d } x_1  =0
   .\nonumber\\
 \end{eqnarray*}
 Namely for each $x_1$   integrating with respect to $x_1^\perp$ we have 
$$  \int_ \L \left  [  
 \frac  {w} { \b (1-\overline m^2)} -    J \star w   
  \right ]      { \rm  d } x_1^{\perp } = 0. $$
   \qed

 Taking into account Lemma  \ref{R.1}   we  write, see  (\ref{O.12}),  
 \begin{equation}\label{O.13} \begin {split}
  &   G_\e (v)=  \e  \left [  \int_{\R \times \L} \sigma(\overline m)\bigl[(\AA v_1)_{x_1}  \bigr]^2{\rm d}x    +\sum_{i\ge2} \int_{\R \times \L}
  \sigma(\overline m)\bigl[ \BB (w_{x_i})  \bigr]^2{\rm d}x   \right ]   \\
  &+       \e  \int_{\R \times \L} \sigma(\overline m)\bigl[(\BB w)_{x_1}  \bigr]^2{\rm d}x 
-   \frac {1}{\e}  \int_{\R \times \L} \sigma(\overline m) [U (v)]^2  {\rm d}x.  
\end {split} \end{equation}
To  show  Proposition  \ref{512}     we  bound from below  the first three  terms of   (\ref{O.13})  and from above the last term of  (\ref{O.13}) in term of   comparable quantities.    
We   estimate  the first  and third  term  of  (\ref{O.13}) in Lemma~\ref{510} and Lemma~\ref{510a}.   The lower bound  for   the second term in  (\ref{O.13}) is easily obtained  
taking into account   that $ \int_{\L \times \R} w_{x_i} \overline m' (x_1) d x=0$  and applying (\ref{G.3})
 \begin{equation}\label{O.15}
     \sum_{i\ge2} \int_{\R \times \L}
  \sigma(\overline m)\bigl[ \BB (w_{x_i})  \bigr]^2{\rm d}x   \ge  \g(L)   \sigma(m_\b)  \sum_{i\ge2} \|w_{x_i}\|^2_2. 
\end{equation}
Then we estimate from above in term of the same quantities  the $U(v)$ term, see Lemma \ref{511}. 

\begin{lm}\label{510} Let $v \in L^2(\R \times \L )$,  $v_{x_1}  \in
L^2( \R \times \L )$, $\int_{\R \times \L } v(x)
\overline m' (x_1) dx =0 $  and $v(x)= v_1(x_1) + w(x)$, see decomposition (\ref{P1}),  then  there exists a positive constant $\g_1:= \g_1 (\b, J) $,
  such that 
  \begin{equation*}
    \int_{\R \times \L}   \sigma(\overline m)  \left [ \left ( \frac {\partial   } {\partial x_1} (\AA  v_1 )\right) \right ]^2   {\rm d}x
  \ge   \s(m_\b)  \g_1\|Pv_1'\|_{L^2(\R \times \L)}^2 
\end{equation*} 
 where  $\AA  $ is the linear operator defined in (\ref{G.2a}),  $ P$ is
the orthogonal projection on  the orthogonal complement of 
$    {\overline m}''  $ of $L^2(\R)$. 
\end{lm}

\noindent  \noindent{\bf Proof:}     We apply    Lemma 3.4 of \cite{CCO1}. 
The  assumption  needed is  $ \int_ \R v_1(x_1)
\overline m' (x_1)  { \rm  d }  x_1 =0 $  which     is indeed satisfied;
see (\ref{NN2}).

We then obtain, from  Lemma 3.4 of  \cite {CCO1},  that there exists a positive constant $\g_1$ depending on $\beta$ and $J$ so that 
 \begin{eqnarray*}  \int_{\R \times \L} \sigma(\overline m)\bigl[(\AA v_1)_{x_1}  \bigr]^2{\rm d}x   &\ge& \s(m_\b)   \g_1  \int_\L   {\rm d} x_1^\perp  \int_\R {\rm d} x_1 [(P v_1')(x_1)]^2 \nonumber\\ & =&  \s(m_\b)  \g_1\|Pv_1'\|_{L^2(\R \times \L)}^2, \nonumber
\end{eqnarray*}
   where  $ P$ is
the orthogonal projection on  the orthogonal complement of 
$ {\overline m}'' $  in $L^2(\R)$, i.e.   
 \begin{equation*}
 Pv_1'= v_1'- \overline m'' \frac { \int_\R   v_1'(x_1) \overline m''(x_1) {\rm d} x_1} {\|\overline m''\|^2_{L^2(\R)}} .
\end{equation*}
\qed

\vskip0.5cm \noindent 

 Next,  we estimate from below the term $\int_{\R \times \L}   \sigma(\overline m)  \left [ \left ( \frac {\partial   } {\partial x_1} (\BB  w )\right) \right ]^2   {\rm d}x$.     When dealing with the  heat equation in our heuristic discussion,  the  corresponding term   was simply dropped. Now we need  to   bound it   from below to  get some positive contribution that may be used to cancel  negative contributions  coming  from the  last term of \eqref {O.13}.  
 The estimate is obtained  by introducing a cut-off function.  Without cut-off we could get an estimate  of the type:
\begin {equation}  \label {eq:t7} \int_{\R \times \L}   \sigma(\overline m)  \left [ \left ( \frac {\partial   } {\partial x_1} (\BB  w )\right) \right ]^2   {\rm d}x \ge       \s (m_\b)  \g^2(L) \| w_{x_1}\|^2_2   -  C \|w\|_2^2.          \end{equation}
The main difference  between this and    \eqref {3.1002}  is that    the term $\|w\|_2^2 $  in  \eqref {eq:t7}  is a priori not small and we do not have a way to control it. 
 
Let  $N \ge 1$  and  $  \phi^2_N (x_1) $,  $x_1 \in \R$    be a smooth cut-off function   so that 
 \begin{equation}\label{O.1}
  \phi^2_N (x_1)= \begin{cases}  0 &  |x_1|  \le N\\  1 &   |x_1|  \ge 2N \end{cases} 
  \end{equation}
and   
 \begin{equation}\label{O.1a}
  |\phi_N (x_1)| \le 1,  \quad  |\phi'_N (x_1)| \le \frac 1 N, \quad \quad  |\phi''_N (x_1)| \le \frac {1} {N^2}.  
\end{equation} 
 The  choice of cut-off  $N$ will  depend  on $L$,   the linear size of the transversal direction to the front, and it will be chosen as     function of $\e$,   see  proof of  Proposition~\ref{512}.  
We have the following. 
\vskip0.5cm 
\begin{lm}\label{510a}  Take  $v \in L^2(\R \times \L )$,  $\int_{\R \times \L } v(x)
\overline m' (x_1) dx =0 $,    $v= v_1 + w$, see decomposition (\ref{P1}), $ w_{x_1} \in L^2(\R \times \L )$  and   $\phi^2_N$ the cut-off  function defined in (\ref{O.1}).  Then  for any  $N\ge 1$,  
  \begin{equation}\label{3.1002}
  \int_{\R \times \L}   \sigma(\overline m)  \left [ \left ( \frac {\partial   } {\partial x_1} (\BB  w )\right) \right ]^2   {\rm d}x  \ge   \frac 1 4    \s (m_\b)  \g^2(L) \|   \phi_N (w)_{x_1}   \|^2_2   -   \|w\|_2^2   \frac{ 1} {N^2}
     D(\b, \g(L))
     \end{equation}
 where  $\BB  $ is the linear operator defined in (\ref{G.2}) and  $D(\b,\g(L))$ is defined in (\ref{O.5}).   
 \end{lm}    

 \noindent \noindent{\bf Proof:}   
 We have that  
\begin{equation} \label {eq:t5a}  \int_{\R \times \L}   \sigma(\overline m)  \left [ \left ( \frac {\partial   } {\partial x_1} (\BB  w )\right) \right ]^2   {\rm d}x  \ge \int_{\R \times \L}   \sigma(\overline m)  \left [  \phi_N   \left ( \frac {\partial   } {\partial x_1} (\BB  w )\right) \right ]^2   {\rm d}x .
\end{equation} 
  Using that  for smooth integrable functions  $g$ and $h$ one has
   $$ \int_{\R} [(gh)']^2= \int_{\R} [gh']^2- \int_{\R} g g''h^2 $$
   we have 
\begin{equation} \label {eq:t4} \begin {split} &   \int_{\R \times \L}   \sigma(\overline m)  \left [  \phi_N   \left ( \frac {\partial   } {\partial x_1} (\BB  w )\right) \right ]^2   {\rm d}x    \\
& = 
\int_{\R \times \L}   \sigma(\overline m)  \left [  \left ( \frac {\partial   } {\partial x_1}  \phi_N       (\BB  w )\right)  \right ]^2   {\rm d}x   + 
\int_{\R \times \L}   \sigma(\overline m)     \phi_N    \phi_N ''  \left [\BB  w   \right ]^2   {\rm d}x.
\end {split} \end{equation}
 By the property of $ \phi_N $, see (\ref{O.1a}), the last term in  \eqref {eq:t4}   is estimated     as following 
 $$\left | \int_{\R \times \L}   \sigma(\overline m)     \phi_N    \phi_N ''  \left [\BB  w   \right ]^2   {\rm d}x \right |   \le \sup |  \phi_N '' | \|\BB  w \|^2_2 \le \frac {1 }{N^2} \|\BB  w \|^2_2  \le   \frac {1 }{N^2}  C (\beta) \|w \|^2_2 ,$$
 where in the last inequality  we  used  the fact   that $\BB$ is a bounded operator in $L^2$.
We then obtain that, see \eqref {eq:t5a},
\begin{equation} \label {eq:t5}  \int_{\R \times \L}   \sigma(\overline m)  \left [ \left ( \frac {\partial   } {\partial x_1} (\BB  w )\right) \right ]^2   {\rm d}x 
 \ge  \int_{\R \times \L}   \sigma(\overline m)  \left [  \left ( \frac {\partial   } {\partial x_1}  \phi_N       (\BB  w )\right)  \right ]^2   {\rm d}x   -  \frac {1 }{N^2}  C(\beta)\|w \|^2_2 .
\end{equation} 
Next  we  estimate  the first term  on the right hand side of   \eqref {eq:t5}.   We write  
$$  \phi_N \BB  w= \BB   ( \phi_N w)+   \phi_N J\star w-  J \star (\phi_N w),  $$
and  apply     the inequality   see Lemma~\ref{A10}  in the Appendix, writing  $\lambda$ in  (\ref{3.B1}) as $\l = \frac 12 $,
$$ (a+b)^2 \ge   \frac 12 a^2-   b^2. $$
We have
 \begin{equation}\label{O.2} \begin {split}
   & \int_{\R \times \L}   \sigma(\overline m)  \left [  \frac {\partial   } {\partial x_1}   \left ( \phi_N        \BB  w \right)  \right ]^2   {\rm d}x    \cr  &\ge 
 \frac 12  \int_{\R \times \L}   \sigma(\overline m)  \left [  \left ( \frac {\partial   } {\partial x_1}        \BB  (\phi_N   w)  \right)  \right ]^2   {\rm d}x \\
  &-     \int_{\R \times \L}   \sigma(\overline m) \left [ J \star (\phi_N   w) -  \phi_N  J \star w \right ]^2.  \end {split}
\end{equation}
Further, we have  
 \begin{equation}\label{S.5} \begin {split} &
 \int_{\R \times \L}   \sigma(\overline m)  \left [  \left ( \frac {\partial   } {\partial x_1}        \BB  (\phi_N   w)  \right)  \right ]^2   {\rm d}x \cr & =  
\int_{\R \times \L}  \sigma(\overline m) \left [   \BB  \left ( \frac {\partial   } {\partial x_1} (w \phi_N  )\right)    + w(x) \phi_N  \frac {\partial   } {\partial x_1} \left ( \frac {1 } {\b (1-\overline m^2)} \right )   \right ]^2   {\rm d}x  \\
 & \ge 
 \frac 12\int_{\R \times \L}  \sigma(\overline m) \left [    \BB  \left ( \frac {\partial   } {\partial x_1}(w\phi_N )\right) \right ]^2   {\rm d}x \\
      &-  \int_{\R \times \L}  \left [  \phi_N w(x)   \frac {\partial   } {\partial x_1} \left ( \frac {1 } {\b (1-\overline m^2)} \right )   \right ]^2.     
\end {split}  \end{equation}
 By  (\ref{NN5}), 
 \begin{equation*} \int_{\L \times \R}  \left ( \frac {\partial   } {\partial x_1} (\phi_N (x_1) w (x)) \right )\bar m'(x_1) {\rm d}x = 
-\int_{\L \times \R}  \phi_N (x_1) w  (x) \bar m''(x_1) {\rm d}x =0,
\end{equation*}
 and   therefore, by (\ref{G.3}) 
 $$ \int_{\R \times \L}  \sigma(\overline m) \left [    \BB  \left ( \frac {\partial   } {\partial x_1}(w\phi_N )\right) \right ]^2   {\rm d}x  \ge \s (m_\b)  \g^2(L) \|  (\phi_N w)_{x_1}   \|^2_2 . $$
 Taking into account that
 $$ \frac {\partial   } {\partial x_1} \left ( \frac {1 } {\b (1-\overline m^2)} \right ) =    \frac {2 \bar  m \bar m'  } {\b (1-\overline m^2)^2} $$
 and $\bar m' $ is exponential decreasing to zero, see (\ref{decay}),  we have that
$$  \int_{\R \times \L}  \left [  \phi_N w(x)   \frac {\partial   } {\partial x_1} \left ( \frac {1 } {\b (1-\overline m^2)} \right )   \right ]^2   \le \|w\|^2_2 C(\beta) e^{-\alpha N}.$$  
  Therefore 
    (\ref{S.5})  can be estimated as follows:     
     \begin{equation}  \begin{split}\label{LL7}  &
\int_{\R \times \L}   \sigma(\overline m)  \left [   \frac {\partial   } {\partial x_1}       \BB   (\phi_N   w)    \right ]^2   {\rm d}x   \cr&  \ge  
 \frac 12 \s (m_\b)  \g^2(L) \|  (\phi_N w)_{x_1}   \|_2^2   -  \|w\|^2_2 C(\beta) e^{-\alpha N }.  
    \end{split} \end  {equation}
Since, see (\ref{O.1a}) 
 $$ \left |\left ( J \star (\phi_N   w)  -  \phi_N  J \star w\right) (x)\right |  \le \int_{\R \times \L}  J (x-y) \left | w (y) [ \phi_N(y_1)-  \phi_N(x_1)  ] \right | dy   \le \frac 1 N  \int_{\R \times \L}  J (x-y) |w (y)| dy, $$
the last term  of  \eqref {O.2} is estimated as the following:
 \begin{equation}\label{Su.1}
  \int_{\R \times \L}   \sigma(\overline m) \left [ J \star (\phi_N   w) -  \phi_N  J \star w \right ]^2  \le   \frac 1 {N^2} \int_{\R \times \L}  (J \star |w|)^2  \le   \frac 1 {N^2} \|w\|^2_2. 
\end{equation}
Taking into account \eqref {LL7} and \eqref {Su.1}  we estimate  (\ref{O.2})  as follows: 
 \begin{equation} \label {eq:t6} \begin {split}
  &   \int_{\R \times \L}   \sigma(\overline m)  \left [      \frac {\partial   } {\partial x_1}  \left ( \phi_N  \BB  w  \right) \right ]^2   {\rm d}x \\
 & \ge     \frac 14  \s (m_\b)  \g^2(L) \|  (\phi_N w)_{x_1}   \|^2_2   -   \|w\|^2_2  C(\beta)   e^{-\alpha N }    -    \frac 1 {N^2} \|w\|^2_2   \\
  & \ge 
\frac 14   \s (m_\b)  \g^2(L) \|  (\phi_N w)_{x_1}   \|^2_2   -    \|w\|^2_2 \left [    C(\beta)  e^{-\alpha N }  +    \frac 1 {N^2} \right ] 
 , \end {split}
 \end{equation}
 Further 
 $$  \|  (\phi_N w)_{x_1}   \|^2_2 =  \|   \phi_N (w)_{x_1}   \|^2 -  \int_{\R \times \L}   \phi_N   \phi_N'' w^2  \ge \|   \phi_N (w)_{x_1}   \|^2_2 -\frac 1  {N^2} \|w\|^2_2 .$$
Finally, from \eqref {eq:t5}  and  \eqref {eq:t6} obtain
\begin{eqnarray}
  &&\int_{\R \times \L}   \sigma(\overline m)  \left [  \phi_N   \left ( \frac {\partial   } {\partial x_1} (\BB  w )\right) \right ]^2   {\rm d}x \nonumber\\
& \ge&   \frac 14  \s (m_\b)  \g^2(L) \|   \phi_N (w)_{x_1}   \|^2_2   -   \|w\|^2_2 \left [ \frac 14   \s (m_\b)  \g^2(L) \frac 1 {N^2} +   C(\beta)   e^{-\alpha N } +   \frac {C(\beta)} {N^2} +     \frac 1 {N^2} \right ]  \nonumber\\
 & \ge&
     \frac 14   \s (m_\b)  \g^2(L) \|   \phi_N (w)_{x_1}   \|^2_2   -   \|w\|^2_2   \frac 1{N^2}
     D(\b, \g(L)), \nonumber
     \end{eqnarray}
     where    
 \begin{equation}\label{O.5}
      D(\b, \g(L)):=  \left [ \frac 14  \s (m_\b)  \g^2(L)   +2 C(\beta)        \right ]. 
     \end{equation}
  \qed 
  
  \vskip0.5cm 
In the next lemma we   estimate  
  $\int\sigma(m)\bigl[U(v)\bigr]^2{\rm d}x $     from above   in two different ways   which will be used in different regimes.

\begin{lm}\label{511}  
Let $v \in L^2(\R \times \L )$,   $\int_{\R \times \L}  v(x) \overline m'(x_1){\rm d}x=0$, $v  \in W^{s+1,2}(\R \times \L)$,  $s > \frac {D} {2}$ , 
  $v=v_1+w$, as in (\ref{P1}),  and $\phi_N$  be the cut-off function defined in (\ref{O.1}).  For the non linear operator  $U(v)$ defined in  (\ref{3.100}) the
following  holds 
\begin{equation}\label{3.102} \begin {split}
  \int\sigma(m)\bigl[U(v)\bigr]^2{\rm d}x  \le &
  ||v'_1||^2_{L^2(\R \times \L) }  C(\b,J) \left \{ \| v\|_{W^{s,2}}^2    +  \| w_{x_1} \|_2^2 N   \right \}  \\ &+ 
    \|w\| ^2_2C(\b,J) \left \{  \| v\|_{W^{s,2}}^2  +  \| v \|_{W^{s+1,2}} ^2 \right \}  \\  &  +   \| v\|_{W^{s,2}}^2    C(\b,J)  \|\phi_N w_{x_1}\|^2_2 .  
 \end {split}   \end{equation}
 Further, assume that  $\|v_1\|_{L^2(\R\times \L )}^2 \le k^2$ and  $\|v'_1\|_{L^2(\R\times \L )}^2 \le k^2$.   Then for  any given $\e_1>0$  there exists    $\l_0= \l_0(\e_1, k)$, see   (\ref{As3}), so that 
 \begin{equation} \label{3.1003} \begin {split} \int\sigma(m)\bigl[U(v)\bigr]^2{\rm d}x  
   &\le 
   \left [    C(\b, J)   +\frac 2 {\l_0}        \right ]  ||v'_1||^4_{L^2(\R\times \L )}   \\
    &+ 
     \| w\|^2_{L^2(\R\times \L )}    C(\b,J,d)  \left [  \| v\|_{W^{s,2}}^2       +         (\frac N {L^d})^2      \| v\|_{W^{s+1,2}}^2   \right ] \\
      &+  8   \e_1  ||P v'_1||^2_{L^2 (\R )}    +   
  2   \| v\|_{W^{s,2}}\|\phi_N w_{x_1}\|^2 ,
 \end {split}   \end{equation}
where $P$ is the orthogonal projection on the orthogonal complement of $\overline m''$ in $L^2(\R)$.
\end{lm}

   \noindent{\bf Proof:} 
 Observe that for some constant $C$  depending only on $\beta$ and $J$  
 \begin{equation}\label{3.AA1}  \begin {split} |U(v)|^2   &\le   \left (  \frac {1} {\b}  \frac{2\overline m}{ (1-\overline m^2)^2}\right )^2 v^2 v^2_{x_1 }  +    
 \left ( \frac {1}{ \b} {1+3\overline m^2 + 2\overline m v\over (1-\overline m^2)^2(1-m^2)}\right)^2 v^4   v^2_{x_1 } \\
 &+  \left (  \frac {1}{\b}   {1+3\overline m^2 + 2\overline m v\over (1-\overline m^2)^2(1-m^2)}   \overline m'  \right )^2 v^4  \\
  & \le    C (\b, J) \bigl(    R(x_1)|v|^4 + |v|^2|v_{x_1}|^2\bigr) 
 \end {split} \end{equation}
where
$R(\cdot)$ is non-negative, exponentially decreasing to zero as $|x_1| \uparrow \infty$ and $\int_{\R}    R(x_1) {\rm d}x_1 =1$.  We start deriving (\ref{3.102}). 
Splitting $v=v_1+w$ as in \eqref {P1}, we have  
$$v^4 =v^2 [v_1+w]^2 \le  2v^2   [v_1^2+w^2].$$
Then since
$\|v\|_\infty   \le  c(d,s) \| v\|_{W^{s,2}}   $, for $s > \frac {D} 2 $,  see Lemma \ref {A1},
  \begin{equation}\label{Ma.1}
 \int_{\R \times \L}  R(x_1)|v|^4{\rm d}x \le   2 c(d,s) \| v\|_{W^{s,2}}^2  \int_{\R \times \L}  R(x_1)[ v_1^2(x_1) +w^2(x)] {\rm d}x .
\end{equation}
We may write
$$v_1(x_1) = v_1(y)+ \int_{y}^{x_1} v'(z){\rm d}z.$$  
 We then multiply both terms by $\overline m' (y)$  and integrate  on the real line.   Since $\int v_1(y)
\overline m' (y) dy =0 $   we have 
 \begin{equation}\label{3.B2}
 v_1(x_1) = \frac {1} {2m_\b} \int_{-\infty}^{\infty} \overline m' (y) \left (
\int_{y}^{x_1}v'_1(z) {\rm d}z  \right ) {\rm d}y 
\end{equation}
and therefore 
\begin{equation}\label{3.B10} \begin {split}
   |v_1(x_1)|  &\le \frac{ 1} {2m_\b} \left (  \int \overline
m'(y)|x_1-y|^{\frac 12 } {\rm d} y
\right )  ||v'_1||_{L^2(\R)}  \\& \le   \left ( \frac{1} {2m_\b} \int \overline
m'(y)|x_1-y|^2  {\rm d} y
\right )^{\frac 14 }  ||v'_1||_{L^2(\R)}     \le C(\b,J)  [ 1+ |x_1|^2]^{\frac 1 4}  
||v'_1||_{L^2(\R)}. \end{split} \end{equation}
 Therefore from (\ref{Ma.1}) and (\ref{3.B10}) we have
 \begin{equation}\label{Ma.10} \begin {split}   
&\int_{\R \times \L}  R(x_1)|v|^4{\rm d}x \le \\    &2 c(d,s)  \| v\|_{W^{s,2}}^2    \left [ ||v'_1||^2_{L^2(\R) } C(\b,J)  \int_{\R \times \L }  R(x_1)(   (1+ |x_1|^2)^{\frac 12}  {\rm d}x +  \int_{\R \times \L}  R(x_1) w^2(x)  {\rm d}x \right ]  \leq   \\
 &\| v\|_{W^{s,2}}^2    C(\b,J,s,d)  \left [  ||v'_1||^2_{L^2(\R \times \L) }   +  \|w\| ^2_{L^2(\R \times \L) }   \right ].  
 \end{split} \end{equation}
 To estimate the contribution from the
last term  term  in \eqref {3.AA1} we   split $v_{x_1}=v_1' + w_{x_1}$ where $v_1$ and $w$ as in \eqref {P1}, obtaining
  \begin{equation}\label{3.A7} \begin {split}     
  \int |v|^2 |v_{x_1}|^2 {\rm d}x  & \le  2  \left [ \int |v|^2 |v'_1 |^2 {\rm d}x + \int |v|^2 |w_{x_1} |^2 {\rm d}x \right ]  \cr   &\le 2   \|v\|^2_{\infty} \|v'_1 \|^2_{L^2(\R \times \L)}  +2 \int |v|^2 |w_{x_1} |^2 {\rm d}x.
 \end{split}  \end{equation}
  Splitting again $v =v_1 + w$  as in \eqref {P1} we estimate the last term of \eqref {3.A7}
 \begin{equation}\label{eq:t9} \begin {split}    \int |v|^2 |w_{x_1} |^2 {\rm d}x  &\le 2 \left[ \int |v_1|^2 |w_{x_1} |^2 {\rm d}x  + 
   \int |w|^2 |w_{x_1} |^2 {\rm d}x  \right ] \cr & \le
   2 \left[ \int |v_1|^2 |w_{x_1} |^2 {\rm d}x  + 
  \|w\|^2_2   \|w_{x_1} \|^2_\infty   \right ] 
     \end{split}  \end{equation}
The first term of \eqref {eq:t9} is estimated by adding and subtracting the cut-off function $ \phi_N^2$.   Taking into account  (\ref{3.B10}), we have 
   \begin{equation} \label{O.9} \begin {split}     
       &   \int |v_1|^2 |w_{x_1} |^2 {\rm d}x   =      \int |v_1|^2 | \phi_N w_{x_1} |^2 {\rm d}x +  \int |v_1|^2  (1-\phi^2_N) |w_{x_1} |^2 {\rm d}x  \\
       &\le 
  \|v_1\|^2_\infty \|\phi_N w_{x_1}\|^2_2  + 
   c(\b,J) \|v_1'\|_{L^2(\R)}^2  \int  (1+ |x_1|^2)^{\frac 1 2}    (1- \phi^2_N) |w_{x_1} |^2 {\rm d}x   \\ 
   &\le    \|v_1\|^2_\infty \|\phi_N w_{x_1}\|^2_2  + 
  \frac {c(\b,J)} {L^d}   \|v_1'\|_{L^2 (R \times \L)}^2   \|w_{x_1} \|_{L^2 (\R \times \L)} ^2 \sqrt {1+N^2}. 
  \end{split}  \end{equation}
  Summarizing the previous estimate we have that    (\ref{3.A7}) is bounded as following 
  \begin{equation} \label{O.3} \begin {split}    
  &    \int |v|^2 |v_{x_1}|^2 {\rm d}x    \le 2   \|v\|^2_\infty  \|v'_1 \|^2_{L^2 (R \times \L)}  + 4   \|w_{x_1} \|_{\infty} ^2 \|w\|^2_{2}      \\
   &  +  4  \|v\|^2_\infty  \|\phi_N w_{x_1}\|^2_2 +
 \frac {C(\b,J)} {L^d }  \|v_1'\|_{L^2 (R \times \L)}^2\||w_{x_1} \|_{L^2 (\R \times \L)} ^2 \sqrt {1+N^2}.  
  \end{split}  \end{equation}
 Hence from \eqref {3.AA1},  (\ref{Ma.10})  and (\ref{O.3})  we have
 \begin{eqnarray}\label{maa.3}  \int\sigma(m)\bigl[U(v)\bigr]^2{\rm d}x   &\leq&
   \| v\|_{W^{s,2}}^2    C(\b,J)  \left [   ||v'_1||^2_{L^2(\R \times \L) }   +  \|w\| ^2_{L^2(\R \times \L) }   \right ]\nonumber\\
   & +&  \| v\|_{W^{s,2}}^2    C(\b,J)  \|\phi_N w_{x_1}\|^2_{L^2 (\R \times \L)}  +  4 \|w_{x_1} \|_{\infty} ^2 \|w\|^2_{L^2 (R \times \L)}    \nonumber\\
    &  + &  \frac {C(\b,J)} {L^d } 
  \|v_1'\|_{L^2 (\R \times \L)}^2  \| w_{x_1} \|_{L^2 (\R \times \L)} ^2 \sqrt {1+N^2}. 
\end{eqnarray}
 We  estimate, see Lemma~\ref{A1} in the Appendix,    
  \begin{equation*}  \|v\|_\infty  \le c(d,s) \| v\|_{W^{s,2}}, \quad s > \frac {D} 2,  \quad \hbox {and} \quad 
  \quad \||w_{x_1} \|_\infty \le c(d,s)
  \| v_{x_1}\|_{W^{s,2}}\le   c(d,s) \| v \|_{W^{s+1,2}}. 
  \end{equation*}
 Hence     \eqref {maa.3}  immediately implies   (\ref{3.102}). 
 We   next  derive (\ref{3.1003}). 
For $\a \in \R$  write  
 \begin{equation}\label{3.B20}
 v'_1= \a \overline m''+g'
\end{equation}
 where
$ 
\int_\R g' (x_1) m''(x_1) {\rm d} x_1 =0  $  so that $Pv'_1= g'$. Note that, as
indicated in our notation,  $ Pv'_1$ is a derivative since 
$v'_1$ and
$  
\overline m'' $ are derivatives.   
Hence, upon integration 
 \begin{equation}\label {3.BB20}
 v_1 = \a \overline m' +g. 
\end{equation} 
The fact that $\int_\R v_1(x_1)\overline m' (x_1) dx_1 =0 $ means that $\|g\|_{L^2 (\R)}$ cannot be too small.
But what we need to know is that $\|g\|_{L^2 (\R)}= \|Pv'_1\|_{L^2 (\R)}$ is not  too small. 
In general, these are simply two different things.
What provides the crucial connection here is that 
$(2m_\b)^{-1} \overline m'(x_1) {\rm d} x_1$ is
a probability measure on $\R$, so that 
$$\int_\R v_1(x_1)\overline m' (x_1) dx_1 =   0 $$ 
implies that
$$ ||g||_\infty  \ge \frac {|\a|} {2m_\b}||\overline m'||_{L^2 (\R )}. $$
Then one may use 
$ \|g\|^2_\infty \le 2  \|g'\|_{L^2 (\R )} \|g \|_{L^2 (\R )} $ to conclude that
$$|\alpha|^2 \le {8m_\b^2\over \|\overline m'\|_{L^2 (\R )}^2}\|g\|_{L^2 (\R )}\|g'\|_{L^2 (\R )} =
{8m_\b^2\over \|\overline m'\|_{L^2 (\R\times \L )}^2}\|g\|_{L^2 (\R\times \L )}\|g'\|_{L^2 (\R\times \L )} 
 .$$
Since
 \begin{equation*}
  v^4 = (v_1+w)^4 \le  4 [ v_1^2+w^2]^2 \le 8 [ v_1^4+w^4] 
 \end{equation*}
we have that
 \begin{equation}\label{LL.20}
   \int_{\R \times \L}  R(x_1)|v|^4{\rm d}x   \le 8  \int_{\R \times \L}  R(x_1)\left [ |v_1|^4 + w^4 \right] {\rm d}x.  
\end{equation} 
For the first term of (\ref{LL.20}) we
  insert the  pointwise bound   
for $v_1 $, see  (\ref{3.B2}),     obtaining, since $R(x_1) $ is rapidly decreasing, from properties 
(\ref{decay}), 
 \begin{equation*}
 \int_{\R \times \L} R(x_1) v^4_1(x_1) {\rm d} x \le ||v'_1||_{L^2(\R)}^4 L^d \int_\R R(x_1) (1+ |x_1|^2) {\rm d} x_1 = \frac {1} {L^d} C(\b, J) ||v'_1||_{L^2(\R\times \L )}^4.   
\end{equation*}
   Next,  
\begin{eqnarray*}  
\int_{\R \times \L}  R(x_1) w^4  {\rm d}x  &\le& \sup_{x_1 \in \R}  R(x_1) \|w\|^2 _{L^\infty (\R \times \L)}  \int_{\R \times \L}    w^2  {\rm d}x   \nonumber\\
 & \le &
C(\b) c(d,s)  \| w\|_{W^{s,2}}^2   \| w\|^2_{L^2},
\end{eqnarray*}
by Lemma~\ref{A1} in the Appendix.   
 For the other term in (\ref{3.AA1}) we write 
 \begin{equation}\label{O.8}
  \int_{\R \times \L}  v^2 (v_{x_1})^2 {\rm d} x    \le  2\int_{\R \times \L}  v^2 (v'_1)^2 {\rm d} x  + 2\int_{\R \times \L}  v^2 (w_{x_1})^2 {\rm d} x. 
  \end{equation}
 We start  estimating the first term of (\ref{O.8})
 \begin{equation}\label{3.C1} \begin {split}
  \int_{\R \times \L}  v^2 (v_1' )^2 {\rm d} x  & =      \int_{\R \times \L}    [    \a \overline m'  +g  +w ]^2   (v_1')^2 {\rm d} x  \cr & \le  2 \int_{\R \times \L}    \left [ \a^2 (\overline
m')^2 +  [g +w]^2
\right ]  (v_1')^2 {\rm d} x. \end {split}
\end{equation}
Moreover from (\ref{3.B20})    we have
 $$    \int_{\R \times \L}  v'_{1} \overline m''  {\rm d} x =
 \a  \|\overline m''\|^2_{L^2(\R \times \L)}$$ 
and therefore  $|\a|$ can be estimated  as 
 \begin{equation}\label{3.C9}
 |\a|=   \frac { \left | \int_{\R \times \L} v'_{1}  \bar m''  {\rm d} x \right | } {   \|\overline m''\|^2_{L^2(\R\times \L)}}   \le  \frac {\|v'_{1}\|_{L^2(\R \times \L)}   }  {   \|\overline m''\|_{L^2(\R \times \L)}}. 
\end{equation}
We obtain from (\ref{3.C1})
 \begin{equation} \label{3.C2} \begin {split}
  &  \int_{\R \times \L}  v^2 (v_1')^2 {\rm d} x    \le 2  \a^2 \int_{\R \times \L}       (\overline
m')^2  (v_1')^2 {\rm d} x   + 2\int_{\R \times \L}    [g  +w ]^2
  (v_1')^2 {\rm d} x   \\ 
  &\leq  
  2 \frac {\|v'_1\|^2_{L^2(\R \times \L)}   }  {   \|\bar m''\|^2_{L^2(\R \times \L)}}  \|v_1'\|^2_{L^2(\R \times \L)}  \|\overline m'\|_\infty^2    +  4\int_{\R \times \L}    g^2  
  (v_1')^2 {\rm d} x   +  4\int_{\R \times \L}   w^2 
  (v_1')^2 {\rm d} x   \\ 
  & \leq 
   2 \frac {\|v'_1\|^4_{L^2(\R \times \L)}   }  {   \|\bar m''\|^2_{L^2(\R \times \L)}}    \|\overline m'\|_\infty^2   +  4\int_{\R \times \L}    g^2 
  (v_1')^2 {\rm d} x   +  4 \|v_1'\|^2_\infty  \|w \|^2_{L^2}.
   \end {split} \end{equation}
   Since $\|g\|^2_\infty \le  2 \|g\|_{L^2 (\R )} \| g'\|_{L^2 (\R )}$    we have 
    \begin{eqnarray}\label{3.C20} &&\int_{\R \times \L}    g^2 
  (v_1')^2 {\rm d} x  \le   \|g \|^2_{L^\infty  }      \int_{\R \times \L}  (v_1')^2 {\rm d} x  \nonumber\\
   & \le&  2 || g||_{L^2 (\R )} || g'||_{L^2 (\R )} \int_{\R \times \L}  (v_1')^2 {\rm d} x\nonumber\\
    & \le&     2 \lambda  \left (     || g||^2_{L^2 (\R )} || g'||^2_{L^2 (\R )}    \right ) + \frac 2 {\l}   \|v_1'\|^4_{L^2(\R \times \L)},  
    \end{eqnarray} 
     for any $\l>0$. 
    Because of      (\ref{3.BB20})    we have
 $$ ||g ||^2_{L^2(\R)} = ||v_1 - \a \overline m' ||^2_{L^2(\R)}.   $$
 Therefore, by (\ref{3.C9}),  
 \begin{equation}\label{3.BB26}  ||g ||^2_{L^2(\R )} =    ||v_1||^2_{L^2(\R )} + |\a |^2 ||\overline
m'||_{L^2(\R )}^2 
  \le    ||v_1||^2_{L^2(\R )}  + \left ( \frac { ||v_1'||_{L^2(\R\times \L )}}   { ||\overline m''||_{L^2(\R\times \L )}  } ||\overline
m'||_{L^2(\R )} \right )^2.   
\end{equation} 
Taking in account  (\ref{3.C20}),   (\ref{3.BB26})  from  (\ref{3.C2})    one
obtains
 \begin{multline}\label{3.C3}
    \int_{\R \times \L}  v^2 (v_1')^2 {\rm d} x    \le   
    2 \frac {1  }  {   \|m''\|^2_{L^2(\R \times \L)}} \sup_{x_1} (\overline m')^2  ||v'_1 ||^4_{L^2(\R \times \L)}  
      + \frac 2 \l     ||v_1' ||^4_{L^2(\R \times \L)}   
   \\ +  2  \l  \left (   \left [     ||v_1||^2_{L^2(\R )} +\left (\frac { ||v_1'||_{L^2(\R\times \L)}} { ||\overline m''||_{L^2(\R\times \L)} } ||\overline
m'||_{L^2(\R)}   \right)^2 \right ]  \| g'\|^2_{L^2 (\R )}  \right )    + 4 \|v_1'\|^2_\infty  \|w \|^2_{L^2(\R \times \L)}.  
\end{multline}
 Assume  that
 \begin{equation*}
  \| v_1\|^2_{L^2(\R \times \L )} \le k^2, \qquad   \|v_1' \|^2_{L^2(\R\times \L)}  \le k^2 .
\end{equation*}
 Then  
\begin{eqnarray*} &&   \left [     ||v_1||^2_{L^2(\R )} +\left (\frac { ||v_1'||_{L^2(\R\times \L)}} { ||\overline m''||_{L^2(\R\times \L)} } ||\overline
m'||_{L^2(\R)}   \right)^2 \right ]  || g'||^2_{L^2 (\R )}  \nonumber\\ 
&&   \le  
  \frac {k^2  } {L^d} \left [ 1 +     1 \frac {\|\overline
m'\|^2_{L^2(\R)} }  { \|\overline m''\|^2_{L^2(\R)} }  \right  ] \| g'\|^2_{L^2 (\R )} =
  \frac {k^2  } {L^{2d}} \left [ 1+    1     \frac {\|\overline
m'\|^2_{L^2(\R)} }  { \|\overline m''\|^2_{L^2(\R)} }   \right  ]  \| g'\|^2_{L^2 (\R \times \L)} .
\end{eqnarray*}
Take   $\l$  in   (\ref{3.C3})  so that
 \begin{equation}\label{As3}
   2\l \frac {k^2  } {L^{2d}} \left \{  1 + \frac {1} { ||\overline m''||^2_{L^2(\R)} } ||\overline
m'||^2_{L^2(\R)}   \right \}  \le \e_1.   
\end{equation}
We denote such $\l$ by $\l_0= \l_0 (\e_1,k)$. 
Then we have 
 \begin{equation}\label{3.C4}  \begin {split} \int_{\R \times \L}  v^2 (v_1')^2 {\rm d} x   & \le   
    2     ||v_1'||^4_{L^2(\R \times \L)}  \left [  \frac {1  }  {   \|m''\|^2_{L^2(\R \times \L)}}   \| \overline m'\|_\infty^2  
      + \frac {1} {\l_0}   \right ]  
  \cr & +   \e_1 || g'||^2_{L^2 (\R\times \L  )}    +  4 \|v_1'\|^2_\infty  \|w \|^2_{L^2(\R \times \L)}. 
\end {split} \end{equation}
 Next we need to estimate the second term of (\ref{O.8}),  $ \int_{\R \times \L}  v^2 (w_{x_1})^2 {\rm d} x$. 
 Splitting $v=v_1+w$ as in \eqref {P1} we get 
  \begin{equation} \label{O.10}  \begin {split}
   \int_{\R \times \L}  v^2 (w_{x_1})^2 {\rm d} x  & \le  2   \int_{\R \times \L}  w^2 (w_{x_1})^2 {\rm d} x  +  \int_{\R \times \L}  v_1^2 (w_{x_1})^2 {\rm d} x \cr &  \le 
2 \|w_{x_1}\|^2_ \infty \|w\|^2_2 +\int_{\R \times \L}  v_1^2 (w_{x_1})^2 {\rm d} x.   
\end {split} \end{equation}
We split the    last term of (\ref{O.10}) applying the cut-off function   $\phi^2_N$ as it was  done previously in     (\ref{O.9}) but we need 
to end up with an estimate where the  $ \|v_1'\|_{L^2(\R \times \L )}^4$ appears. 
    Denote $h_N(x_1)=  (1+ |x_1|^2)^{\frac 1 2}    (1- \phi^2_N(x_1))$.  We therefore, see    (\ref{O.9}),  have
     \begin{equation} \label{O.90}  \begin {split}   & \int |v_1|^2 |w_{x_1} |^2 {\rm d}x   \\  & \le 
       \|v_1\|^2_\infty \|\phi_N w_{x_1}\|^2_{L^2(\R \times \L)}  + 
   C(\b,J)  \|v_1'\|_{L^2(\R)}^2  \int  h_N(x_1) |w_{x_1} |^2 {\rm d}x. 
 \end {split} \end{equation}
   Note that $h_N (\cdot)$ is smooth and has  support in $[-2N,2N]$. 
 Integrating by part we have 
   $$ \int h_N(x_1) |w_{x_1} |^2 {\rm d}x   =
    - \int   w \left [ h_N' (x_1)w_{x_1} +   h_N(x_1)w_{x_1 x_1} \right ]  {\rm d}x.  $$ 
   By Schwartz inequality we then get 
  \begin{equation*}
    \int h_N(x_1) (w_{x_1} )^2 {\rm d}x  \le \|w\|_2    \left \{  \sup |h_N' (x_1)|  \|w_{x_1}\|_2 +  \sup|h_N (x_1)| \|w_{x_1x_1} \|_2 \right \}. 
   \end{equation*}
   We immediately estimate
   $$  \sup |h_N' (x_1)| \le C,$$
    $$  \sup |h_N (x_1)| \le \sqrt { 1+ 4N^2} \le 3 N. $$
  Summarizing  from    (\ref{O.90})  we obtain  
 \begin{equation*}   
 \begin {split}    
       & \int |v_1|^2 (w_{x_1} )^2 {\rm d}x    \cr &  \le \|v_1\|^2_\infty \|\phi_N w_{x_1}\|^2_2  + C(\b,J)   \frac N {L^d}\|v_1'\|_{L^2 (\R\times \L)}^2   \|w\|_2      \left \{  \|w_{x_1}\|_2  + \|w_{x_1x_1} \|_2  \right \} .   
       \end {split}  
        \end{equation*}
    Since  $ab \le  \frac 12   a^2 + \frac 12  b^2 $   we have
   \begin{equation*} 
    \begin {split}  & \int |v_1|^2 |w_{x_1} |^2 {\rm d}x   \cr&  \le\|v_1\|^2_\infty \|\phi_N w_{x_1}\|^2_{L^2}  +   C(\beta,J) \left [ \frac 12   \|v_1'\|_{L^2(\R \times \L )}^4    +  \frac 12
        \|w\|^2_2  (\frac {N } {L^d} ) ^2  \left [ \|w_{x_1}\|_2 + \|w_{x_1x_1} \|_2 \right ]^2 \right ].   
    \end {split}  
      \end{equation*}
   Summing up all the estimates,  (\ref{O.9}), (\ref{3.C4}) we have 
  \begin{equation*}   \begin {split}     
  \int\sigma(m)\bigl[U(v)\bigr]^2{\rm d}x  &\le     \frac {8} {L^d} C(\b, J) ||v'_1||_{L^2(\R\times \L )}^4  +  8 C(\b,J,d)    \| w\|_{W^{s,2}}^2   \| w\|^2  \\
&+   
  2 \frac {1  }  {   \|\bar m''\|^2_{L^2(\R \times \L)}} \|\overline m'\|_\infty^2  ||v'_1 ||^4_{L^2(\R \times \L)}  
      + \frac 2 {\l_0}     ||v_1'||^4_{L_2 (\R \times \L)}   \\
 & +  8  \e_1 || g'||^2_{L^2 (\R\times \L )}    +  4 \|v_1'\|^2_\infty  \|w \|^2_2  \\
 & +  
 2 \|w_{x_1}\|_\infty \|w\|^2_2   +
  2 \|v_1\|^2_\infty \|\phi_N w_{x_1}\|^2 _2     \\
   & + 
 C(\beta, J)  \left [      \|v_1'\|_{L^2 (\R\times \L)}^4    +   
        \|w\|^2_2     \left(\frac {N} {L^d}\right)^2  \left [ \|w_{x_1}\|_2^2  + \|w_{x_1x_1} \|_2^2  \right ] \right ]. 
         \end {split}    \end{equation*}
Recalling that $ g' = Pv'$, see after \eqref {3.BB20},
 and estimating  $$ \| v_1 \|_\infty \le  \| v \|_\infty   \le c(d,s)\| v \|_{W^{s,2}},  \qquad 
   \| v'\|_\infty + \|w_{x_1} \|_\infty   \le  \| v_{x_1}\|_\infty  \le  c(d,s)\| v \|_{W^{s+1,2}} $$
 we get (\ref{3.1003}).     \qed

\vskip0.5cm
\noindent {\bf Proof  of Proposition~\ref{512} } Writing  $G_\e (v) $ as in  (\ref{O.13}),  applying   Lemma~\ref{510},  Lemma~\ref{510a} and \eqref {O.15}
we have that
 \begin{eqnarray*}  G_\e (v)  & \ge &
 \e  \left [  \frac 14   \s (m_\b)  \g^2(L) \|   \phi_N (w)_{x_1}   \|_2^2 -     \|w\|^2_2 \frac {1} {N^2}   D(\b, \g(L)) \right ]    \nonumber\\  
 &  + &
\e  \s(m_\b)  \g_1\|Pv_1'\|_{L^2 (\R \times \L)}^2   + \e   
\g(L)   \sigma(m_\b)  \sum_{i\ge2} \|w_{x_i}\|^2_2   \nonumber\\
&    -& \frac {1}{ \e} \int\sigma(\overline
m)\bigl[U(v)\bigr]^2{\rm d}x ,
\end{eqnarray*}
 where   $U(v)$ is  defined in  (\ref{3.100}).
Suppose that  
 \begin{equation}\label{3.BB22} 
 ||Pv'_1||_{L^2(\R \times \L)}^2 >\frac 12 ||v'_1||_{L^2(\R \times \L)}^2  
\end{equation}
and
 \begin{equation}\label{NN.1}
    \| v\|_{W^{s+1,2}}^2 \le \e_0, \quad s > \frac {D} {2}.  
     \end{equation}
   From (\ref{3.102}) of  Lemma~\ref{511}
     we obtain 
    \begin{equation*}    \begin {split}
      G_\e (v)  & \ge   \e   \left [  \frac 14   \s (m_\b)  \g^2(L) \|   \phi_N (w)_{x_1}   \|_2^2 -    \|w\|^2_2 \frac {1} {N^2}   D(\b, \g(L) ) \right ]      
\cr &+ 
\e  \s(m_\b)  \g_1\|Pv_1'\|_{L^2 (\R \times \L)}^2   + \e   
\g(L)   \sigma(m_\b)  \sum_{i\ge2} \|w_{x_i}\|^2_2  \cr 
  &- 
 \frac {1}{\e}   \left\{   ||v'_1||^2_{L^2(\R \times \L) }  C(\b,J) \left [ \| v\|_{W^{s,2}}^2    + \||w_{x_1} \|_2 ^2 N   \right ] \right.   
\cr  & + \left .
    \|w\| ^2_2 C(\b,J) \left [ \| v\|_{W^{s,2}}^2  +  \| v \|_{W^{s+1,2}} ^2 \right]  \right. \cr
     &+   \left .   \| v\|_{W^{s,2}}^2    C(\b,J)  \|\phi_N w_{x_1}\|^2_2   \right \}  . 
\end {split}  
\end{equation*}
 To show that $  G_\e (v) \ge 0$ under assumption (\ref{3.BB22})  it  is enough to choose our parameters so that    the following 
 three inequalities are satisfied: 
 \begin{equation}  \label{EE.1}
    \left [ \e  \frac 14   \s (m_\b)  \g^2(L)   -  \frac {1}{\e}    \| v\|_{W^{s,2}}^2    C(\b,J)  \right ]  \|   \phi_N (w)_{x_1}   \|_2^2      \ge 0, 
   \end{equation}
 \begin{equation}\label{EE.2}  \e   
\g(L)   \sigma(m_\b)  \sum_{i\ge2} \|w_{x_i}\|^2  -   \|w\|^2    \left \{    \e            \frac {1} {N^2}   D(\b, \g(L) )    + \frac {1}{\e}  C(\b,J)\left [  \| v\|_{W^{s,2}}^2       +     \| v \|_{W^{s+1,2}} ^2 \right  ] \right \}    \ge  0  
\end{equation}
 \begin{equation}\label{EE.3}
  \e  \s(m_\b)  \g_1\|Pv_1'\|_{L^2 (\R \times \L)}^2  -  \frac {1}{\e}   \|v_1'\|_{L^2 (\R \times \L)}^2  C(\b,J) \left \{  \| v\|_{W^{s,2}}^2 +   \|w_{x_1} \|_2 ^2 N \right \} \ge 0.  
 \end{equation}
  To satisfy   (\ref{EE.1}),    under the assumptions (\ref{NN.1}),      we need  
  \begin{equation}\label{NN.2}
      \g^2(L)  \ge C \frac {\e_0} {\e^2}
    \end{equation}
    for some positive constant $C$.     By  the  Poincar\'e inequality (see the  similar estimate in (\ref{Pa.10})),  we have
 \begin{equation}\label{SS.1}
    \|  \nabla^{\perp}   w\|^2         \ge   \frac {c(d)}  {L^{2 }} \|  w\|^2.
   \end{equation}  
 We can then  satisfy   (\ref{EE.2}),  if 
 \begin{equation*}
   \e   
\g(L)   \sigma(m_\b)  \frac {c(d)}  {L^{2 }}    -        \left \{ \e       \frac {1} {N^2}   D(\b, \g(L) )    + \frac {1}{\e}     C(\b,J)\left [  \| v\|_{W^{s,2}}^2       +     \| v \|_{W^{s+1,2}} ^2 \right  ]     \right ]    \ge  0.  
\end{equation*}
Under the assumptions    (\ref{3.BB22})
it is enough to  require
\begin{equation}\label{EE.4a}       \g(L)     \frac {1}  {L^{2 }}  \ge     \frac {1} {  N^2}  , \qquad 
\g(L)   \frac {1}  {L^{2 }}    \ge  C\frac {\e_0} {\e^2},      \end{equation}  
for some positive constant $C$.
 By (\ref{3.BB22}) to fulfill  (\ref{EE.3}) we need to require
 \begin{equation*}
  \e  \s(m_\b)  \g_1 \frac 1 2   -  \frac {1}{\e}    C(\b,J) \left \{  \| v\|_{W^{s,2}}^2 +    \|w_{x_1} \|_{L^2 (\R \times \L)} ^2 N \right \} \ge 0, 
  \end{equation*}
 which means  
 \begin{equation}\label{FF.1}
   \g_1   \ge   \frac { \e_0 } {\e^2} N.  
\end{equation}
 Choose  the  cut-off $N= \e^{-a} $
with $a = a(L)>0$  so that the first requirement in      (\ref{EE.4a}) holds. 
Choose  then   $\e_0= \e^{2+ a +b}$ with $ b=b(L)$ so that 
   (\ref{NN.2}),  the second condition in  (\ref{EE.4a}) and (\ref{FF.1})  hold.  Let  $r=r(L) \ge   a +b$ and we get \eqref {3.B5},  that is  $   G_\e (v) \ge 0$,  under condition \eqref {3.BB22}.

   Next, suppose    (\ref{3.BB22}) is false, i.e. : 
 \begin{equation}\label{3.B22}
 ||Pv'_1||_{L^2(\R \times \L)}^2 \le \frac 12 ||v'_1||_{L^2(\R \times \L)}^2.  
\end{equation}
Then from (\ref{3.B20}) 
 $$ \|v'_1-g'\|_{L^2(\R)} \le |\a |  \|\overline m''\|_{L^2(\R)} $$
and applying   (\ref{3.B1}) with $\l=  \frac 13$
 and (\ref{3.B22})
one  obtains  
 \begin{equation*}
  |\a |^2  ||\overline m''||^2_{L^2(\R)} \ge  ||v'_1-g'||^2_{L^2(\R)} \ge \frac 13
||v'_1||_{L^2(\R)}^2-\frac 12 ||g'||_{L^2(\R)}^2   \ge \frac {1}{12} ||v'_1||_{L^2(\R)}^2.  
\end{equation*}
Therefore 
 \begin{equation}\label{3.B230}
  ||v'_1||_{L^2(\R)}^2 \le    12  ||\overline m''||^2_{L^2(\R)}  |\a |^2 .
  \end{equation}
Since $ v$  orthogonal to $\overline m' $,  implies  $ v_1$  and $w$ orthogonal to $\overline m' $
we have, see (\ref{3.BB20}), 
 \begin{equation}\label{3.B21}
  \frac {\alpha}  {2 m_\b}  = \frac {1} {2 m_\b } \int_{\R} g
\overline m'{\rm d} x
\le ||g||_{\infty}.  
\end{equation} 
Inequality $||g||^2_\infty \le 2||g||_2    ||g'||_2 $,  
    (\ref{3.B230}) and  (\ref{3.B21}) imply 
 \begin{equation}\label{3.B25}
  ||v'_1||_{L^2 (\R)}  \le C(\b, J) ||g||_{L^2 (\R)}^{\frac 12}  ||g'||_{L^2 (\R)}^{\frac 12}.
  \end{equation}
 Recall, see \eqref {3.B20},  that  $\|Pv'_1\|_2= \|g'\|_2$.   
 We apply estimate      (\ref{3.1003}) with $ \|v_1\|^2_{L^2 (\R \time \L)} \le k^2$,  $\|v_1'\|^2_{L^2 (\R \time \L)} \le k^2$ and for $\e_1>0$.  The actual values of $k$ and $\e_1$ will be chosen later.  
 Thus it is enough to show  
 \begin{equation}  \label{3.B26}   \begin {split}  G_\e (v)    \ge & 
\e  \left [   \frac 14  \s (m_\b)  \g^2(L) \|   \phi_N (w)_{x_1}   \|_2^2 -     \|w\|^2_2 \frac {1} {N^2}   D(\b, \g(L)) \right ] \\     
   + &
\e  \s (m_\b)  \g_1\|Pv_1'\|_{L^2 (\R \times \L)}^2   + \e   
\g(L)   \sigma(m_\b)  \sum_{i\ge2} \|w_{x_i}\|^2_2    \\
      -& 
 \frac {b_0}{\e}   ||v'_1||^4_{L^2(\R\times \L )}  -  \frac {a_0}{\e}    \| w\|^2_{L^2(\R\times \L )} 
      \\   -& \frac {1}{\e}    8  \e_1  \|Pv'_1\|^2_{L^2 (\R \times \L)}   
   -  \frac {1}{\e}  
  2 \|v_1\|^2_{\infty} \|\phi_N w_{x_1} \|^2_2
   \ge 0, 
   \end {split}
   \end{equation} 
  where  we  denoted 
  \begin{equation}\label{GG.5}
   a_0= \left [   C(\b,J,d)   \| w\|_{W^{s,2}}^2     + \left [  4\|v_1'\|^2_\infty +2 \|w_{x_1}\|_\infty \right ]  +  \frac {1}  { 2 m_\b^2}     
         (\frac {N} {L^d})^2  \left [ \|w_{x_1}\|_2+ \|w_{x_1x_1} \|_2\right ]^2 \right ],
              \end{equation}
              and
         $$ b_0= \left [  \frac {1} {L^d} C(\b, J) +2 \frac {1}  {   \|m''\|^2_{L^2(\R \times \L)}} \sup_{x_1} (\overline m')^2  +\frac {2} { \l_0 (\e_1, k)} +   \frac {1 } { 2 m_\b^2}       \right ].  $$
          We estimate
         $$\|v_1'\|_\infty  \le  c(d,s) \| v\|_{W^{s+1,2}} , \qquad   \|w_{x_1x_1} \|_2  \le   \| v\|_{W^{2,2}} \le  \| v\|_{W^{s+1,2}}  , \quad \|w_{x_1}\|_\infty \le  c(d,s) \| v\|_{W^{s+1,2}}.  $$ 
                Therefore, \eqref {GG.5},  by assumptions (\ref{mag1}),  is bounded by 
\begin{equation} \label{FF.2}
      a_0 \le       C  N^2   \| v\|_{W^{s+1,2}}    
   \end{equation} 
   We need to require, see (\ref{3.B26}),   the following 
 three conditions:
\begin{equation} \label{GG.1}
    \left [ \e  \frac 14    \s (m_\b)  \g^2(L) -   \frac {1}{\e}  
  2 \|v_1\|^2_\infty \right ]   \ge 0,  
\end{equation}
  
 \begin{equation}\label{GG.2} \begin {split} 
  & \e   
\g(L)   \sigma(m_\b)  \sum_{i\ge2} \|w_{x_i}\|^2  \\ & -  
       \| w\|^2_{L^2(\R\times \L )}  \left [ \e \frac 14   \frac {1} {N^2}   D(\b, \g(L) )   +   \frac {1}{\e}    a_0  \right ]     
        \ge 0,
  \end {split}     \end{equation} 
  \begin{equation}\label{GG.3} \begin {split}    
  & \e  \s(m_\b)  \g_1\|Pv_1'\|_{L^2 (\R \times \L)}^2   \\ 
  & - 
  \frac {1}{\e}    b_0 ||v'_1||^4_{L^2(\R\times \L )}  
 -   \frac {1}{\e}    8   \e_1  || g'||^2_{L^2 (\R \times \L)}   \ge 0. 
  \end {split}     \end{equation} 
   To satisfy (\ref{GG.1}) taking into account that
   $$   \|v_1\|^2_\infty  \le c(d,s) \|v\|^2_ {W^{s,2}} \le  \e_0$$ 
   we need to require
   \begin{equation*}  \left [ \e   \frac 14    \s (m_\b)  \g^2(L)   -   \frac {1}{\e}  
  \e_0 \right ]   \ge 0.  
 \end{equation*} 
    
  To fulfill (\ref{GG.2}),  taking into account  (\ref{SS.1})  and \eqref {FF.2},  we need to require
 \begin{equation*}
   \e   \g(L)   \sigma(m_\b)  \frac {1} {L^2} c(d)    -  
        \left [  \e \frac 14 \frac {1} {N^2}   D(\b, \g(L))   +    \frac {1}{\e}    C N^2 \e_0 \right ]     
        \ge 0, 
        \end{equation*} 
   therefore 
     \begin{equation*}
      \g(L)     \ge   C(\b,J,d,L) \left [  \frac {1} {  N^2}   + \frac {1} {\e^2}  C N^2 \e_0 \right ].
     \end{equation*} 
      This  forces the choice we have made of $ \e_0$ 
 as $ \e_0= \e^{2+2a+b} $.    To satisfy (\ref{GG.3}), taking into account (\ref{3.B25})    
 we   require 
\begin{equation*}  \e  \s(m_\b)  \g_1\|Pv_1'\|_{L^2 (\R \times \L)}^2   -  \frac {1}{\e}   ||g'||^2_{L^2 (\R \times \L)}   \left \{  b_0 \|g\|^2_{L^2 (\R \times \L)}   +     8   \e_1  \right \}   \ge 0.  
\end{equation*} 

Note that   $$\|Pv_1'\|_{L^2 (\R \times \L)}^2 = \|g'\|_{L^2 (\R \times \L)}^2\ . $$
We now seek to 
bound  $ ||g||^2_{L^2 (\R \times \L)}$  as in (\ref{3.BB26}). 
We would then  require, in terms of order of magnitude, 
 \begin{equation*}
    \g_1   -  b_0  \frac 1 {\e^2}\e_0   - 8  \frac 1 {\e^2} \e_1    \ge 0. 
   \end{equation*} 
  We then choose $\e_1=\e_0$ and $k^2= \e_0$    when applying \eqref {3.AA1} of  Lemma \ref {511}.  Recall  that     $ b_0  \le ( C + \frac 1 {\l_0}) $, and  $\l_0\simeq  \frac {\e_1 } {\k^2}  =    1 $, for the choice done.
   We  therefore get
  \begin{equation*}
      \g_1   -    \frac C {\e^2}\e_0   - 8  \frac 1 {\e^2} \e_0    \ge 0. 
   \end{equation*} 
  Taking $ \e_0  = \e^{2+ r(L)}$ with $r(L)=2a(L) +b(L)$ we get  the thesis. \qed

  \vskip0.5cm
  \noindent 
{\bf Proof of Theorem~\ref{51}:}   Applying   Lemma~\ref{DD1}, see   (\ref{S.S1}),  we have
 \begin{eqnarray*}  
   \II \big (m( t)\big )  &\ge  &    (1-2 \e)  \int_{\R \times \L} \sigma(\overline m)\bigl[(\BB v)_{x_1}  \bigr]^2{\rm d}x    \nonumber\\
  &+& \left [   \e  \left (     \sum_{i\ge1} \int_{\R \times \L}
  \sigma(\overline m)\bigl[(\BB v)_{x_i}  \bigr]^2{\rm d}x \right )-
 \frac {1}{\e}  \int_{\R \times \L} \sigma(\overline m) [U (v)]^2  {\rm d}x \right ] \nonumber\\
  & +&   (1-3 \e)   \sum_{i\ge2} \int_{\R \times \L}
  \sigma(\overline m)\bigl[(\BB v)_{x_i}  \bigr]^2{\rm d}x.  
  \end{eqnarray*}
 Proposition  \ref{512}  then  delivers the thesis. \qed

\section{Bound   on  the dissipation rate of the free energy in terms of the excess free energy}

In this section we establish
a  bound on the rate ${\cal I}(m(t))$ at which the excess free energy 
${\cal F}(m(t)) - {\cal F}(\md)$ is dissipated  in term of 
${\cal F}(m(t)) - {\cal F}(\md)$ itself, working under the hypothesis that
\begin{equation}\label{much}
{\cal I}(m(t)) << \bigl[{\cal F}(m(t)) - {\cal F}(\md)\bigr]\ .
\end{equation}
On the other hand, when  (\ref{much}) is not satisfied,
there is ample dissipation, as explained in the introduction.      
Denote by
\begin{equation}\label{G.8}
  \phi( v(t))= L^d  + \int_{\R \times \L} x_1^2 (\BB_{a(t)} v(t))^2  { \rm  d } x. 
  \end{equation}
      The main
result of this section is the following.

\begin{thm}\label{51b} Let $m(\cdot,t)$ be a solution of
(\ref{1.1}) and set $ m(\cdot,t)= \overline m_{a(t)}(\cdot) + v(\cdot, t)$ where $a(t)$ is chosen so that minimizes $\|m(t)- \overline m_a\|^2_{L^2(\R \times \L)}$.    For any
  $\e>0$ small enough, 
    there is   $\d_1=\d_1 (\e,d,\b,J,L) >0 $ and $\e_1= \e_1(\e,\b,J)$  
   so that   at all time $t$ for which  $ \|v (t)\|_ {W^{s+1,2} } \le \d_1  $,  $|a(t)| \le 1$, where    $s > \frac {D}2$, 
   \begin{equation}\label{may21}
     C(\beta) \|v (t)\|^2    \le  \frac 4 {9 L^2} , 
    \end{equation}  
  
   and 
\begin{equation}\label{anti4}  {\cal I} (m(t)) \le \e_1 \bigl[{\cal F}\big(m(t))-{\cal F}\big(\overline
m     \big)\bigr] 
  \end{equation}
     we have that        
\begin{equation}\label{1.03}
{d\over dt}\bigl[{\cal F}\big(m(t))-{\cal F}\big(\overline
m     \big)\bigr]
\le -  9 (1-\s (m_\b))^2  (1+\e) \frac {\bigl[{\cal
F}\big(m( t))-{\cal F}\big(\overline m  
\big)\bigr]^2}
 {     \phi (t)   } .
\end{equation}
\end{thm} 
 
\medskip

The proof of Theorem \ref{51b}, given at the end of this section, is based on Theorem~\ref{51}  and an application of the following constrained version
of Weyl's uncertainty principle  proved in   Section 2 of  \cite{CCO2}.

\begin{thm}\label{W1}  Let $\psi(x)$ be a function on the real line
such that
\begin{equation}\label{cond1}
\int |\psi'(x)|^2{\rm d}x < \infty\qquad{\rm and}\qquad 
\int |x\psi(x)|^2{\rm d}x < \infty
  \end{equation}
and such that either
\begin{equation}\label{cond2}
\psi(0) = 0
  \end{equation}
or
\begin{equation}\label{cond3}
\int\psi(x){\rm d}x = 0\ .
  \end{equation}
Then
\begin{equation}\label{constrainedup}
\biggl(\int |\psi'(x)|^2{\rm d}x\biggr)
\biggl(\int |x\psi(x)|^2{\rm d}x\biggr) \ge {9\over 4}
\biggl(\int |\psi(x)|^2{\rm d}x\biggr)^2.
  \end{equation}
\end{thm}
\medskip
Notice that under (\ref{cond1}), $\psi$ is 
integrable  and well--defined at $0$,
so (\ref{cond2}) and (\ref{cond3}) make sense. 
\vskip0.5cm
   Recall that    $m (t)= \bar m_{a(t)} +v (t)$ and   $v= v_1+w$ as in   \eqref {P.1}.  We will apply  Theorem \ref {W1}   to $v_1 $,  but the   argument used   in  the one dimensional setting, see
 \cite {CCO2},   does not suffice. Namely we get an extra  term, see the last term of (\ref{anti8}),  due to the multidimensionality  of the problem.

\begin{lm}\label{anti1}   Let  $m (t)= \bar m_{a(t)} +v (t)$, $v= v_1+w$ as in   \eqref {P.1} and $ |a(t)| \le 1$, 
$$ \int_R \left ( x_1  v_1(x_1)\right )^2 dx_1 <\infty,$$ 
 \begin{equation}\label{anti2}
 \int_R v_1(x_1) (\overline m)'_{a(t)} (x_1)    { \rm  d } x_1  =0. 
  \end{equation}
For any $\e>0$, there exists $ \d_1= \d_1 (\e)$  and   $ \e_1= \e_1 (\e,\d_1)$ so that when  $ \|v \|_{W^{s,2}} \le \d_1$ and  
\begin{equation*}
   {\cal I} (m) \le \e_1 \bigl[{\cal F}\big(m)-{\cal F}\big(\overline
m     \big)\bigr], 
  \end{equation*}
 we have
\begin{eqnarray}\label{anti8}
  \int_{\R }  \bigl[(\AA v_1)_{x_1}  \bigr]^2{\rm d}x_1  &\ge&   \frac {(1-\e)^3} {(1+\e)^2}
   \frac 9 4 \frac 
   {\left ( \int_{\R }    (\AA v_1)^2   { \rm  d } x_1 \right )^2  }  { \int_{\R }  x_1^2  (\AA v_1)^2   { \rm  d } x_1 +1}  \nonumber\\
   &-&    \e_1^2  \frac 94     \frac {1} { \e^3 L^{2d} }  \frac {\left (   C(\beta,J)  \left[{\mathcal F}\big(m(t))-{\mathcal F}(\overline m)\right]\right )^2} 
    { \int_{\R }  x_1^2  (\AA v_1)^2   { \rm  d } x_1 +1}. 
\end{eqnarray}
\end{lm}
 
   \noindent{\bf Proof:}   
The proof of the lemma when  $m $ is  antisymmetric  in the $x_1$ variable is   a straightforward application of  (\ref{constrainedup}).  In such a case,    $a(t)=0$ for all $t\ge0$ and    $ (\AA v_1)(0)=0$.  By Theorem \ref {W1} one gets  
 \begin{equation*} 
 \int_{\R }  \bigl[(\AA v_1)_{x_1}  \bigr]^2{\rm d}x_1   \ge 
           \frac 9 4 \frac { \left ( \int_\R  (\AA v_1)^2 {\rm d}x_1\right )^2 }  { \int_\R x_1^2  (\AA v_1)^2  {\rm d}x_1}  
  \end{equation*}
    and (\ref{anti8}) holds for $\e_1=0$ and $\e=0$.    Without this symmetry condition, the proof is more involved. 
 The  argument used in this case in  the one dimensional setting, see
\cite {CCO2},   requires further elaboration:  We get an extra  term, see the last term of (\ref{anti8}),  due to the multidimensionality  of the problem.  

  We introduce the smearing operator
  $${\cal S}v_1(x_1)  = {1\over 2m_\b}\overline m'  \star v_1(x_1)\ .$$
Notice that ${\cal S}$ is a contraction on $L^2(\R)$, and it commutes with
differentiation. Hence,
\begin{equation}\label{DD3}
\| (\AA v_1)_{x_1} \|_{L^2 (\R)} \ge \|{\cal S}\bigl({\cal A}v_1\bigr)_{x_1}\|_{L^2 (\R)} =
\|\bigl({\cal S}{\cal A}v_1\bigr)_{x_1}\|_{L^2 (\R)}\ . 
  \end{equation}
Further, note that
 $${\cal S}(\AA v_1)(a(t)) = 
  \frac{1}{ 2m_\b}    (\overline  m' \star    \AA v_1)(a(t))   = \frac{1}{ 2m_\b} \int_{\R} (\overline m)'_{a(t)}(x_1)  (\AA v_1)(x_1)   d  x_1=0 $$
  by (\ref{anti2}).   Hence the
constrained uncertainty principle applies with the result that, see (\ref{DD3}), 
\begin{equation}\label{falklands}
 \|\  ({\cal A}v_1)_{x_1} \|_{L^2 (\R)}^2  \ge \|\bigl({\cal S}{\cal A}v_1\bigr)_{x_1}\|_{L^2 (\R)}^2 \ge {9\over 4}
{\|{\cal S}{\cal A}v_1\|_{L^2 (\R)}^4\over \|(x_1-a(t)){\cal S}{\cal A}v_1\|_{L^2 (\R)}^2}\ .
  \end{equation}
We now need to remove ${\cal S}$. In the numerator  we have for all $\e > 0$,
$$  
\|{\cal S}{\cal A}v_1\|_{L^2(\R)}^2 = 
 \|{\cal A}v_1 + ({\cal S}{\cal A}v_1 - {\cal A}v_1)\|_{L^2(\R)}^2 \ge 
 (1-\e)\|{\cal A}v_1\|_{L^2(\R)}^2 - {1\over \e}\|({\cal S}{\cal A}v_1 - {\cal A}v_1)\|_{L^2(\R)}^2. $$
Applying Lemma \ref{TT2}    one can show that
\begin{equation*}
\|({\cal S}{\cal A}v_1 - {\cal A}v_1)\|_{L^2(\R)}^2 \le C(\b,J) \|(\AA v_1)_{x_1}\|_{L^2(\R)}^2 .
  \end{equation*}
 Applying Theorem~\ref{51},   
    we have  in particular that 
$$\|(\AA v_1)_{x_1}\|_{L^2(\R)}^2 \le \frac 1 {L^d}C(\beta,J)  \II (m(t) $$
and by (\ref{anti4}) we have 
$$\|(\AA v_1)_{x_1}\|_{L^2(\R)}^2 \le     \e_1  \frac 1 {L^d}C(\beta,J)  \bigl[{\cal F}\big(m(t))-{\cal F}\big(\overline
m     \big)\bigr].  $$ 
Therefore
\begin{equation*}
  \|{\cal S}{\cal A}v_1\|_{L^2(\R)}^2   \ge   (1-\e)\|{\cal A}v_1\|_2^2 - {1\over \e} \e_1  \frac 1 {L^d}C(\beta,J)  \bigl[{\cal F}\big(m(t))-{\cal F}\big(\overline
m     \big)\bigr].  
  \end{equation*} 

Applying inequality (\ref{3.B1}) with $\l=1-\e$ we have  
\begin{equation}\label{problem3}
  \|{\cal S}{\cal A}v_1\|_{L^2(\R)}^4    \ge   (1-\e)^3 \|{\cal A}v_1\|_2^4 - {1\over \e^3}   \left (\e_1  \frac 1 {L^d}C(\beta,J)  \bigl[{\cal F}\big(m(t))-{\cal F}\big(\overline
m     \big)\bigr]\right )^2 .
  \end{equation} 
  To remove ${\cal S}$ from the denominator, write
\begin{equation}\label{borneo}
\int(x_1 - a(t))^2\bigl({\cal S}{\cal A}v_1\bigr)^2{\rm d}x_1\le
(1+\e)\int x_1^2\bigl({\cal S}{\cal A}v_1\bigr)^2{\rm d}x_1 +
\biggl({1+\e\over \e}\biggr) a(t)^2\|{\cal S}{\cal A}v_1\|_{L^2(\R)}^2.
  \end{equation}
By Minkowski's inequality and   
commuting convolution with multiplication by $x_1$, one has
$$\|x_1{\cal S}{\cal A}v_1\|_{L^2(\R)} \le \|{\cal S}x_1{\cal A}v\|_{L^2(\R)} + 
\|\tilde{\cal S}{\cal A}v_1\|_{L^2(\R)}$$
where
$\tilde{\cal S}$ denotes convolution by $(2 m_\b)^{-1}x_1\overline m'(x_1)$. 
Clearly
$\tilde{\cal S}$ is bounded on $L^2$ with norm no greater than 
$(2 m_\b)^{-1}\|x_1\overline m'\|_1$. And since ${\cal S}$ is a contraction on
$L^2$, one has
$$\|x_1{\cal S}{\cal A}v_1\|_2 \le \|x_1{\cal A}v\|_2 + 
(2 m_\b)^{-1}\|x_1\overline m'\|_1\|{\cal A}v_1\|_2\ .$$
Thus, for all $\e>0$,
\begin{equation}\label{sumatra}
\|x_1{\cal S}{\cal A}v_1\|_{L^2(\R)}^2 \le (1+\e)\|x_1{\cal A}v_1\|_{L^2(\R)}^2 +
\biggl({1+\e\over \e}\biggr)(2 m_\b)^{-2}\|x_1\overline m'\|_1^2\|{\cal A}v_1\|_{L^2(\R)}^2\ .
  \end{equation}
Combining (\ref{sumatra}) and (\ref{borneo}),  recalling the hypothesis that
$|a(t)|\le 1$,  and $ {\cal A}$ is a bounded operator one has
\begin{equation}\label{singapore}
\int_\R (x_1 - a(t))^2\bigl({\cal S}{\cal A}v_1 \bigr)^2{\rm d}x_1 \le
(1+\e)^2\|x_1{\cal A}v_1\|_{L^2(\R)}^2 + {1\over 2}\ 
  \end{equation} 
when   $||v_1 ||_{L^2(\R)} $ is sufficiently small.
 Combining  
  (\ref{falklands}), (\ref{problem3}) and (\ref{singapore}), we  obtain  the final result.
\qed 

\vskip0.5cm \noindent 
 \noindent {\bf  Proof of Theorem~\ref{51b}:}  
Since  $  \cal F $ is decreasing along the solution of (\ref{1.1}), see
Theorem  \ref{43}, we have 
\begin{equation*}
   \frac {d {\cal F}} {dt} \big (m( t)\big )=    - {\cal  I }
\big (m(t)\big ), 
  \end{equation*} 
where $ {\cal  I }
\big (m\big )$ is defined in (\ref{8.1}). 
Applying   Theorem \ref{51}, denoting $m = \overline m+v $,     splitting $v$ as in (\ref{P1})   we have    
\begin{eqnarray*} 
 {\cal  I }
\big (\overline m+v \big )  &\ge   &    (1-3 \e)  \left [  \int_{\R \times \L} \sigma(\overline m)\bigl[(\AA v_1)_{x_1}  \bigr]^2{\rm d}x     +     \int_{\R \times \L} \sigma(\overline m)\bigl[(\BB w)_{x_1}  \bigr]^2{\rm d}x  \right ] \nonumber\\
 & + &   (1-3 \e)   \sum_{i\ge2} \int_{\R \times \L}
  \sigma(\overline m)\bigl[(\BB w)_{x_i}  \bigr]^2{\rm d}x.  
  \end{eqnarray*}
  In particular, since $\sigma(\overline m) \ge \sigma(m_\b)$ 
\begin{eqnarray*}
 {\cal  I }
\big (\overline m+v \big ) & \ge   &    (1-3 \e)  \sigma(m_\b)\left [  \int_{\R \times \L}  \bigl[(\AA v_1)_{x_1}  \bigr]^2{\rm d}x        \right ] \nonumber\\
 & + &   (1-3 \e)   \sigma(m_\b) \sum_{i\ge2} \int_{\R \times \L}
 \bigl[(\BB w)_{x_i}  \bigr]^2{\rm d}x. 
 \end{eqnarray*}
 Taking into account that  for each fixed $x_1 \in \R$,   $\int_{ \L}  (\BB w)(x_1, x^\perp)   d x^\perp =0$   we  apply  to the last term the Poincar\'e inequality,   see  \eqref {M.u1},  
obtaining 
\begin{equation}\label{Pa.10}  \begin {split}  \|  \nabla^{\perp} \BB w\|^2  & =   \int _\R   { \rm  d } x_1   \left ( \sum_{i\ge2}   \int_{ \L}  { \rm  d } x^\perp |  (\BB w(x_1, x^\perp))_{x_i} |^2 \right )      \\ &  \ge
\frac {c(d)}  {L^{2}} \int_\R  \int_{ \L} | \BB w(x_1, x^\perp) |^2  dx  = \frac {c(d)}  {L^{2 }} \|\BB w\|^2. 
\end {split} \end{equation}
 Then
\begin{equation*}
  {\cal  I }
\big (\overline m+v \big )  \ge     (1-3 \e) \sigma(m_\b) \left [  \int_{\R \times \L}  \bigl[(\AA v_1)_{x_1}  \bigr]^2{\rm d}x         
   +      \frac {c(d)}  {L^{2 }} \|\BB w\|^2_2 \right ]   .
  \end{equation*}
Applying Lemma~\ref{anti1},  we get
\begin{equation*}
  \begin {split}  {\cal  I }
\big (\overline m+v \big )      &\ge     (1-3 \e)\sigma (m_\b)  \left [  \frac 9 4 \frac {(1- \e)^3} {(1+\e)^2}  \frac {\left ( \int_{\R\times \L}    (\AA v_1)^2 \right )^2 }  { { \int_{\R \times \L }  x_1^2  (\AA v_1)^2   { \rm  d } x +L^d}} +       \frac {c(d)}  {L^{2 }}     \|\BB w\|^2_2   \right ] - \e_1^2 \RR  \\
 &   \ge      (1-3 \e)
    \sigma (m_\b)   \frac 9 4  \frac {(1- \e)^3} {(1+\e)^2} \frac 1    {    \int_{\R\times \L }  x_1^2  (\AA v_1)^2   { \rm  d } x  +L^d} \left [   \left ( \int_{\R\times \L}    (\AA v_1)^2 \right )^2   +     \frac 4 { 9 L^2}      \|\BB w\|^2_2      \right ] - \e_1^2 \RR 
 \end {split}     
 \end{equation*}
     where  
\begin{equation}\label{DO1}
 \RR =  \frac 94     \frac {  C(\beta,J)} { \e^3}  \frac {\left (    \bigl[{\cal F}\big(m(t))-{\cal F}\big(\overline
m     \big)\bigr]\right )^2}     { \int_{\R \times \L }  x_1^2  (\AA v_1)^2   { \rm  d } x +L^d}.   
  \end{equation}
 Our aim is to prove that
\begin{equation*}
\left [  \left ( \int_{\R\times \L}    (\AA v_1)^2 \right )^2    +     \frac 4 { 9 L^2}        \|\BB w\|^2_2   \right ] \ge 
 \left [  \ \bigl[{\cal F}\big(m(t))-{\cal F}\big(\overline
m     \big)\bigr]^2   \right ].  
  \end{equation*}
By orthogonality
\begin{equation*}
  \|\AA v_1\|^2_2   + \|\BB w\|^2_2 =  \|\BB v\|^2_2 .
  \end{equation*}
Then 
$$  \|\BB v\|^4_2 = \left ( \|\AA v_1\|^2_2 + \|\BB w\|^2_2 \right )^2 =  \|\AA v_1\|^4_2 + ( \|\BB w\|^2_2+ 2 \|\AA v_1\|^2_2)  \|\BB w\|^2_2. $$
Suppose   that
\begin{equation}\label{G.10}
 \|\BB w\|^2_2+ 2 \|\AA v_1\|^2_2  \le    \frac 4 { 9 L^2}
  \end{equation}
then 
\begin{equation*}
 \|\AA v_1\|^4_2  +    \frac 4 { 9 L^2}  \|\BB w\|^2_2  \ge   \|\BB v\|^4_2.   
  \end{equation*}
This implies that  \begin{equation*}
   {\cal  I }
\big (\overline m+v \big )  \ge 
  (1-3 \e)\sigma (m_\b)    \frac 9 4 \frac {(1- \e)^3} {(1+\e)^2} 
    \frac {  \| \BB v \|^4_2 } {  \int_{\R  \times \L}  x_1^2  (\AA v_1)^2   { \rm  d } x +L^d }  - \e_1^2 \RR.     
  \end{equation*}
To compare  $ \|\BB v\|^2_2 $ with  $[{\cal F}\big(m(t))-{\cal F}\big(\overline
m     \big)\bigr] $ under assumption \eqref {anti4}  we apply Lemma~\ref{TT3}  stated and proven below  and we obtain
\begin{equation}\label{mn1}
 \|\BB v\|^2_2 \ge    2  \tilde \a (1- 2\e)  [{\cal F}\big(m(t))-{\cal F}\big(\overline
m     \big)\bigr], 
  \end{equation}
 where 
 \begin{equation}\label{atildedef}
 \tilde\alpha = {1\over \beta(1-m_\beta^2)} - 1  = \frac { 1- \sigma(m_\beta) } { \sigma(m_\beta)}.
  \end{equation}
Taking into account (\ref{mn1}) we have
\begin{equation*}
  {\cal  I }
\big (m(t)\big )    \ge (1-\s(m_\b))^2  9   \frac 1 { \s(m_\b)}     \frac {(1- \e)^3} {(1+\e)^2} \frac {  \left [{\cal F}\big(m(t))-{\cal F}\big(\overline
m     \big) \right ]^2     }    { \int_{\R \times \L}  x_1^2  (\AA v_1)^2   { \rm  d } x  +L^d}  - \e_1^2 \RR.
  \end{equation*}
 Recalling the definition  of  $ \RR$, see  (\ref{DO1}),  
 choosing $\e_1$ small enough  so that
 \begin{equation*}
  \frac 94     \frac {\e_1} { L ^d\e^3}      C(\beta,J)    \le  \e 
  \end{equation*}
 we get    
 \begin{equation}\label{G.16c}
    {\cal  I }
\big (m(t)\big )    \ge (1-\s(m_\b))^2  9   \frac 1 { \s(m_\b)}     \frac {(1- 2\e)^3} {(1+\e)^2} \frac {  \left [    {\bigl[{\cal F}\big(m(t))-{\cal F}\big(\overline
m     \big)\bigr]^2 }    \right ] }    { \int_{\R \times \L}  x_1^2  (\AA v_1)^2   { \rm  d } x  +L^d}.  
   \end{equation}
  By  Lemma  \ref {R.1}    we have that
  $$  \int_{\R \times \L}  x_1^2  (\AA v_1)^2   { \rm  d } x    \le \int_{\R\times \L}      x_1^2      ( \BB v)^2  { \rm  d } x   , $$
from (\ref{G.16c}), taking into account (\ref{G.8})  we get (\ref{1.03}). 
Next we verify that   requirement (\ref{G.10}) is indeed satisfied  under  assumptions (\ref{may21}).
Namely $   \|\BB w\|^2_2+ 2 \|\AA v_1\|^2_2 = \|\BB v\|^2_2+  \|\AA v_1\|^2_2$ and 
\begin{multline*}   \|\AA v_1\|^2_2=  \|\BB  v_1\|^2_2  =   \int_{\R}  {\rm d} x_1\left ( \frac 1 {L^d} \int_\L \BB v(x_1, x_1^{\perp}) {\rm d}x_1^{\perp} \right )^2   \le \\  \left ( \frac 1 {L^d} \right)^2   \int_R  {\rm d} x_1\left (  \left ( \int_\L  (\BB v) ^2(x_1, x_1^{\perp}) {\rm d}x_1^{\perp} \right)^{\frac 12} L^{\frac d2} \right )^2 =
\frac 1  {L^d}     \int_{\L \times \R}  (\BB v)^2(x) {\rm d}x. 
\end{multline*}
Then  
$$    \|\BB w\|^2_2+ 2 \|\AA v_1\|^2_2 \le  \|\BB v\|^2_2 [1 + \frac 1 {L^d}] \le   C(\beta) [1 + \frac 1 {L^d}]  \| v\|^2_2,   $$ 
since  $  \|\BB v\|^2_2 \le  C(\beta)  \| v\|^2_2 $.  We get (\ref{G.10}).
 \qed 
\vskip0.5cm \noindent 
 Next we compare  $ \|\BB v\|^2_2 $ with  $[{\cal F}\big(m(t))-{\cal F}\big(\overline
m     \big)\bigr] $.   
\vskip0.5cm \noindent 
\begin{lm}\label{TT3} 
 Take   $v \in W^{1,2} (\R \times \L)$,  
 $\int_{\R \times \L} v \overline m'= 0$ and $m=\bar m +v$.  For any $ \e>0$ there exists  $ \e_1= \e_1 (\e,L, \b, J) $ so that for 
  \begin{equation}\label{anti40}   {\cal I} (m(t)) \le \e_1 \bigl[{\cal F}\big(m(t))-{\cal F}\big(\overline
m     \big)\bigr],   
  \end{equation}
 
 \begin{equation}\label{CE1}
 \|\BB v \|^2_{L^2 (\R \times \L)} \ge   2  \tilde \a (1- 2\e)  [{\cal F} (\bar m +v)-{\cal F} (\overline
m)\bigr] 
  \end{equation}
and 
 \begin{equation}\label{CE2}
 \|\BB v \|^2_{L^2 (\R \times \L)} \le   2  \tilde \a (1+ 2\e)  [{\cal F} (\bar m +v)-{\cal F} (\overline
m)\bigr],  
  \end{equation}
where $\tilde \a$ is defined in (\ref{atildedef}). 
\end{lm}
 
 \noindent{\bf Proof:}     
 We have that
$$    \|\BB v\|^2_2 = \langle \BB v, \BB v\rangle= \tilde \a  \langle v,\BB v\rangle +    \langle (\BB- \tilde \a)  v, \BB v> = \tilde \a   \langle v,\BB v\rangle +    \langle v, (\BB- \tilde \a) ( \BB v)\rangle. $$
  Therefore  $$    \|\BB v\|^2 \ge   \tilde \a  \langle v,\BB v\rangle  - |    \langle   v, (\BB- \tilde \a) ( \BB v) \rangle|. $$
 By Lemma \ref{TT1} 
  we have
   $$ | \langle    v, (\BB- \tilde \a) ( \BB v) \rangle|  \le   \|v\|_{L^2 (\R \times \L)}  \|(\BB- \tilde \a) ( \BB v) \|_{L^2 (\R \times \L)}\le    \|v\|_{L^2 (\R \times \L)}  \|    \|\nabla ( \BB v)\|_{L^2 (\R \times \L)} .  $$ 
  By Theorem~\ref{51}  and assumption (\ref{anti40}) we have 
  $$ C (\beta) \|\nabla ( \BB v)\|^2_{L^2 (\R \times \L)}   \le  \II (v) \le \e_1 [{\cal F}\big(m(t))-{\cal F}\big(\overline
m     \big)\bigr].$$  
Under this condition 
 $$ |\langle   v, (\BB- \tilde \a) ( \BB v) \rangle |  \le   \|v\|_{L^2 (\R \times \L)} [ \e_1 K( \b J){\cal F}\big(m(t))-{\cal F}\big(\overline
m     \big)\bigr]^{\frac 12}. $$
Further, by Lemma \ref{A2} 
$$  \|v\|^2_{L^2 (\R \times \L)} \le \frac 4 {\g (L)} [{\cal F}\big(m(t))-{\cal F}\big(\overline
m     \big)].  $$
Therefore
$$ |\langle    v, (\BB- \tilde \a) ( \BB v) \rangle|  \le  [ \e_1 K(\b J) \frac 4 {\g (L)} ]^{\frac 12}  [{\cal F}\big(m(t))-{\cal F}\big(\overline
m     \big)].   $$
Take  $\e_1$ small enough  so that 
$$ [ \e_1 K(\b J) \frac 4 {\g (L)} ]^{\frac 12} \le \e$$
Then  we get (\ref{CE1}). 
Arguing in a similar way we get (\ref{CE2}). 
\qed

\section{ Moment estimates}

In this section we control  the evolution of
\begin{equation}\label{phidef}
\phi(t) = L^d + \int_{\R \times \L } \sigma(\overline m_{a(t)}) |x_1\left({\cal B}_{a(t)} v (t)\right)|^2{\rm d}x   
\end{equation}
in term  of the free energy functional $\FF$.   As we discussed in the introduction   it is important to  have the right constant  multiplying  the free energy.
 In the next theorem we show two estimates. The first estimate (\ref{phibound1})     does  not quantify the constant  multiplying  the free energy and holds under less  restrictive assumptions.  To show the second estimate, see (\ref{epsn5}),  we need   that the dissipation ${\cal I}(m(t))$ is small compared to the  excess free energy, see  (\ref{epsn5}).
 For proving the  main result    we need both of them. 

\vskip0.5cm \noindent 
\begin{thm}\label{41}  Let $m(\cdot,t)$ be a solution of
(\ref{1.1}).  For any $\e > 0$, $L>0$   there are constants
$\kappa_0(\b,J, \e,L)$,    $\d_0(\b,J,\e,L)$ and $\e_1(\b,J,\e,L)$
such that for all $t$ with   
 $\| v(t)\|_{W^{s,2}} \le  \kappa_0$,   $s>\frac {D} 2$, $||v (t)||_2 <\d_0$ and $|a(t)| \le 1$ there exists a positive constant  $ B= B(\kappa_0,L,d, \beta, J)$        
\begin{equation}\label{phibound1}
{ {\rm d}\over {\rm d}t}\phi (t) \le B
\bigl[{\cal F}(m(t)) - {\cal F}(\overline m)\bigr]\ .
\end{equation}
Further  if 
\begin{equation}\label{epsn5}
{\cal I}(m(t)) \le \e_1\bigl[{\cal F}(m(t)) - {\cal F}(\md)\bigr]
\end{equation}
then  
\begin{equation}\label{phibound}
{ {\rm d}\over {\rm d}t}\phi (t) \le 
(1+\e)4(1-\sigma(m_\b))^2
\bigl[{\cal F}(m(t)) - {\cal F}(\overline m)\bigr]\ .
\end{equation}
\end{thm}

The proof of Theorem  \ref{41}  is based on several intermediate results.
We start  deriving   the 
     full non--linear evolution
for $v$   inserting $m(t) = \overline m_{a(t)}+v (t)$  into (\ref{1.1}).   Taking into account that
 $\overline m$ is a stationary solution of (\ref{1.1}), i.e  
\begin{equation*}
\nabla  \overline m  - 
  {\beta (1-\overline m ^2) (J\star
\nabla \overline m )}   =0,  
\end{equation*}
 we obtain:
\begin{eqnarray}\label{vevoA}
{\partial v\over \partial t} &=& 
  \nabla\cdot \left (   \nabla v- \b (1- \overline m^2) J \star  \nabla v \right )\nonumber\\ 
&+& \beta \nabla\cdot \left ( v(v+2\overline m) J \star  \left (  \nabla v +  \nabla \overline m \right ) \right ) + \dot a(t)\overline m' \nonumber\\
 &=&
 \nabla\cdot \left (   \sigma(\overline m)  \nabla ({\cal B}v) \right ) \nonumber\\
&+&\beta \nabla\cdot \left ( v^2J \star \nabla \overline m \right ) +  \beta \nabla\cdot   \left ( (v(v+2\overline m)J\star \nabla v  \right ) \nonumber\\
&+& \dot a(t)\overline m'.
\end{eqnarray}
Differentiating   (\ref{phidef})   produces    terms involving $\dot a(t)$.   We    estimate these by applying    Theorem~\ref{5}.

\begin{lm} Let $v$ be a solution of (\ref{vevoA}). 
Then for any $\epsilon > 0$   there are   constant $\delta =
\delta(\e, \beta,J)>0$ and $ \kappa= \kappa (\delta, \e, \beta,J)  >0$,  such that for
all $t$ with $\|v(t)\|_2 \le \delta $ and  $\| v(t)\|_{W^{s+1,2}} \le \kappa$ for $s>  \frac {D} 2$ 
\begin{equation}\label{lem1bound} \frac{{\rm d}}{{\rm d}t}\phi  (t) \le 
2\int {\cal B}\left  ( \sigma(\overline m)x_1^2{\cal B}v\right ){\partial v\over \partial t}
{\rm d}x  +
\e  \bigl[{\cal F}(\overline m+v) - {\cal F}(\overline m)\bigr]\ .
\end{equation}
\end{lm}
\bigskip

\noindent{\bf Proof:} 
Since ${\cal B}$ is self adjoint,
\begin{eqnarray}\label{phirate}
\frac{{\rm d}}{{\rm d}t}\int\sigma(\overline m)|x_1{\cal B }v|^2{\rm d}x &=&
2\int {\cal B}\bigl( \sigma(\overline m)x_1^2{\cal B}v\bigr){\partial v\over \partial t}
{\rm d}x\nonumber\\
&+&\dot a(t)2\biggl(\b\int \overline m\overline m'x_1^2|{\cal B}v|^2 {\rm d}x+ 
\int \sigma(\overline m)x_1^2({\cal B}v){2\overline m\overline m'\over \b(1-\overline m^2)^2}v{\rm d}x
\biggr).
\end{eqnarray}
By the exponential decays properties  of $\overline m$, see (\ref{decay}), the boundedness of ${\cal B}$ on $L^2 (\R \times \L)$
and Theorem~\ref{5}, which says that $|\dot a(t)| \le D(\kappa,\beta, J)\|v(t)\|_2$,
one clearly has 
\begin{eqnarray}
\dot a(t)2\biggl(\b\int \overline m\overline m'x_1^2|{\cal B}v|^2 {\rm d}x &+&
\int \sigma(\overline m)x_1^2({\cal B}v){2\overline m\overline m'\over \b(1-\overline m^2)^2}v{\rm d}x
\biggr) \nonumber\\ 
&\le&
|\dot a(t)| C \left [ \|\overline m\overline m'x_1^2\|_\infty \int  |{\cal B}v|^2 {\rm d}x  +\| \sigma(\overline m)x_1^2 { \overline m\overline m'\over \b(1-\overline m^2)^2}\|_\infty 
\int  ({\cal B}v) v{\rm d}x \right ]  \nonumber\\
 &\le &
C (\kappa, \beta, J) \|v(t)\|_2   \bigl[{\cal F}(\overline m+v) - {\cal F}(\overline m)\bigr] 
\nonumber
\end{eqnarray}
where $C$ is a constant depending only on $\b$, $J$  and $\kappa$ that changes from line to line. In the last inequality we applied   Lemma~\ref{TT3} and Lemma~\ref{A2}given  the Appendix.  \qed

\vskip0.5cm \noindent 
We will separately estimate the linear and nonlinear contributions from
(\ref{vevoA}) to (\ref{lem1bound}). Since ${\cal B}\overline m' = 0$, 
the term containing $\dot a$ in (\ref{vevoA}) makes no contribution to (\ref{phirate}).

The basic manipulation, to be done repeatedly in the rest of the proof, is
to commute differentiation
and multiplication by $x_1$ with ${\cal B}$. Therefore we define
\begin{equation*}
g(x_1) =  ( \frac 1 {\sigma (\overline m (x_1))})'= {2\bar m (x_1) \bar m'(x_1)\over \beta (1-\overline m^2(x_1))^2}, \qquad x_1 \in \R
\end{equation*}
and observe that  
\begin{equation}\label{diffcom}
 \nabla \left ({\cal B}w\right) =  e_1 gw + {\cal B} (\nabla w),
 \end{equation}
where $e_1$ is the  $D$  unit vector in the $x_1-$ direction, $ e_1 =(1,0, \dots 0)$.
In (\ref{diffcom}) and in the following  we denote by  ${\cal B} (\nabla w)$
  the $D$ vector with components  ${\cal B} (w_{x_j})$,  $j=1, \dots , D$
 and by $ \| \nabla w  \|^2_2 = \sum_{i=1}^{D}   \|\frac { \partial  w } {\partial x_i} \|^2_2 $. 
Furthermore, define the convolution operator ${\cal C}$ by
\begin{equation*}
{\cal C}w(x) = \int_{\R \times \L} J(y)y_1w(x-y){\rm d}y. 
\end{equation*}
Observe that for any function $w$,
\begin{equation}\label{mulcom}
x_1\bigl({\cal B}w\bigr) = {\cal B}(x_1w) + {\cal C}w
\end{equation}
where $x_1w$ denotes the function with values $x_1 w(x_1, x_1^\perp)$.  
By Young's inequality ${\cal C}$ is bounded on all $L^p$ with operator norm
\begin{equation}\label{cbnd}
\|{\cal C}\| \le \int_{\R \times \L}  |x_1J(x)|{\rm d}x \ .
\end{equation}
We need the following technical lemma.

\begin{lm}\label{63}  For  $w \in L^2(\R \times \L)$,  $w_{x_i} \in L^2(\R \times \L)$, $i=1, \dots,  D$ we have 
\begin{equation}\label{trade1}
\|\sigma(\overline m)x_1 \nabla \bigl({\cal B}w\bigr)\|_2  \le 
{\|\bigl(\sigma (\overline m)x_1\overline m'\bigr)'\|_{L^2(\R \times \L)}\over \|\overline m'\|_{L^2(\R \times \L)} }\|w\|_2   +
\gamma(L)^{-1/2}\|{\cal B}^{1/2}\bigl(\sigma(\overline m)x_1 \nabla \bigl({\cal B}w\bigr) \bigr)\|_2
\end{equation}
where $\gamma (L)$ is the spectral gap (\ref{G.3}) of ${\cal B}$. 
Further    there is
a finite constant $C>0$ depending only on $\beta$ and $J$ such that
whenever $|a(t)| \le 1$,
\begin{equation}\label{trade2}\|J\star (x_1 \nabla w)\|_2 \le C \left (
\|w\|_2+ \|{\cal B}^{1/2}  (\sigma(\overline m )x_1 \nabla \bigl ({\cal
B }w  \bigr)\|_2 \right ). 
\end{equation}
\end{lm}
\bigskip

\noindent{\bf Proof:} Let $P$ denote the orthogonal projection onto
the span of $\overline m'$; i.e., the null space of ${\cal B}$. Then
\begin{eqnarray}
P\bigl(\sigma(\overline m)x_1 \bigl({\cal B}w\bigr)_{x_1}\bigr) &=&
{1\over \|\overline m'\|_2^2}
\langle\overline m',\sigma(\overline m)x_1 \bigl({\cal B}w\bigr)_{x_1}\rangle_{L^2}\overline m'  \nonumber\\
&=& -{1\over \|\overline m'\|_2^2}\langle\bigl(\sigma(\overline m)x_1\overline m'\bigr)',\bigl({\cal B}w\bigr)\rangle_{L^2}
\overline m',\nonumber
\end{eqnarray}
and
$$  P\bigl(\sigma(\overline m)x_1 \bigl({\cal B}w\bigr)_{x_j}\bigr) =
 {1\over \|\overline m'\|_2^2}
\langle\overline m',\sigma(\overline m)x_1 \bigl({\cal B}w\bigr)_{x_j}\rangle_{L^2}\overline m' =0, \quad j \neq 1. $$
Hence, by the Schwarz inequality and  the  fact that   $ {\cal B}$ is bounded, we get 
\begin{equation}\label{schwarzie}
\|P\bigl(\sigma(\overline m)x_1\bigl({\cal B}w\bigr)_{x_1} \bigr)\|_2 \le
{\|\bigl(\sigma(\overline m)x_1\overline m'\bigr)'\|_2\over \|\overline m'\|_2}\|w\|_2\ .
\end{equation}
Next, 
\begin{eqnarray}
\|P^\perp\bigl(\sigma(\overline m)x_1 \nabla \bigl({\cal B}w\bigr) \bigr)\|_2 &=&
\|{\cal B}^{-1/2}{\cal B}^{1/2}P^\perp\bigl(\sigma(\overline m)x_1 \nabla\bigl({\cal B}w\bigr) \bigr)\|_2 \nonumber\\
&\leq&\gamma(L)^{-1/2}\|{\cal B}^{1/2}P^\perp\bigl(\sigma(\overline m)x_1 \nabla \bigl({\cal B}w\bigr) \bigr)\|_2 \nonumber\\
&=&\gamma(L)^{-1/2}\|P^\perp{\cal B}^{1/2}
\bigl(\sigma(\overline m)x_1 \nabla \bigl({\cal B}w\bigr) \bigr)\|_2 \nonumber\\
&\leq&\gamma (L)^{-1/2}\|{\cal B}^{1/2}
\bigl(\sigma(\overline m)x_1\nabla \bigl({\cal B}w\bigr)\bigr)\|_2\ .\nonumber
\end{eqnarray}
Hence, the Minkowski inequality and (\ref{schwarzie}) yield (\ref{trade1}).
To  prove (\ref{trade2})  we 
define the operator ${\cal D}$ by
\begin{equation}\label{const1}
{\cal D}w = {1\over \b(1-m_\b^2)}w - J\star w\ .
\end{equation}
 Fourier transforming, one sees
that ${\cal D}$ is bounded with a bounded inverse since
$\b(1-m_\b^2) < 1$. 
Note that
\begin{equation}\label{TT.1}
 {\cal D}w =  \BB w - \tilde g w,  
 \end{equation}
where 
\begin{equation}\label{mars2}
\tilde g (x_1) = {1\over \b(1 - \overline m^2(x_1))} - {1\over \b(1 - m_\b^2)}, \quad x_1 \in \R. \end{equation}
Also, ${\cal D}$ commutes with convolution by $J$, and
\begin{equation}\label{mars1}
x_1{\cal D}w = {\cal D}(x_1w) + {\cal C}w,   
\end{equation}
 as with ${\cal B}$ in (\ref{mulcom}). Hence,
\begin{multline}
J\star x_1 \nabla w = {\cal D}^{-1}J\star\bigl({\cal D} (x_1 \nabla w) \bigr) = {\cal
D}^{-1}J\star\bigl(x_1{\cal D} (\nabla w)  - {\cal C}(\nabla w)   \bigr) =\nonumber\\-{\cal D}^{-1}{\cal
C}\bigl(\nabla J \star w\bigr) + {\cal D}^{-1}J\star\bigl(x_1{\cal B}\nabla w - x_1\tilde g \nabla w\bigr),  
\end{multline}
  where we have
used   that convolution with $J$ commutes with
${\cal C}$ and that $J\star \nabla w = \nabla J \star w$. Next, applying  (\ref{diffcom}),
$${\cal D}^{-1}J\star\bigl(x_1 {\cal B}(\nabla w)  - x_1\tilde g \nabla w \bigr) = {\cal
D}^{-1}J\star\bigl(x_1 \nabla \bigl({\cal B}w\bigr) \bigr)  -  {\cal
D}^{-1}J\star (x_1g w e_1+ x_1\tilde g  \nabla w\bigr).$$ 
We have that
$$ \| {\cal
D}^{-1}J\star (x_1g w e_1+ x_1 \tilde g  \nabla w \bigr)\|_2   \le C(\b,J) \|w\|_2,$$
where  
we  used the rapidly decay of $g$ and  $\tilde g$ and that  $J\star (x_1 \tilde g   \nabla w)  = \nabla J \star   (x_1 \tilde  g w)  +  J \star   \nabla (x_1 \tilde g) w $.  Thus,
\begin{eqnarray}
\|J\star(x_1 \nabla w)\|_2 
&\leq&\|{\cal D}^{-1}{\cal C}\left(\nabla J \star
w\right)\|_2 \nonumber\\ 
&+&
\|{\cal D}^{-1}J\star\bigl(x_1 \nabla \bigl({\cal B}w\bigr) \bigr)\|_2   \nonumber\\
 & +& \| {\cal
D}^{-1}J\star (x_1g w e_1+ x_1 \tilde g  \nabla w \bigr) \|_2  
\nonumber\\
&\leq&C(\b,J)\bigl(\|w\|_2 e_1+ \|\sigma(\overline m)x_1 \nabla \bigl({\cal B}w\bigr)\|_2\ . \nonumber
\end{eqnarray} 
 Now application of (\ref{trade1})
yields
(\ref{trade2}). \qed
\medskip

\noindent Next we estimate the nonlinear contribution from (\ref{vevoA}) to (\ref{lem1bound}).
\medskip
\begin{lm}\label{64}  Let $v$ be a solution of (\ref{vevoA}). 
Then for any $\epsilon > 0$ there are  constants $\d =\delta(\beta,J,\e,L)>0$ 
and $\kappa = \kappa(\beta,J,\e,L )>0$ 
such that for
all $t$ with $\|v(t)\|_2 \le \delta$,  $\| v(t)\|_{W^{s,2}} \le  \kappa $, and $|a(t)| \le 1$,

\begin{eqnarray*}
 \frac{{\rm d}}{{\rm d}t}\phi (t) &\leq&
2  \int {\cal B}\bigl( \sigma(\overline m)x_1^2{\cal B}v\bigr) \nabla\cdot \left (   \sigma(\overline m)  \nabla ({\cal B}v) \right )   
{\rm d}x  \nonumber\\
&+& \e   \bigl[{\cal F}(\overline m+v) - {\cal F}(\overline m)\bigr]
+ \e\|{\cal B}^{1/2}\bigl(\sigma(\overline m)x_1 \nabla \bigl({\cal B}v\bigr)\|_2^2.
\end{eqnarray*}
\end{lm}

\noindent{\bf Proof:} We separately estimate the contribution of the
two nonlinear terms in (\ref{vevo}) to  (\ref{lem1bound}), beginning with the more
difficult of the two:
\begin{equation}\label{badone}
2\int {\cal B}\bigl( \sigma(\overline m)x_1^2{\cal B}v\bigr)\beta \nabla   \cdot \bigl(v(v+2\overline m)J \star \nabla v \bigr)
{\rm d}x.
\end{equation}
Now integrating by parts and applying (\ref{diffcom}) to (\ref{badone})
yields
$$  \int {\cal B}\bigl( \sigma(\overline m)x_1^2{\cal B}v\bigr)\beta \nabla   \cdot \bigl(v(v+2\overline m)J \star \nabla v \bigr)
{\rm d}x = A_1 +A_2 $$
with 
\begin{eqnarray*} A_1&=& -2\int g   \bigl( \sigma(\overline m)x_1^2{\cal B}v\bigr) e_1
\beta\bigl(v(v+2\overline m)  \nabla J \star  v \bigr){\rm d}x \nonumber\\ 
&-& 2\int {\cal B}\bigl( \sigma(\overline m)'x_1^2{\cal B}v\bigr)
\beta\bigl(v(v+2\overline m)e_1 \nabla J *v\bigr){\rm d}x,  \nonumber
\end{eqnarray*} 
\begin{eqnarray*}  A_2 &=&
-4\int {\cal B}\bigl( \sigma(\overline m)x_1 e_1{\cal B}v\bigr)
\beta\bigl(v(v+2\overline m)J*\nabla v\bigr){\rm d}x\nonumber\\
&-&2\int {\cal B}\bigl( \sigma(\overline m)x_1^2 \nabla \left (  {\cal B}v \right ) \bigr)
\beta\bigl(v(v+2\overline m)J*\nabla v\bigr){\rm d}x
\end{eqnarray*}
where the first term in $A_1$   comes from the first term of  (\ref{diffcom}), and the remaining term of $A_1$ together with the term $A_2$ 
  from differentiating the product $\sigma(\overline m)x_1^2{\cal B}v$.
We have also used the fact that $J\star  \nabla v  =  \nabla J \star v$.
Because of (\ref{decay}),
 by Lemma~\ref{A2} in the appendix and the inequality 
$\|v\|_\infty  \le c(d,s) \|v \|_{W^{s,2}} $,  for $s> \frac {D} 2 $
\begin{equation}\label{temp}
 |A_1| \le C\|v\|_\infty\|v\|_2^2   \le  C c(d,s) \|v \|_{W^{s,2}} \|v\|_2^2  \le   C c(d,s) \|v \|_{W^{s,2}}  {1 \over   \g (L)}  \bigl[{\cal F}(\overline m+v) - {\cal F}(\overline m)\bigr] 
 \end{equation}
where $C$ is a constant depending only on $\beta$ and $J$. 
Then  for any $\e>0$   there are constants
 $\delta>$ 
and $\kappa>0$ 
such that for
all $t$ 
the quantity in (\ref{temp}) is no greater than
\begin{equation}\label{temp2}
  |A_1| \le    {\e\over 3}  \bigl[{\cal F}(\overline m+v) - {\cal F}(\overline m)\bigr]\ .
  \end{equation}
To estimate $A_2$  we need to commute an $x_1$ past 
${\cal B}$. Applying (\ref{mulcom}), these become
\begin{eqnarray}A_2 &=&
 4\int {\cal C}\bigl( \sigma(\overline m){\cal B}v\bigr) e_1
\beta \bigl(v(v+2\overline m)J* \nabla v  \bigr){\rm d}x
+2\int {\cal C}\bigl( \sigma(\overline m)x_1 \nabla \bigl({\cal B}v\bigr)\bigr)
\beta\bigl(v(v+2\overline m)J*\nabla v \bigr){\rm d}x\nonumber\\
&-&4\int {\cal B}\bigl( \sigma(\overline m){\cal B}v\bigr) e_1
\beta\bigl(v(v+2\overline m)x_1 J*\nabla v \bigr){\rm d}x
-2\int {\cal B}\bigl( \sigma(\overline m)x_1 \nabla \bigl({\cal B}v\bigr) \bigr)
\beta\bigl(v(v+2\overline m)x_1 J*\nabla v \bigr){\rm d}x. \nonumber
\end{eqnarray}
Now, it is exactly the convolution by $J$ in ${\cal B}$ that doesn't commute
with multiplication by $x_1$ so that
$$x_1J\star w = J\star (x_1 w) + {\cal C}w$$
so that the integrals above can be partially rewritten as
\begin{eqnarray} A_2 &=&
 4\int {\cal C}\bigl( \sigma(\overline m){\cal B}v\bigr)  
\beta \bigl(v(v+2\overline m) e_1 \nabla J*   v  \bigr){\rm d}x
+2\int {\cal C}\bigl( \sigma(\overline m)x_1 \nabla \bigl({\cal B}v\bigr)\bigr)
\beta\bigl(v(v+2\overline m)J*\nabla v \bigr){\rm d}x\nonumber\\
&-&4\int {\cal B}\bigl( \sigma(\overline m){\cal B}v\bigr)
\beta\bigl(v(v+2\overline m) e_1 {\cal C} (\nabla v)\bigr){\rm d}x
-2\int {\cal B}\bigl( \sigma(\overline m)x_1 \nabla \bigl({\cal B}v\bigr)\bigr)
\beta\bigl(v(v+2\overline m){\cal C} (\nabla v\bigr){\rm d}x\nonumber\\
&-&4\int {\cal B}\bigl( \sigma(\overline m){\cal B}v\bigr)
\beta\bigl(v(v+2\overline m) e_1 J*(x_1 \nabla v )\bigr){\rm d}x\nonumber\\
&-&2\int {\cal B}\bigl( \sigma(\overline m)x_1 \nabla \bigl({\cal B}v\bigr)\bigr)
\beta\bigl(v(v+2\overline m)J*(x_1 \nabla v )\bigr){\rm d}x.\nonumber
\end{eqnarray}
Clearly there is a constant $C$ depending only on $\beta$ and $J$ so that
$$\|{\cal C} (\nabla v) \|_2 \le C\|v\|_2,$$
and hence the four terms containing ${\cal C}$ may be estimated, as in (\ref{temp})   by 
$$ \|v\|_\infty\|v\|_2^2  \le   c(d,s) \|v \|_{W^{s,2}}  {1 \over   \g (L)}  \bigl[{\cal F}(\overline m+v) - {\cal F}(\overline m)\bigr]. $$
Hence   there are  constants $\kappa$ and $ \delta$ so that 
\begin{equation}\label{temp3}
  C c(d,s) \|v \|_{W^{s,2}}  {1 \over   \g (L)}  \bigl[{\cal F}(\overline m+v) - {\cal F}(\overline m)\bigr] \le  {\e\over 3  }  \bigl[{\cal F}(\overline m+v) - {\cal F}(\overline m)\bigr]
\end{equation}
for all $t$ with $\|v(t)\|_2 \le \delta$, $\|v \|_{W^{s,2}}  \le \kappa$ and $|a(t)| \le 1$. 

Next, by the Schwarz inequality, and then (\ref{trade2}) of Lemma \ref {63},
\begin{eqnarray}\label{temp4}
-4 \int {\cal B}\bigl( \sigma(\overline m){\cal B}v\bigr)
\beta\bigl(v(v+2\overline m) e_1 J*(x_1 \nabla v )\bigr){\rm d}x   &\leq&
C\|v\|_\infty\|v\|_2\|J*(x_1 \nabla v )\|_2  \nonumber\\ 
&\leq& C\|v\|_\infty\|v\|_2\nonumber\\
&\times& \left [  \|v\|_2 + \|{\cal B}^{1/2}\sigma(\overline m)x_1 \nabla \bigl({\cal
B}v\bigr)\|_2 \right ]
\end{eqnarray}
and
\begin{eqnarray}\label{temp5}
-2 \int {\cal B}\bigl( \sigma(\overline m)x_1 \nabla \bigl({\cal B}v\bigr)\bigr)
\beta\bigl(v(v+2\overline m)J*(x_1 \nabla v )\bigr){\rm d}x &\leq&
C\|v\|_\infty\|{\cal B}^{1/2}\sigma(\overline m)x_1\nabla \bigl({\cal B}v\bigr)\|_2
\|J*(x_1 \nabla v )\|_2\nonumber\\
&\leq& C\|v\|_\infty\|{\cal B}^{1/2}\sigma(\overline m)x_1\nabla \bigl({\cal B}v\bigr) \|_2\nonumber\\
&\times&
\left [  \|v\|_2 + \|{\cal B}^{1/2}\sigma(\overline m)x_1\nabla \bigl({\cal B}v\bigr)\|_2\right ].
\end{eqnarray}
 Hence the sum of the two terms in (\ref{temp4}) and (\ref{temp5}) is no greater than
$$ C\|v\|_\infty\ \left [ \|v\|_2^2 +\|{\cal B}^{1/2}\sigma(\overline m)x_1\nabla \bigl({\cal B}v\bigr)\|_2 ^2\right ] $$
and now decreasing $\delta$ and $\kappa$ as necessary, we obtain
as before from $\|v\|_\infty  \le c(d,s) \|v \|_{W^{s,2}}  $   and
Lemma~\ref{A2} in the Appendix  that this is no greater than
\begin{equation}\label{temp7}
 |A_2| \le {\e\over 3 } \bigl[{\cal F}(\overline m+v) - {\cal F}(\overline m)\bigr]
+ \e\|{\cal B}^{1/2}\sigma(\overline m)x_1 \nabla \bigl({\cal B}v\bigr)\|_2^2 
\end{equation}
for all $t$ with $\|v(t)\|_2 \le \delta$,  $\|v \|_{W^{s,2}}  \le \kappa$ and $|a(t)| \le 1$.
Thus the estimate on (\ref{badone}) follows from
(\ref{temp2}), (\ref{temp3}) and (\ref{temp7}).

It remains to estimate the contributions to   (\ref{lem1bound}) from the other of the two
non-linear terms in (\ref{vevoA}), namely 
\begin{equation*}
 -2\int  \nabla \left ({\cal B}\left ( \sigma(\overline m)x_1^2{\cal B}v \right )\right ) 
\bigl( v^2 J\star \nabla \overline m\bigr){\rm d}x.
\end{equation*}
Proceeding as above,
though with with much less effort, one obtains  that  this term  is bounded by
\begin{equation*}
\|v\|_\infty\bigl(C \|v\|_2^2 + 
\|{\cal B}^{1/2}\sigma(\overline m)x_1\nabla \bigl({\cal B}v\bigr)\|_2 ^2\bigr)
\end{equation*}
where the extra factor of $\|v\|_\infty$ comes
from the nonlinearity.  Using once more the inequality
$\|v\|_\infty  \le c(d,s) \|v \|_{W^{s,2}} $, one sees that for $\delta$ sufficiently small,
one can combine the above estimates, once more using Lemma~\ref{A2} in the Appendix,   to obtain 
the proof of the lemma. \qed

\begin{thm}\label{(65)} Let $v$ be a solution of (\ref{vevoA}).  For any
$\epsilon > 0$    there are  constants $\d = \delta(\beta,J,\e,L)>0$  and
$\kappa = \kappa(\beta,J,\e,L )>0$  such that for all $t$ with $\|v(t)\|_2 \le
\delta$, $\|v \|_{W^{s,2}} \le \kappa$, for $s>\frac {D} 2$  and $|a(t)| \le 1$,

\begin{eqnarray}\label{lem6bound}
\frac{{\rm d}}{{\rm d}t}\phi (t) &\leq&
-4\int {\cal B}\bigl( \sigma(\overline m){\cal B}v\bigr)
x_1 \sigma(\overline m)  e_1 \cdot  \nabla ({\cal B}v) {\rm d}x 
-2\|{\cal B}^{1/2} \left (\sigma(\overline m)x_1 \nabla \bigl({\cal B}v\bigr) \right )\|_2^2 \nonumber\\
&+&I_1 + I_2 + I_3 + I_4 \nonumber\\
&+& \e  \bigl[{\cal F}(\overline m+v) - {\cal F}(\overline m)\bigr] +
\e\|{\cal B}^{1/2}\left ( \sigma(\overline m)x_1 \nabla \bigl({\cal B}v\bigr)\right) \|_2^2 
\end{eqnarray}
where
\begin{eqnarray*}
I_1 &=&    -2 \int g (x_1) \bigl( \sigma(\overline m)x_1^2{\cal B}v\bigr)
e_1 \sigma(\overline m)   \nabla ({\cal B}v)  {\rm d}x,   \nonumber\\
I_2 &=&  - 2   \int {\cal B} \left (     \sigma(\overline m)'   x_1^2  {\cal B}v \right )   e_1
 \sigma(\overline m) \nabla ({\cal B}v) {\rm d}x,\nonumber\\
 I_3 &=&
-4 \int {\cal C}\bigl( \sigma(\overline m){\cal B}v\bigr) e_1
\sigma(\overline m)  \nabla ({\cal B}v)  {\rm d}x,
 \nonumber\\
I_4 &=& -2\int {\cal C}\bigl( \sigma(\overline m)x_1 \nabla ({\cal B}v) \bigr)
\sigma(\overline m)  \nabla ({\cal B}v)  {\rm d}x. \nonumber\\
\end{eqnarray*}
\end{thm}

\noindent{\bf Proof:} Denote by  $A$ 
\begin{equation*}
A = 2 \int {\cal B}\bigl( \sigma(\overline m)x_1^2{\cal B}v\bigr) \nabla\cdot \left (   \sigma(\overline m)  \nabla ({\cal B}v) \right )  
 = - 2 \int \nabla \left (   {\cal B}\bigl( \sigma(\overline m)x_1^2{\cal B}v\bigr)\right )     \sigma(\overline m)  \nabla ({\cal B}v). 
 \end{equation*}
 By Lemma~\ref{64} the only term to take care to get (\ref{lem6bound}) is $A$. 
Now applying (\ref{diffcom})   yields
\begin{equation}\label{intss}
A = 
-2 \int g (x_1) \bigl( \sigma(\overline m)x_1^2{\cal B}v\bigr) e_1
\sigma(\overline m)  \nabla ({\cal B}v)  {\rm d}x  
-2  \int {\cal B} \left ( \nabla \left (   \sigma(\overline m)x_1^2{\cal B}v \right )\right ) \sigma(\overline m)  \nabla ({\cal B}v)  {\rm d}x.    
\end{equation}
Further differentiating the product $\sigma(\overline m)x_1^2{\cal B}v$ we have 
\begin{eqnarray}
 -2 \int {\cal B} \left ( \nabla \left (   \sigma(\overline m)x_1^2{\cal B}v \right )\right ) \sigma(\overline m)  \nabla ({\cal B}v)  {\rm d}x   &=&
-2  \int {\cal B} \left (    ( \sigma(\overline m) )'    x_1^2  {\cal B}v \right )  
e_1  \sigma(\overline m)(\nabla ({\cal B}v)   {\rm d}x  
\nonumber\\   
&-&4\int {\cal B}\bigl( \sigma(\overline m)x_1  {\cal B}v\bigr)
e_1\sigma(\overline m)  \nabla ({\cal B}v){\rm d}x\nonumber\\
&-&2\int {\cal B} \left (  \sigma(\overline m)x_1^2  \nabla ({\cal B}v) \right ) 
\bigl( \sigma(\overline m)(\nabla ({\cal B}v) \bigr){\rm d}x. 
\nonumber
\end{eqnarray}
 
 Denote  
\begin{equation*}
I_1=  -2 \int g (x_1) \bigl( \sigma(\overline m)x_1^2{\cal B}v\bigr)
e_1 \sigma(\overline m)  \nabla ({\cal B}v)  {\rm d}x, 
\end{equation*}
\begin{equation*}
 I_2=  - 2   \int {\cal B} \left (    ( \sigma(\overline m) )'    x_1^2  {\cal B}v \right )  
e_1 \sigma(\overline m) \nabla ({\cal B}v)  {\rm d}x. 
\end{equation*}
 We obtain, see    
(\ref{intss}),   
\begin{eqnarray*}
A &=& I_1 + I_2 \nonumber\\
&-& 4\int {\cal B}\bigl( \sigma(\overline m)x_1  {\cal B}v\bigr) e_1
\sigma(\overline m)  \nabla ({\cal B}v){\rm d}x\nonumber\\
&-&2\int {\cal B} \left (  \sigma(\overline m)x_1^2  \nabla ({\cal B}v) \right ) 
\bigl( \sigma(\overline m)(\nabla ({\cal B}v) \bigr){\rm d}x\ .
\end{eqnarray*}
Next, to exploit the positivity of ${\cal B}$, we need
to distribute the factors of $x_1$ symmetrically in the last integral.
To do this, apply (\ref{mulcom}) to account for commuting multiplication by $x_1$
with ${\cal B}$. We also do this in the other integral, so that
the same function $\sigma(\overline m)x_1\nabla ({\cal B}v)$ is produced there as well.
The result is
\begin{eqnarray}A &=& I_1 + I_2 \cr
&-&4\int {\cal C}\bigl( \sigma(\overline m){\cal B}v\bigr) e_1
\sigma(\overline m)  \nabla ({\cal B}v)  {\rm d}x\nonumber\\
&-&2\int {\cal C}\bigl( \sigma(\overline m)x_1 \nabla ({\cal B}v) \bigr)
\sigma(\overline m)  \nabla ({\cal B}v)  {\rm d}x\nonumber\\
&-&4\int {\cal B}\bigl( \sigma(\overline m){\cal B}v\bigr)
x_1 e_1\sigma(\overline m)  \nabla ({\cal B}v) {\rm d}x\nonumber\\
&-&2\int {\cal B}\bigl( \sigma(\overline m)x_1\nabla ({\cal B}v) \bigr)
\sigma(\overline m) x_1 \nabla ({\cal B}v){\rm d}x\ .\nonumber
\end{eqnarray}
Now denote the first two terms after $I_1$ and $I_2$; i.e., those containing ${\cal C}$,
by $I_3$ and $I_4$ respectively. Then, by Lemma 4.4, the result is proved. \qed
\bigskip

 \vskip1.cm
 \noindent{\bf Proof of Theorem  \ref{41}: }  
The starting point for  proving both (\ref{phibound1}) and (\ref{phibound}) is   Theorem~\ref{(65)}.  
 We start proving (\ref{phibound1}).
 The first two terms in (\ref{lem6bound}) are the key
to the analysis. They correspond to the two terms produced in (\ref{heatex})
when similar estimates were performed on the heat equation 
as an illustration of the method.
To see this more easily, introduce the
following notations:
\begin{equation}\label{fdef}
f =   e_1 \sigma(\overline m){\cal B}v 
\end{equation}
 and 
\begin{equation}\label{hdef} 
h = \sigma(\overline m)x_1 \nabla ({\cal B}v)\ ,
\end{equation}
$$  \langle f, h\rangle = \sum_{i=1}^{D} \int_{\R \times \L} f_i(x) h_i(x) dx .$$ 
Notice that $f_i=0$ for all $i\ge 2$.
These  first two terms in (\ref{lem6bound}) can be written  as following:
\begin{equation}\label{mag3} \begin {split}
-4\langle f,{\cal B }h\rangle -2\langle h,{\cal B}h\rangle & = 
-\langle h+2f,{\cal B}(h+2f)\rangle - \langle h,{\cal B}h\rangle +
4\langle f,{\cal B}f\rangle  \cr &\le
-\langle h,{\cal B}h\rangle  + 4\langle f,{\cal B}f\rangle. 
\end {split} \end{equation}  
 The next step is to estimate each of the $I_j$ appearing in   (\ref{lem6bound}) in terms of
$\|v\|_2^2$, using the negative
term in (\ref{mag3}) to absorb contributions from $\nabla v$.

First, using the Schwarz inequality, and then the arithmetic--geometric mean
inequality,
$$
I_1 \le 2\|g \sigma(\overline m)x_1^2{\cal B}v\|_2
\|\sigma(\overline m) \nabla ({\cal B}v)  \|_2  \le 
\lambda \|g \sigma(\overline m)x_1^2{\cal B}v\|^2_2 +
{1\over \lambda}\|\sigma(\overline m)\nabla ({\cal B}v)\|^2_2$$
for any $\lambda > 0$. Now choose $\lambda$ so large that the estimate
(\ref{trade1}) of Lemma ~\ref{63}  gives
$${1\over \lambda}\|\sigma(\overline m)\nabla ({\cal B}v)\|_2^2 \le 
{1 \over 4}\left (\|v\|_2^2 + \langle h,{\cal B}h\rangle \right )$$
where $h$ is given in (\ref{hdef}).
The choice of $\l$ depends  also on  $L$. 
One obtains
a constant $C$ depending
  on $\beta$,  $J$ and $L$  such  that
\begin{equation}\label{I1bnd}
I_1 \le 
{1\over 4}\langle h,{\cal B} h\rangle + C\|v\|_2^2\ .
\end{equation}
It is easier to deal with $I_2$. Schwarz and (\ref{decay}) suffice to establish
that there is a constant $C$ depending
only on $\beta$,  and certain finite moments of $\overline m'$ so that
\begin{equation}\label{I2bnd}
I_2 \le 
{1\over 4}\langle h,{\cal B} h\rangle + C\|v\|_2^2\ .
\end{equation}
To bound $I_3$, we will integrate by parts. Note that using (\ref{fdef})
$$I_3 = -4\int\bigl({\cal C}f\bigr)\sigma(\overline m)
 \nabla ({\cal B}v) {\rm d}x
= 4\int \nabla ({\cal C}f\bigr) \sigma(\overline m)
\bigl({\cal B}v\bigr){\rm d}x + 4\int\bigl({\cal C}f\bigr)\sigma(\overline m)'
\bigl({\cal B}v\bigr){\rm d}x.$$
Using this, (\ref{cbnd}) and the rapid decay of $\sigma(\overline m)'$ coming
from (\ref{decay}), there is clearly a constant $C$ depending only on $\beta$ and $J$
so that
\begin{equation}\label{I3bnd}
I_3 \le 
C\|\sigma(\overline m){\cal B}v\|_2^2\ .
\end{equation}
Finally, to bound $I_4$, we use (\ref{hdef}) and again integrate by parts:
$$I_4 = -2\int\bigl({\cal C}h\bigr)\sigma(\overline m)
\nabla \bigl({\cal B}v\bigr) {\rm d}x =
2\int \nabla\bigl({\cal C}h\bigr) \sigma(\overline m)
\bigl({\cal B}v\bigr){\rm d}x + 2\int\bigl({\cal C}h\bigr)\sigma(\overline m)'
\bigl({\cal B}v\bigr){\rm d}x.$$
Now proceeding as with $I_3$, one obtains a constant $C$ depending only on
$\beta$ and $J$ so that
\begin{equation}\label{I4bnd}
I_4 \le 
C\|\sigma(\overline m){\cal B}v\|_2\langle h,{\cal B} h\rangle^{1/2}\le
{1\over 4}\langle h,{\cal B} h\rangle + 
4C^2\|\sigma(\overline m){\cal B}v\|_2^2\ .
\end{equation}
Then, from  (\ref{I1bnd}),  (\ref{I2bnd}), (\ref{I3bnd}) and (\ref{I4bnd}) 
we have
\begin{equation}\label{mag2}
 I_1+I_2+I_3+I_4 \le  {3\over 4}\langle h,{\cal B} h\rangle + C\|v\|_2^2.  
 \end{equation}
>From Theorem~\ref{(65)}, taking into account (\ref{mag3})  and (\ref{mag2}) we have
\begin{equation*}
\begin {split}
\frac{{\rm d}}{{\rm d}t}\phi (t)  &  \le  
- \langle h,{\cal B}h\rangle  + 4\langle f,{\cal B}f\rangle      
+  {3\over 4}\langle h,{\cal B} h\rangle + C\|v\|_2^2  
\cr & +\e  \bigl[{\cal F}(\overline m+v) - {\cal F}(\overline m)\bigr]   +
\e\|{\cal B}^{1/2}\bigl(\sigma(\overline m)x_1 \nabla \bigl({\cal B}v\bigr)\|_2^2. 
\end {split} \end{equation*}
Take $\e <\frac 14$, and using the fact that
${\cal B}$ is bounded, with a bound depending only on $\beta$  $J$ and $L$,
 Lemma~\ref{A2},
we have  (\ref{phibound1}).
\medskip 
 To get  (\ref{phibound})   we   estimate 
  $I_1$ through $I_4$  under the assumption (\ref{epsn5}).  We have 
\begin{eqnarray*}
I_1 &=& 
-2\langle g\sigma^2(\overline m)x_1^2  {\cal B}v,   ({\cal B}v)_{x_1} \rangle_{L^2}  \nonumber\\
&\leq&2\|g\sigma^2(\overline m)x_1^2\|_\infty\|{\cal B}v\|_2\|  ({\cal B}v)_{x_1}  \|_2 \nonumber\\
&\leq&2\|g\sigma^2 (\overline m)x_1^2\|_\infty \left ( \frac {2 \tilde \a \e_1 (1+ 2\e) } { (1-3\e)
\s (m_\b)}  \right )^{   1/2}\bigl[{\cal F}(\overline m+v) - {\cal F}(\overline m)\bigr]
\end{eqnarray*}
 where we used (\ref{CE2})  ( $\|\BB v \|^2_{L^2 (\R \times \L)} \le   2  \tilde \a (1+ 2\e)  [{\cal F}\big(m(t))-{\cal F}\big(\overline
m     \big)\bigr]$)   of Lemma~\ref{TT3},  (\ref{f2})  of Theorem~\ref{51} in the last step, together with the assumption (\ref{epsn5}).
Note that $\|g\sigma(\overline m)^2x_1^2\|_\infty$ is bounded by    a constant 
depending  only  on $\b$ and  $J$  by (\ref{decay}) and the hypothesis that $ |a(t)|
\le 1$, since (\ref{decay}) implies that $g$ is a rapidly decaying
bump function centered on $a(t)$.  
Other $L^\infty $ estimates  involving $x_1$ will be treated in the same way
without further mention. This is the only use made of $ |a(t)|
\le 1$.
Similarly,
\begin{eqnarray*}
I_2 &=& -2\langle \sigma(\overline m){\cal B}\bigl(\sigma(\overline m)'x_1^2{\cal B}v\bigr) e_1,
\nabla \bigl({\cal B}v) 
\rangle_{L^2}  \nonumber\\
&\leq&2\|\sigma(\overline m)\|_\infty\|\|{\cal B}\|\|\sigma(\overline m)'x_1^2\|_\infty
\|{\cal B}v\|_2\|\nabla \bigl({\cal B}v) 
\|_2 \nonumber\\
&\leq& 2\|\sigma(\overline m)\|_\infty\|\|{\cal B}\|\|\sigma(\overline m)'x_1^2\|_\infty
\left ( \frac {2 \tilde \a \e_1 (1+ 2\e) } { (1- 3\e)
\s (m_\b)}  \right )^{ 1/2} \bigl[{\cal F}(\overline m+v) - {\cal F}(\overline m)\bigr]
\end{eqnarray*}
In the same way, one obtains similar bounds for $I_3$ and $I_4$ and
then since all of the $\|\cdot\|_\infty$ terms are bounded {\it a--priori}
in terms of $\b$ and $J$, there is a constant $C$ depending only on
$\b$ and $J$ such that
$$I_1+I_2+I_3+I_4 \le C \e_1\bigl[{\cal F}(\overline m+v) - {\cal F}(\overline m)\bigr]\ .$$
Choosing $\e_1 = \e/C$, one has
\begin{equation}\label{halfway}
I_1+I_2+I_3+I_4 \le \e\bigl[{\cal F}(\overline m+v) - {\cal F}(\overline m)\bigr]\ .
\end{equation}
Recalling  notations (\ref{fdef})  and 
 (\ref{hdef}),   Theorem~\ref{(65)} and (\ref{halfway}) one has 
\begin{eqnarray}\label{tempbnd42}
\frac{{\rm d}}{{\rm d}t}\phi(t) &\le& 
-4\langle f,{\cal B}h\rangle_{L^2} -(2-\e)\langle h,{\cal B}h\rangle_{L^2}\nonumber\\
&+& 2 \e\bigl[{\cal F}(\overline m+v) - {\cal F}(\overline m)\bigr]\ .  
\end{eqnarray} 

Now, since ${\cal B}$ is  non negative,
\begin{eqnarray} \label{MM1}
-4\langle f,{\cal B}h\rangle_{L^2} -(2-\e)\langle h,{\cal B}h\rangle_{L^2}   
&=&-\langle (2-\e)^{1/2}h+2(2-\e)^{-1/2}f,{\cal B}
\bigl((2-\e)^{1/2}h+2(2-\e)^{-1/2}f\bigr)\rangle_{L^2}\nonumber\\ 
&+&
4(2-\e)^{-1}\langle f,{\cal B}f\rangle_{L^2} \nonumber\\
 & \le& 4(2-\e)^{-1}\langle f,{\cal B}f\rangle_{L^2}.  
 \end{eqnarray}

To bound this in terms of the excess free energy, 
one makes repeated use of (\ref{subs2}) of Lemma~\ref{TT2} 
together with the self adjointness and boundedness of ${\cal B}$,
to replace factors of $\sigma(\overline m)$ with factors of $\sigma(m_\b)$:
\begin{eqnarray}\label{MM2}
\langle f,{\cal B}f\rangle_{L^2} & = &
 \langle \sigma(\overline m){\cal B}v,{\cal B}\sigma(\overline m){\cal B}v\rangle_{L^2}\nonumber\\
&=& \langle \sigma (m_\b){\cal B}v,{\cal B}\sigma(\overline m){\cal B}v\rangle_{L^2}  +
\langle [ \sigma(\overline m){\cal B}v  - \sigma (m_\b){\cal B}v ]  ,{\cal B}\sigma(\overline m){\cal B}v\rangle_{L^2} \nonumber\\
 & \le &
 \sigma(m_\b)\langle {\cal B}v,{\cal B}\sigma(\overline m){\cal B}v\rangle_{L^2}
+ C\| \nabla \bigl({\cal B}v\bigr)\|_2\|{\cal B}v\|_2 \nonumber\\
&=&\sigma(m_\b)\langle {\cal B}^2v,\sigma(\overline m){\cal B}v\rangle_{L^2}
+ C\| \nabla \bigl({\cal B}v\bigr)\|_2\|{\cal B}v\|_2 \nonumber\\
&\leq& \sigma(m_\b)^2\langle {\cal B}^2v, {\cal B}v\rangle_{L^2}
+ 2 C\| \nabla \bigl({\cal B}v\bigr)\|_2\|{\cal B}v\|_2  \nonumber\\
&\leq&\sigma(m_\b)^2\tilde\alpha\langle {\cal B}v,{\cal B}v\rangle_{L^2}
+ 3C\| \nabla \bigl({\cal B}v\bigr)\|_2\|{\cal B}v\|_2  \nonumber\\
&=&\sigma(m_\b)^2\tilde\alpha^2\langle v,{\cal B} v\rangle_{L^2}
+ 4C\| \nabla \bigl({\cal B}v\bigr)\|_2\|{\cal B}v\|_2,   
\end{eqnarray}
where $C$ is constant derived from those in the cited lemmas.
Combining estimates (\ref{tempbnd42}), (\ref{MM1})  and (\ref{MM2}) we have 
\begin{equation*}
\frac{{\rm d}}{{\rm d}t}\phi(t) \le 
-\sigma(m_\b)^2\tilde\alpha^2 4(2-\e)^{-1} \langle v,{\cal B} v\rangle_{L^2}
+ 4C\| \nabla \bigl({\cal B}v\bigr)\|_2\|{\cal B}v\|_2  +
 2\e\bigl[{\cal F}(\overline m+v) - {\cal F}(\overline m)\bigr]\ .
 \end{equation*} 
By  the hypotheses (\ref{epsn5}),  (\ref{CE2}) of Lemma~\ref{TT3} ,  (\ref{f2})  of and Theorem~\ref{51} 
  we have
$$4C\| \nabla \bigl({\cal B}v\bigr)\|_2\|{\cal B}v\|_2  \le
\e \bigl[{\cal F}(\overline m+v) - {\cal F}(\overline m)\bigr]$$
for $\e_1$ sufficiently small.
By Lemma~\ref{A2} of  the Appendix 
$\langle v,{\cal B} v\rangle_{L^2} \le 
\frac 2  {(1-\e)}\bigl[{\cal F}(\overline m+v) - {\cal F}(\overline m)\bigr]$ for $\delta$ and $\kappa$
sufficiently small. 
Redefining $\e$, one has
\begin{equation*}
\frac{{\rm d}}{{\rm d}t}\phi(t) \le 
(1+\e)4\tilde \alpha^2\sigma(m_\b)^2
\bigl[{\cal F}(\overline m+v) - {\cal F}(\overline m)\bigr]\ .
\end{equation*}
which is the desired result since $\tilde \alpha^2\sigma^2(m_\b)  =
(1-\sigma(m_\b))^2$.

   \section {Proof of the main results  }

  We begin this section by proving several lemmas concerning the $L^1$ norm. The
first of these will be used in the proof of Theorem \ref {main1}  to control $ |a-
a(t)|$ when $a$ is given by \eqref {conser}   and $a(t)$  denotes the shift from the origin of the closest  
front to $m(t)$, see    \eqref {madef}.

\bigskip
\begin{lm}\label{61} 
  Let $w $ be a function such that   $w  -\rm
sgn(x_1)m_{\b} $,  $x_1 \in \R$, is integrable   and   $b$ be fixed by the
condition  
\begin{equation}\label{fixesa1} \int_{\R \times \L}  \bigl(w(x) - \bar m_b(x)\bigr){\rm d}x = 0. \end {equation}
Then for any $c$
\begin{equation}\label{notoofar1}|b-c| \le 
\frac{1}{ 2m_\b} \frac 1 {L^d} \int _{\R \times \L}  \bigl|w(x) - \bar m_c (x)\bigr|{\rm d}x.  
\end {equation} 
In particular, for any solution $m(t)$ of \eqref {1.1}  and any $t$ such
that $m(t)-\rm sgn(x_1)m_\b $ is integrable, 
\begin{equation}\label{notoofar}|a(t) - a| \le 
\frac{1}{ 2m_\b}  \frac 1 {L^d}   \int _{\R \times \L} \bigl|m(x,t) - \bar m_{a(t)}(x)\bigr|{\rm d}x  
\end {equation} 
where $a$ is fixed by the condition that
\begin{equation*}
  \int \bigl(m(x,0) - \bar m_a(x)\bigr){\rm d}x = 0.
\end {equation*} 
\end {lm}

\bigskip
\noindent{\bf Proof:}\ First, since $\int_{\R} \bar m_0'(x_1){\rm d}x_1 = 2m_\b > 0$,
there is exactly one $b$ such that \eqref {fixesa1} holds. Next,
adding and subtracting $\bar m_{c}$, one sees
$$\int_{\R \times \L} \bigl(w(x) - \bar m_{c}(x)\bigr){\rm d}x = 
\int_{\R \times \L} \bigl(\bar m_b(x) - \bar m_{c}(x)\bigr){\rm d}x.$$
Also, it is clear that
$$\int_{\R \times \L} \bigl(\bar m_b(x) - \bar m_{c}(x)\bigr){\rm d}x = 2m_\b(b-c) L^d $$
and \eqref {notoofar1}  easily follows. \qed
\bigskip

\begin{lm}\label{5.2a}    Let $w$ be any function such that
\begin{equation*} 
 \int_{\R \times \L}  |w(x)|^2(1+ x_1^2){\rm d}x < \infty\ .
 \end {equation*}    
For any $0 <\d < 1$ so that
$$C(\d, L) = \biggl(\int_{\R \times \L} (1 + x_1^2)^{-(1+\d)/2}{\rm d}x\biggr)^{1/2} < \infty $$
we have 
$$ \|w\|_1 \le C(\d, L )\|(1 + x_1^2)^{1/2}w\|_2^{(1+\d)/2}\|w\|_2^{(1-\d)/2}.$$
\end {lm}

\bigskip
\noindent{\bf Proof:}  Let $p = (1+\d)/2$  and observe that  
\begin{equation}  \label {l2p1} \begin {split}  
\int_{\R \times \L}  |w(x)|{\rm d}x = 
&\int_{\R \times \L}  (1 + x_1^2)^{-p/2}(1 + x_1^2)^{p/2}|w(x)|{\rm d}x   \cr \le 
\biggl(&\int_{\R \times \L}  (1 + x_1^2)^{-p}{\rm d}x\biggr)^{1/2}
\biggl(\int_{\R \times \L} (1 + x_1^2)^{p}|w(x)|^2{\rm d}x\biggr)^{1/2}. 
\end {split}  \end {equation} 
Jensen's inequality, for   $p<1$,  implies
\begin{equation}  \label {l2p2} \frac{1}{ \|w\|_2^2}\int_{\R \times \L} (1 + x_1^2)^{p}|w(x)|^2{\rm d}x \le
\biggl({\frac{1}{ \|w\|_2^2}}\int_{\R \times \L} (1 + x_1^2)|w(x)|^2{\rm d}x \biggr)^p. 
\end {equation} 
The result easily follows from \eqref {l2p1} and  \eqref {l2p2}.  \qed 
 
\bigskip
Since  $ \phi (t)$ is defined in term of  moments of ${\cal B}v$ instead of $v$,  see \eqref{phidef},
we need one more lemma to apply the previous one.
\bigskip
\begin{lm}\label{5.3a}   There is a finite constant $C$
depending only on $\b$ and $J$ so that for all $t$ such that
$|a(t)| \le 1$ and $\|v(t)\|_2 \le 1$,
\begin{equation}  \label {orata1} \|(1+x_1^2)^{1/2}v(t)\|_2^2 \le C(\b,J)\phi(t)
\end {equation}
and
\begin{equation}  \label {orata2} \bigl(\phi(t) -L^d\bigr) \le C(\b,J)\|x_1 v(t)\|_2^2.
\end {equation}
\end {lm}

\bigskip
\noindent{\bf Proof:} Let ${\cal D}$ be the operator defined  in \eqref {const1}
$${\cal D}w = {1\over \b(1-m_\b^2)}w - J\star w\ .$$
As we have pointed out in  Section 5  this operator is bounded
and has a bounded inverse on $L^2 (\R \times \L)$. Then, using once more
the  rules for commuting convolution and multiplication by $x_1$, see \eqref {mars1},  we have 
\begin{equation*}  
\begin {split}   
\|x_1v\|_2 \le &\|{\cal D}^{-1}\|\|{\cal D}x_1v\|_2    \cr =
&\|{\cal D}^{-1}\|\|x_1{\cal D}v - (x_1J)\star v\|_2  \cr \le 
&\|{\cal D}^{-1} \| \bigl( \|x_1{\cal D}v\|_2 + \|(x_1J)\|_1\| v\|_2\bigr) \cr \le 
& C \bigl(\|x_1{\cal D}v\|_2 + \|v\|_2\bigr)
\end {split} 
\end {equation*}
for some constant $C$ depending only on $\b$ and $J$.
Next, taking into account \eqref {TT.1}
\begin{equation*}   \begin {split}   
\|x_1{\cal D}v\|_2  = &\|x_1{\cal B}v - x_1\tilde g v\|_2 \cr \le  
 &\|x_1{\cal B}v\|_2 + \|x_1\tilde g\|_\infty\| v\|_2,  
 \end {split} 
\end {equation*}
where, recall \eqref {mars2}, 
$$\tilde g(x_1) = \sigma(m_\b)^{-1} - \sigma(\bar m(x_1))^{-1}, \qquad x_1 \in \R.$$
 Since
   the hypothesis  $|a(t)| \le 1$,
$\|x_1\tilde g\|_\infty \le C$ for some constant $C$
depending only on $\b$ and $J$.  Thus we have 
\begin{equation*}
 \|x_1{\cal D}v\|_2 \le   C\bigl(\|x_1{\cal B}v\|_2+ \|v\|_2\bigr).  
\end {equation*}
 Finally,
\begin{equation*} 
\begin {split}   
&\sigma(m_\b)\|x_1{\cal B}v\|_2^2 = \int_{\R \times \L}  x_1^2(\sigma(m_\b) - \sigma(\bar m))
\bigl({\cal B}v\bigr)^2{\rm d}x + \int_{\R \times \L}  x_1^2\sigma(\bar m)
\bigl({\cal B}v\bigr)^2{\rm d}x  \cr \le 
&\|x_1^2(\sigma(m_\b) - \sigma(\bar m))\|_\infty\|{\cal B}v\|_2^2 +
\int_{\R \times \L}   x_1^2 \sigma(\bar m)
\bigl({\cal B}v\bigr)^2{\rm d}x \end {split}
\end {equation*}
and the sup norm is again bounded by some constant 
$C(\b,J)$ depending only on $\b$ and $J$ by
\eqref {decay} and the hypothesis  $|a(t)| \le 1$.
Combining these estimates, one easily obtains
$$\int_{\R \times \L} (1+x_1^2)v^2 dx  \le C \biggl(
\int_{\R \times \L} x_1^2\sigma(\bar m)\bigl({\cal B}v\bigr)^2{\rm d}x + \|v\|_2^2\biggr)$$
which yields \eqref {orata1} since $\|v(t)\|_2 \le 1$ by hypothesis.
The proof of  \eqref {orata2} simply reverses the above steps. With $C$
changing from line to line, one easily obtains
\begin{equation*}    \begin {split}   
\int_{\R \times \L}  x_1^2 \sigma(\bar m)
\bigl({\cal B}v\bigr)^2{\rm d}x \le 
&C\bigl(\|x_1{\cal B}v\|_2^2    + \|v\|_2^2\bigr) \cr  \le
&C\bigl( \|x{\cal D}v\|_2^2   + \|v\|_2^2\bigr) \le 
 C\bigl( \|{\cal D}x_1v\|_2^2  + \|v\|_2^2\bigr)  \cr  \le
&C\bigl( \|x_1v\|_2^2  + \|v\|_2^2\bigr) \end {split}
\end {equation*}
and this complete the proof. 
\qed

\vskip0.5cm
\noindent{\bf { Proof of Theorem  \ref {main1} : }} First, 
fix $\e>0$, and then choose $\delta_1$,  $\kappa_1$ and $\e_1$
small enough so that both   the following three estimates hold
under the condition that
\begin{equation}  \label {pivotal1}{\cal I}(m(t)) \le \e_1\bigl[{\cal F}(m(t)) - {\cal F}(\md)\bigr]
\end {equation}
for all $t$ such that  $|a(t)|\le 1$, $\|v(t)\|_2 \le \delta_1$ and $\|v(t)\|_{W^{s+1,2}} 
\le
\kappa_1$,  $s > \frac {D}2$:
\begin{equation}  \label {onekeyend} \frac{{\rm d}}{{\rm d}t}\bigl[{\cal F}\big(m(t))-{\cal F}\big(\bar m  \big)\bigr]
\le -9(1-\e)(1-\sigma(m_\b))^2  \frac {\bigl[{\cal
F}\big(m( t))-{\cal F}\big(\bar m  
\big)\bigr]^2}
 {{ \phi ( t )}  } 
\end {equation}
and
\begin{equation}  \label {phiboundend} \frac{{\rm d}}{{\rm d}t}\phi(t) \le 
(1+\e)4(1-\sigma(m_\b))^2
\bigl[{\cal F}(\bar m+v) - {\cal F}(\bar m)\bigr].
\end {equation}
This is possible by Theorems  \ref {51b} and  \ref {41}.
By  \eqref {phibound1} of Theorem \ref {41} and  Lemma \ref {A2} in the appendix   decreasing $\d_1>0$ and
$\kappa_1>0$ if need be, we have for a finite constant  $B$ and   $c(\kappa_1)$  
\begin{equation}  \label {phiboundendvar}\frac{{\rm d}}{{\rm d}t}\phi(t) \le 
 B\bigl[{\cal F}(\bar m+v) - {\cal F}(\bar m)\bigr]
 \end {equation}
 and
\begin{equation}  \label {lemmabnd}{1\over 4} \g (L) \|v\|_2^2 \le \bigl[{\cal F}(\bar m+v) - {\cal F}(\bar m)\bigr]
\le  c(\kappa_1)\|v\|_2^2.
 \end {equation}
   Next define $\delta_0$ by 
\begin{equation}  \label {donedef} \delta_0 = \frac { \delta_1 \sqrt { \gamma (L)}}{4(c(\kappa_1)+1)}
 \end {equation}
where  $\g(L)$ and $c(\kappa_1)$ are the constants in \eqref {lemmabnd}. 
  Theorem  \ref {smoo}    applied with the values of 
$\delta_0$, $\delta_1$ and $\kappa_1$ fixed above, 
guarantees the existence of an $\e_0>0$ and a $t_0$ so that when the initial
data satisfies
$\|m_0 - \bar m\|_2 \le \e_0$, the solution to \eqref {1.1} satisfies
\begin{equation}  \label {initcond1}\|v(t_0)\|_2 \le \delta_0
 \end {equation}  and 
$$\|v(t)\|_{W^{s,2}} \le \kappa_1$$ 
for all $t\ge t_0$ such that $\|v(t)\|_2 \le \delta_1$. 
We have from Theorem \ref {smoo}    that  
$$\int_{\R \times \L}  \bigl (  x_1  (m(x,t_0) - \bar m_0(x)\bigr)^2{\rm d}x \le 2c_0$$
where $c_0$ is the constant specified in the hypotheses of Theorem  \ref {main1}.
Clearly then,
$$\|x_1v(t_0)\|^2_2 \le 2\bigl(\|x_1(m(t_0) - \bar m_0) \|^2_2 + 4m_\b a(t_0) L^d \bigr) 
.$$ 
By  Theorem \ref {5}     we may suppose, further decreasing $\d_1$ if
need be, that $4m_\b a(t_0) L^d \le c_0$. Then
$$\|x_1v(t_0)\|_2^2 \le 5c_0$$
and hence, by \eqref {orata2} of Lemma \ref  {5.3a}, 
\begin{equation}  \label {initcond2} \phi(t_0) \le \tilde c_0
 \end {equation}
where $\tilde c_0$ is a finite constant depending only on $c_0$, $\b$,  $J$ and $L$. 
 Hence, writing $f(t) = \bigl[{\cal F}(m(t)) - {\cal
F}(\md)\bigr]$,
  we  have to control on the values of both $f(t_0)$ and $\phi(t_0)$
through \eqref {initcond1} and \eqref {initcond2}.

The time $t_0$
is the time we have to wait for the smoothing properties of 
the equation \eqref {1.1} to regularize our data enough that the estimates
above all hold, and it fixes the left end of the interval on which we shall
work. To fix the right end, which we shall eventually show to be $+\infty$,
define
$$T_0 = \min\{\ 
\inf\{\ t>t_0\ |\ \|v(t)\|_2 \ge \delta_1/2\ \}\ ,
\ \inf\{\ t>t_0\ |\    |a(t)| \ge 1\ \}\ \}     \ .$$
Then, uniformly on the interval $(t_0,T_0)$, both of  the estimates 
\eqref {phiboundendvar} and  \eqref {lemmabnd} holds. Moreover for   
 those $t$ in  $(t_0,T_0)$  such
that \eqref {pivotal1} holds, one also has  \eqref {onekeyend} and \eqref
{phiboundend}. Hence  we have the following alternative:

One the one hand, in case 
\begin{equation*} 
{\cal I}(m(t)) \le \e_1\bigl[{\cal F}(m(t)) - {\cal F}(\md)\bigr]
 \end {equation*}
\begin{equation}  \label {oneway} \begin {split}   
&\frac{{\rm d}}{{\rm d}t}f(t)
\le -\tilde A{f(t)^2\over  \phi(t)} \cr
&\frac{{\rm d}}{{\rm d}t}\phi(t) \le 
\tilde B f(t) \end {split}
 \end {equation}
where $\tilde A$ and $\tilde B$ by
\begin{equation*}  
 \begin {split}    &\tilde A = 9(1-\e)(1-\sigma(m_\b))^2\cr &\tilde B =
4(1+\e)(1-\sigma(m_\b))^2.
 \end {split}
 \end {equation*}
  On the other hand,    in case 
\begin{equation*} 
 {\cal I}(m(t)) \ge {\e_1\over2} \bigl[{\cal F}(m(t)) - {\cal
F}(\md)\bigr],
 \end {equation*}
\begin{equation*}  
 \begin {split} 
&\frac{{\rm d}}{{\rm d}t}f(t)
\le -{\e_1\over2}f(t) \cr
&\frac{{\rm d}}{{\rm d}t}\phi(t) \le Bf(t).
 \end {split}
 \end {equation*}
  In the application of the system of differential inequalities \eqref {system}, it
is the {\it ratio} of the constants $A$ and $B$ that determines 
the exponent $q$, see Theorem 5.1 of   \cite{CCO1}.
Indeed,
$$q = {(A/B)\over (A/B)+1}.$$
 The values of $A$ and $B$ themselves can be changed, 
keeping this ratio fixed, simply by rescaling the
time $t$. Therefore we define
$$A = {\tilde A\over \tilde B}B $$
and observe that
$${\e_1\over2}f(t) = {\e_1\over2 A f(t)}Af(t)^2 \ge 
{\e_1\over2 A f(t)}A{f(t)^2\over \phi}$$
since $\phi(t) \ge 1$ by definition. Now, by \eqref  {lemmabnd} we may further
decrease
$\d_1>0$ if need be to ensure that 
$${\e_1\over2 A f(t)} \ge 1\ .$$
Doing so, we have that in case 
$${\cal I}(m(t)) \ge {\e_1\over2} \bigl[{\cal F}(m(t)) - {\cal F}(\md)\bigr]$$
\begin{equation}  \label {otherway2} \begin {split} 
&\frac{{\rm d}}{{\rm d}t}f(t)
\le -A{f(t)^2\over \phi(t)} \cr
&\frac{{\rm d}}{{\rm d}t}\phi(t) \le Bf(t)\
 \end {split}
 \end {equation}
 where
$${A \over B} = {\tilde A\over \tilde B}\ .$$
Now suppose that at $t_0$,
$${\cal I}(m(t_0)) > {\e_1\over2}
\bigl[{\cal F}(m(t_0)) - {\cal F}(\md)\bigr]\ .$$
Define
$$t_1 = \inf\{\ t > t_0\ |\ {\cal I}(m(t)) \le {\e_1\over2}
\bigl[{\cal F}(m(t)) - {\cal F}(\md)\bigr]  \}\ ,$$
\begin{equation*}   \begin {split} 
&t_2 = \inf\{\ t > t_1\ |\ {\cal I}(m(t)) \ge \e_1
\bigl[{\cal F}(m(t)) - {\cal F}(\md)\bigr]  \}\ ,\cr
&t_3 = \inf\{\ t > t_2\ |\ {\cal I}(m(t)) \le \frac {\e_1} 2
\bigl[{\cal F}(m(t)) - {\cal F}(\md)\bigr]  \}\ ,\end {split}
 \end {equation*}
 and so forth.  We follow the usual convention that if there is no 
$t < T_0$ satisfying
the condition, the infimum is set to be $T_0$. Notice that since ${\cal I}(m(t))$ and ${\cal F}(m(t))$ are  continuous function of $t$,  $t_3>t_2>t_1> t_0$.
The sequence of times $t_j$ can have no limit point except possibly
$T_0$, since at such a limit point, the continuous function  ${\cal I}(m(t))$ 
 would take on two values.

If at $t_0$,
$${\cal I}(m(t_0)) \le \frac {\e_1} 2
\bigl[{\cal F}(m(t_0)) - {\cal F}(\md)\bigr]\ ,$$
one would define
$$t_1 = \inf\{\ t > t_0\ |\ {\cal I}(m(t)) \ge \e_1
\bigl[{\cal F}(m(t)) - {\cal F}(\md)\bigr]  \}\ ,$$
and then proceed as above with the opposite alternation.

In either case, one produces a sequence of intervals $[t_j, t_{j+1}]$ on which 
\eqref {oneway} and \eqref {otherway2} hold in successive alternation.
On each of these intervals, we may apply Theorem 5.1 of \cite{CCO1}.
To put all of these estimates together in a transparent way, we
rescale the intervals on which \eqref {otherway2} holds. Supposing that
\eqref {otherway2} holds on $[t_0,t_1]$, define
$$ s(t) = {A\over \tilde A}(t-t_0)\qquad{\rm and}\qquad 
s_1 = {A\over \tilde A}(t_1-t_0)$$
for $t_0 < t < t_1$, 
$$ s(t) = s_1 + (t-t_1)\qquad{\rm and}\qquad 
s_2 = s_1 +(t_2-t_1)$$
for $t_1 < t < t_2$, 
$$ s(t) = s_2 + {A\over \tilde A}(t-t_2)\qquad{\rm and}\qquad 
s_3 = {A\over \tilde A}(t_3-t_2)$$
for $t_2 < t < t_3$, and so forth in alternation.
It follows that 
\begin{equation*} 
  \begin {split} 
&{{\rm d}\over {\rm d}s}f(s)
\le -\tilde A{f(s)^2\over \phi(s)} \cr
&{{\rm d}\over {\rm d}s}\phi(s) \le \tilde Bf(s)
\end {split}
 \end {equation*}
for all $s$ with $0 \le s \le s(T_0)$. By Theorem 5.1 of \cite{CCO1},
\begin{equation*} 
 \begin {split} 
f(s) &\le f(0)^{1-q} \phi(0)^q 
\biggl({\phi(0)\over f(0)} + (\tilde A +\tilde B)s\biggr)^{-q}\cr
\phi(s) &\le
f(0)^{1-q}\phi(0)^{q}
\biggl({\phi(0)\over f(0)} + (\tilde A +\tilde B)s\biggr)^{1-q}
\end {split}
 \end {equation*}
where 
$$q = {\tilde A\over \tilde A + \tilde B}$$
and where $f(0)$ and $\phi(0)$ are bounded by \eqref {initcond1} and \eqref {initcond2}.
Now, for any $\d>0$, we can choose $\e$ so that
\begin{equation*} 
{\tilde A\over \tilde A + \tilde B} = {9\over 13} - \d\ .
\end {equation*}
 
We shall now show that for $\d$ small enough, $|a(t)| \le 1/2$ for
all $t\le T_0$.
Then by Lemmas \ref {5.2a},  \ref {5.3a}  and  estimate \eqref {lemmabnd}, 
\begin{equation}  \label {landofoz}   \begin {split} 
\|v(s)\|_1 \le 
&C(\d, L )\|(1+x_1^2)^{1/2}v(s)\|_2^{(1+\d)/2}\|v(s)\|_2^{(1-\d)/2}  \cr  \le
&C(\d,L)C(\b,J)^{(1+\d)/4}(\frac 4 {\g(L)}) ^{(1-\d)/4}
\phi(s)^{(1+\d)/4}f(s)^{(1-\d)/4} \cr \le
&C(\d,L)C(\b,J)^{(1+\d)/4}(\frac 4 {\g(L)}) ^{(1-\d)/4}  f(0)^{(1-q)/4 }\phi(0)^{q/4}
\biggl({\phi(0)\over f(0)} + (\tilde A +\tilde B)s\biggr)^{(1-2q +\d)/4}.\cr
\end {split}
 \end {equation}
The right hand side is decreasing
for $\d < 5/26$, and we now choose $\d$ to be at least this small.
Moreover, the value at $s=0$ can be made arbitrarily small by decreasing
$\d_1$. We now do so, if need be, to ensure that
$$\|v(s)\|_1 \le m_\b/2$$ for all $s \le s(T_0)$.
Hence,  for  $a$, as in \eqref {conser},  by Lemma  \ref {61},    we have 
$$|a(t) - a| \le 1/4  \frac 1 {L^d} $$
for all $t\le T_0$. But then by Lemma  \ref {61} again, this implies that
$|a(t)| < \frac 1 {2L^d}$ for all $t \le T_0$.
Hence if $T_0 < \infty$, it is because $\|v(T_0)\|_2  = \delta_1/2$.
But since \eqref {lemmabnd} is still valid with the same constants   on 
the closed interval
$[t_0,T_0]$, and since the excess free energy is monotone decreasing, we have
\[    \begin {split} 
 {\delta_1^2  \over 4} & =  \|v(T_0)\|_2^2  \cr \le
&\frac 4 {\g(L)} \bigl(\FF(\bar m + v(T_0)) - \FF(\bar m)\bigr) \le 
\frac 4 {\g(L)} \bigl(\FF(\bar m + v(t_0)) - \FF(\bar m)\bigr)  \cr \le
&c (\kappa_1)\frac 4 {\g(L)} \|v(t_0)\|_2^2 \le c (\kappa_1) \frac 4 {\g(L)} \delta_0^2.
\end {split}
\]
This contradicts \eqref {donedef}, and hence $T_0 < \infty$ is not possible.
We now clearly have \eqref {result111}  since
$$s(t) \ge \min\biggl\{{A\over \tilde A}\ , \ 1\biggr\}(t-t_0)\ .$$
Also from this and \eqref {landofoz}, we have
$$\|m(t) - \bar m_{a(t)}\|_1 \le c_2(1 + c_1 t)^{-(5/52 - \d)}\ .$$
But
\[    \begin {split}  \|m(t) - \bar m_a\|_1       \le &
 \|m(t) - \bar m_{a(t)}\|_1 +\|\bar m_a - \bar m_{a(t)}\|_1  \cr=
&\|m(t) - \bar m_{a(t)}\|_1 + 2m_\b L^d |a - a(t)|  \cr \le
&2 \|m(t) - \bar m_{a(t)}\|_1
\end {split}
\]
  by \eqref {notoofar}.  Hence \eqref {result222} follows as well.
\qed

   \section {Proof of  Smoothing estimates } 
 
  The main goal of this section is to deduce the regularity properties stated in  Theorem   \ref {smoo}
for the  derivatives  of   the solution $m(t)$ of  equation \eqref {1.1}  that starts sufficiently close in the $L^2$ norm to
$\bar m_b$ for some $b$:

The proof depends on several intermediate results concerning the evolution
of $m(t) - \bar m_b$ and its derivatives for {\it fixed} $b$. To
simplify the notation, we will write $\bar m$ instead of $\bar m_b$, so one
should keep in mind that in this section $\bar m$ is not necessarily the 
antisymmetric, increasing instanton.

Let $v = m -\bar m$. Notice we are {\it  not} assuming that $\bar m = \bar m_{a(t)}$, so this definition of $v(t)$
differs slightly from the one used in the rest of the paper. However, for most of this section, it is the most convenient notation. It saves us from keeping  $\dot a(t)$ terms throughout the many calculations that follow.
{}From the evolution equation
$${\partial m\over \partial t} = \nabla \cdot \left(\nabla m - \beta(1-m^2)J*\nabla m\right)$$
and the eigenvalue equation
$$\beta(1-\bar m^2)J*\bar m' = \bar m'$$
we deduce the following evolution equation for $v$
 
\begin{equation}  \label {vevo} \begin {split}   
\frac{\partial v}{\partial t}
 &= \nabla \cdot \left(\nabla v - \beta(1-\bar m^2)\nabla J*v\right)  
\cr & + \beta \nabla \cdot e_1\left(v(v+2\bar m)J*\bar m'\right) 
+ \beta \nabla \cdot\left(v(v+2\bar m)\nabla J*v \right).  
\end {split}  \end {equation} 
Here and in what follows, $e_1$ denote the unit vector in the $x_1$ direction. 
Define
\begin{equation*} 
\Psi := \beta \left(e_1 J\star \bar m' + \nabla J\star v\right)  \qquad{\rm and}\qquad 
\Phi :=  2\bar m \Psi\ . 
\end {equation*} 
 We can write \eqref {vevo} as
\begin{equation}  \label {vevo2}  
\frac{\partial v}{\partial t}
=\nabla \cdot \left(\nabla v - \beta(1-\bar m^2)\nabla J*v\right)  +
\nabla\cdot \left(v^2\Psi\right) +  \nabla\cdot \left(v\Phi\right)\ . 
  \end {equation} 
Since $\Psi$ and hence $\Phi$ depend on $v$, both the second and the third terms on the right in 
\eqref {vevo2} are nonlinear in $v$. However, because of the convolution with $J$, the dependence on $v$
that enters through these terms is harmless as far as smoothness of $v$ is concerned.  
For any multindex $\alpha$,
 denote by  $D^{\alpha}$  the corresponding differentiation operator. 
  Since both
$\|m\|_\infty \le 1$ and  $\|\bar m\|_\infty \le m_\beta \le 1$, $\|v\|_\infty \le 2$,    we have
$$
\|D^{\alpha}(J*v)\|_\infty = \|(D^{\alpha}J)*v\|_\infty \le  2\|(D^{\alpha}J)\|_\infty\ ,
$$
 independent of $v$. Then, since $\bar m$ is smooth, there exist finite constants $C_\alpha$ depending only on 
$J$ and $\alpha$ so that
\begin{equation}  \label {evo3}  
\|D^{\alpha}\Psi\|_\infty \le  C_\alpha \qquad{\rm and}\qquad 
\|D^{\alpha}\Phi\|_\infty \le  C_\alpha\ .  \end {equation} 
 Our first goal is to study the smoothing properties of  \eqref {vevo2}. We shall show
that on any interval of time on which $\|v(t)\|_2$ stays bounded, solutions immediately
develop derivatives of all orders even if the initial data is not smooth. To use
this, we need to know that $\|v(t)\|_2$ stays bounded in some interval of the origin.
Later of course we shall see that if $\|v(0)\|_2$ is small enough, 
this holds globally in time. 
 In the next Lemma  and in what follows, $C$ will denote a constant depending 
only on $J$ and $\beta$ but otherwise changing from line to line.

\begin{lm}    Let $v$ be a solution of \eqref {vevo}. Then there is a finite constant
$C$  depending only on $J$ and $\beta$ so that for all $t>0$, 
$$\|v(t)\|_2^2 \le  e^{Ct}\|v(0)\|_2^2\ .$$
 \end {lm}

\noindent{\bf Proof:}  
{}From  \eqref {vevo2} we have
\begin{equation*}   
\begin {split}   
& \frac{d }{d t}\| v\|_2^2  
= 2\int v{\partial v\over \partial t}{\rm d}x   = \cr &
-2\| \nabla v\|_2^2 + 2\beta\int \nabla v\cdot  \left((1-\bar m^2) \nabla J* v \right){\rm d}x  
\cr & -  2\int \nabla v\cdot \left(v^2\Psi + v\Phi\right){\rm d}x  
\cr & \le   -2\| \nabla v\|_2^2  + C\|\nabla v\|_2\|v\|_2
 \end {split}  
 \end {equation*} 
where in the last line we have used the bound $\|v\|_\infty \le 2$. 
Completing the square leads to
${\displaystyle {d\over dt}\| v\|_2^2 \le C\| v\|_2^2}$, and the result follows directly. \qed

\medskip
\begin{lm}\label{lm1}     Let $v$ be a solution of \eqref {vevo} and suppose that for some 
finite $\delta$ and positive $T_\delta$,
$$\|v(t)\|_2^2 \le \delta\qquad{\rm for \ all}\quad t \le T_\delta\ .$$
Then
$$\| \nabla v(t)\|_2^2 \le \frac{\delta}{2t} +  C \delta   \qquad{\rm for \ all}\quad t \le T_\delta$$
where $C$ is a constant depending only on $J$ and $\beta$.
 \end {lm}

\noindent{\bf Proof:} We begin with the $L^2$ norm of the first derivatives.
\begin{equation}  \label {difineq} \begin {split} &  \frac{ d }{ d t}\| \nabla v\|_2^2 
= -2\int \Delta v{\partial v\over \partial t}{\rm d}x \cr &= 
-2\| \Delta v\|_2^2 + 2\beta\int \Delta v\nabla\cdot \left((1-\bar m^2)\nabla J*v\right){\rm d}x 
\cr & +  2\int \Delta v \nabla\cdot \left( v^2\Psi + v\Phi\right){\rm d}x.  
\end {split}  \end {equation} 
By the Schwarz inequality, this is no more than
$$ 
-2\| \Delta v\|_2^2 + 2 \| \Delta v\|_2\left[
\beta\|\nabla\cdot  (1-\bar m^2)\nabla J* v)\|_2 + \|\nabla\cdot v^2\Psi\|_2
 + \|\nabla  \cdot (v\Phi)\|_2\right]
$$
Now, 
\begin{equation}  \label {evo8}
\|\nabla \cdot \left((1-\bar m^2)\nabla J* v\right)\|_2 
 \le 
\|2\bar m' \nabla  J*v\|_2 + \| (\Delta J)* v\|_2  
 \le 
C\|v\|_2 . 
\end {equation} 
 Then by \eqref {evo3} and the {\it a-priori} estimate  $\|v\|_\infty \le 2$, 
$$ \|\nabla\cdot v^2\Psi\|_2  \le  C(\|\nabla v\|_2 + \|v\|_2)
\qquad{\rm and}\qquad 
 \|\nabla  \cdot (v\Phi)\|_2 \le  C(\|\nabla v\|_2 + \|v\|_2)\ .
 $$
Combining this with \eqref {evo8}, we obtan
$$\frac {d} {dt} \| \nabla v\|_2^2 \le
-2\| \Delta v\|_2^2 + \| \Delta v\|_2C(\|\nabla v\|_2 + \|v\|_2). $$
We use half of our dissipative term $-2\| \Delta v\|_2^2$ to eliminate
reference to $\nabla v$ and $\Delta v$ in the positve part of this bound. To do so,
note that by the Schwarz inequality
$$ \| \nabla v\|_2^2 = -\int (\Delta v) v{\rm d}x \le \| \Delta v\|_2\| v\|_2\ .$$
Therefore,
$\| \Delta v\|_2\| \nabla v\|_2 \le  \| \Delta v\|_2^{3/2}\|  v\|_2^{1/2}$.
By the arithmetic--geometric mean inequality, 
$$\|\Delta v\|_2^{3/2}\|v\|_2^{1/2} \le {3\epsilon\over 4}\|\Delta v\|_2^2 + \frac 1 {4 \e^3}
\|v\|_2^2\  $$
 for any $\epsilon > 0$. Even more simply
${\displaystyle \|\Delta v\|_2\|  v\|_2 \le {\epsilon \over 2}\|\Delta v\|_2^2 + 
{1 \over 2\epsilon}\|  v\|_2^2}$, and thus, 
$${d\over dt}\| \nabla v\|_2^2 \le
\left( -2 + {\epsilon 5C\over 4}  \right)  \| \Delta v\|_2^2 + {3C\over 4\epsilon}\|v\|_2^2.
$$
Again by the Schwarz inequality
$$\| \Delta v\|_2^2 \ge {\| \nabla v\|_2^4\over \| v\|_2^2}\ .$$
Using this, and choosing 
$\epsilon$ so that $5C\epsilon \le 4$,
one finally obtains
$${d\over dt}\| \nabla v\|_2^2 \le
-{\| \nabla v\|_2^4\over \| v\|_2^2} + C\|v\|_2^2\ .$$

Now by hypothesis, for all times $t$ under consideration, we have the bound
$\| v\|_2^2 \le \delta$. Letting $x(t)$ denote the value of $\| \nabla v \|_2^2$
at time $t$, we then have the differential inequality
\begin{equation*}  
\frac  d  {dt} x \le -\frac{x^2}{\delta} + C\delta\ . 
\end {equation*} 
Introducing $y = 1/x$, one obtains a differential inequality of the form
$${d\over dt}y \ge {1\over \delta}  - C\delta y^2.$$
Now let $y_*$  $$y_* = {1\over  \sqrt {2 C}\delta} \ ,$$
so that  for $0\le y\le y_*$, $y'\ge \frac 1 {2 \delta} $, and for
$y_*\le y\le \sqrt 2 y_{*}$, we have at least that $y'\ge 0$. This means
that $y$ increases with rate at least $1/2\delta$ until $y_*$ is reached.
At this point it is still increasing, and it continues to increase until
$ \sqrt 2 y_{*}$, and it never again passes below this value, and hence 
never again below $y_*$ either.  Therefore
$$y(t) \ge \min\{t/2\delta\ ,\ y_*\}\qquad{\rm for\ all}\quad t  >  0\ ,$$
and hence
$$x(t) \le \max\{2\delta/t\ ,\  \sqrt {2C} \delta \} \le 2\delta/t +  \sqrt {2C} \delta 
\qquad{\rm for\ all}\quad t > 0\ .$$
This proves 
the stated assertion about $\| \nabla v\|_2^2$.  \qed

\medskip
 
We next consider the second derivatives, where a new feature emerges.

\medskip
\begin{lm}  Let $v$ be a solution of \eqref {vevo} and suppose that for some 
$\delta>0$ and $T_{\delta}>0$,
$$\| v(t)\|_2^2 + \| \nabla v(t)\|_2^2 \le \delta  \qquad{\rm for \ all}\quad t \le T_\delta\ .$$
Then
$$\| \Delta v(t)\|_2^2 \le \frac{ \delta }{2t} +   C \delta \qquad{\rm for \ all}\quad t \le T_{\delta}$$
where $C$ is a constant depending only on $J$ and $\beta$. 
 \end{lm}

\noindent{\bf Proof:} 
\begin{equation}  \label {evo20} \begin {split}    
\frac{d }{dt}\| \Delta v\|_2^2 
&=2\int [(-\Delta)^{2}v]{\partial v\over \partial t}{\rm d}x   = 
-2\| \nabla \Delta v\|_2^2  
 +  2\beta\int [\nabla \Delta v]\cdot \nabla {\rm div} \left((1-\bar m^2)\nabla J* v\right){\rm d}x  \cr &
+  2 \int [\nabla \Delta v]  \cdot \nabla {\rm div} \left[
\left(v^2\Psi\right)  +  \left(v\Phi\right) \right]{\rm d}x.   
\end {split}  \end {equation}   
By the Schwarz inequality, this is no more than
\begin{equation} \label {evo21}
-2\| (-\Delta)^{3/2}v\|_2^2 + C \| (-\Delta)^{3/2}v\|_2\left[
\|\Delta((1-\bar m^2)J*\nabla v)\|_2 
+\|\Delta\left(v^2\Psi\right)\|_2 
+\|\Delta\left(v\Phi\right)\|_2\right].
 \end {equation}   
 The estimation of this proceeds as before, but with one new feature:
Now there is a contribution of the form 
$$\Delta v^2 = 2v\Delta v + 2|\nabla v|^2\ .$$
As before, we can use the bound $\|v\|_\infty \le 2$ to conclude that
$\|v\Delta v\|_2 \le 2 \|\Delta v\|_2$. However
$$\||\nabla v|^2\|_2 =   \|\nabla v\|_4^{2}\ .$$
Thus,  using the elementary estimate $\|\nabla v\|_2 \le \|v\|_2 + \|\Delta v\|_2$, we bound the quantity in
\eqref {evo21} by
 \begin{equation}  \label {evo22} \begin {split}  &
-2\| (-\Delta)^{3/2}v\|_2^2 +    C\| (-\Delta)^{3/2}v\|_2 \left[
\|\Delta ( (1-\bar m^2)J*\nabla v)\|_2 
+\|\Delta\left(v^2\Psi\right)\|_2 
+\|\Delta\left(v\Phi\right)\|_2\right] \le \cr & 
-2\| (-\Delta)^{3/2}v\|_2^2 +  C  \| (-\Delta)^{3/2}v\|_2 \left(\|v\|_2  + \|\Delta v\|_2 + 
\||\nabla v|^2\|_2\right)\ . 
\end {split}  \end {equation}  
 To handle $\||\nabla v|^2\|_2$, we compute
\begin{equation*} 
\begin {split}    & 
\int |\nabla v|^4{\rm d}x 
 =  \int \nabla v \cdot \nabla v |\nabla v|^2 {\rm d}x 
 \cr & =  -\int v[(\Delta v)|\nabla v|^2 - 2 D^2v(\nabla v,\nabla v)]{\rm d}x  \cr &
 \le  C\|D^2v\|_2\||\nabla v|^2\|_2 
 \le  C\|\Delta v\|_2\left(\int |\nabla v|^4{\rm d}x\right)^{1/2}\ .
 \end {split}  
  \end {equation*}  
  That is, 
\begin{equation}  \label {evo29} \|\nabla v\|_4^{2}  \le C\|\Delta v\|_2\ . 
  \end {equation}   
 Using this in \eqref {evo22}, our estimate for right hand side of 
\eqref {evo20} becomes 
$$ \frac{d }{d t} \| \Delta v\|_2^2 \le -2\| (-\Delta)^{3/2}v\|_2^2+ C  \| (-\Delta)^{3/2}v\|_2(\|v\|_2 + \|\Delta v\|_2).$$ 
 Now by Schwarz,
\begin{equation}  \label {evo25}
\|\Delta v\|^2_2 = \langle (-\Delta)^{3/2}v, (-\Delta)^{1/2}v\rangle \le 
\| (-\Delta)^{3/2}v\|_2\|\nabla v\|_2. 
 \end {equation}  
Therefore
$$\frac d {dt} \| \Delta v\|_2^2 \le -2 \| (-\Delta)^{3/2}v\|_2 +
\| (-\Delta)^{3/2}v\|_2C\left(\| (-\Delta)^{3/2}v\|_2^{1/2}\|\nabla v\|_2^{1/2} + \|v\|_2\right)\ .$$
By the same type of arithmetic-geometric mean argument we made earlier, we obtain
$$\frac d {dt} \| \Delta v\|_2^2 \le - \| (-\Delta)^{3/2}v\|_2^2 +
C\left(\|\nabla v\|_2^2 + \|v\|_2^2\right)\ ,$$
 
and then by \eqref {evo25},
\begin{equation*}  
  \begin {split} & 
\frac d {dt}  \| \Delta v\|_2^2 
 \le  - \frac {\|\Delta v\|^4_2} { \|\nabla v\|_2^2}+ C\left(\|\nabla v\|_2^2 + \|v\|_2^2\right)
 \cr & 
 \le  -\frac {\|\Delta v\|^4_2} {\delta }+ C\delta. 
 \end {split} 
  \end {equation*}   
The analysis of this differential inequality proceeds exactly as with \eqref {difineq}  in the previous lemma.  \qed
  
\medskip

Up to this point, our analysis had not depended in any significant way on the dimension. To proceed to higher smoothness estimates, we need to take the dimension into account: One last new feature enters in adapting our strategy for proving smoothness to higher derivatives.

\begin{lm}    For dimension $D \le 3$, let $v$ be a solution of \eqref {vevo} and suppose that for some 
$  \delta>0$ and $T_{\delta}>0$,
$$\| v(t)\|_2^2 + \| \nabla v(t)\|_2^2  + \| \Delta v(t)\|_2^2 \le \delta  \qquad{\rm for \ all}\quad t \le T_\delta\ .$$
Then
$$\| \nabla \Delta v(t)\|_2^2 \le \frac{ \delta }{2t} + C\delta \qquad{\rm for \ all}\quad t \le T_{\delta}$$
where $C$ is a constant depending   on $J$, $\beta$ and the dimension.
 \end{lm}

\noindent{\bf Proof:}   
\begin{equation*}  
\begin {split} &
\frac{d }{d t}\| \nabla \Delta v\|_2^2 
=-2\int (-\Delta)^{3}v{\partial v\over \partial t}{\rm d}x   = 
-2\| (-\Delta)^{2}v\|_2^2\cr &
 + 2\beta\int [(-\Delta)^{2}v]\Delta{\rm div}\left((1-\bar m^2)\nabla J* v\right){\rm d}x  
\cr &+2 \int [(-\Delta)^{2}v] \left[
\Delta{\rm div}\left(v^2\Psi\right)  +  \Delta{\rm div}\left(v\Phi\right) \right]{\rm d}x
 \end {split} 
  \end {equation*}  
 and again by the Schwarz inequality, this is no more than
\begin {equation*}  
-2\| (-\Delta)^{2}v\|_2^2 + 2\beta \| (-\Delta)^{2}v\|_2\left[
\|\Delta {\rm div}((1-\bar m^2)J*\nabla v)\|_2 
+\|\Delta {\rm div}\left(v^2\Psi\right)\|_2 
+\|\Delta {\rm div}\left(v\Phi\right)\|_2\right] .
 \end {equation*} 
 As before, in estimating $\|\Delta {\rm div}((1-\bar m^2)J*\nabla v)\|_2$, we may let all derivatives fall on $J$ to obtain the bound $C\|v\|_2$. Also, by  \eqref {evo3}, we have that 
$$\|\Delta {\rm div}\left(v\Phi\right)\|_2 \le 
C\left[ \|v\|_2+ \|\nabla v\|_2 + \|\Delta v\|_2 + \|\nabla\Delta v\|_2\right]\ .$$
However, to estimate $\|\Delta {\rm div}\left(v^2\Psi\right)\|_2$, we need a bound on 
$$\|v^2\|_2 +  \|v|\nabla v|\|_2 +   \| |\nabla v|^2\|_2 +  \|v |\nabla \Delta v|\|_2 + \| |\nabla v| |\Delta v| \|_2 \ $$
for the first, second and fourth terms we may use the {\it a-priori} bound $\|v\|_\infty \le 2$. 
We have already estimated the third term in \eqref {evo29}. The term that forces us to make dimension dependent
estimates is $\| |\nabla v| |\Delta v| \|_2$. The strategy that led to \eqref {evo29}does not work here.
Instead, we use the Sobolev embedding estimate
\begin{equation*}  
\|\nabla v\|_\infty \le C(\|\nabla\Delta v\|_2 + \|\nabla v\|_2)\ , 
 \end {equation*}  
 valid in dimensions $2$ and $3$. 
We then use the bound $\|\Delta v\|_2 \le \delta \le 1$ to obtain
$$\|\Delta {\rm div}\left(v^2\Psi\right)\|_2 \le 
C\left[ \|v\|_2+ \|\nabla v\|_2 + \|\Delta v\|_2 + \|\nabla\Delta v\|_2\right]\ .$$
Finally, we estimate
$$\|\nabla \Delta v\|_2^2 = -\int (-\Delta)^{2}v  \Delta v d x \le \|(-\Delta)^{2} v\|_2\|\Delta v\|_2 \le \delta 
\|\Delta^2 v\|_2$$
to obtain 
$$\frac{ d }{ d t}\| \nabla \Delta v\|_2^2 \le  
-2\| (-\Delta)^{2}v\|_2^2  + C[ \delta +  \| (-\Delta)^{2}v\|_2]\ .$$
The  analysis of this differential inequality proceeds as before. \qed

We now come to the general case.

\begin{lm}\label{lm4}         For dimension $D \le 3$, let $v$ be a solution of \eqref {vevo} and suppose that for some 
$  \delta>0$ and $T_{\delta}>0$, and $k\in \N$
$$\sum_{j=0}^k \| (-\Delta)^{j /2}v(t)\|_2^2  \le \delta  \qquad{\rm for \ all}\quad t \le T_\delta\ .$$
Then
$$\| (-\Delta)^{(k+1)/2}v(t)\|_2^2 \le \frac{ \delta }{2t} + C\delta \qquad{\rm for \ all}\quad t \le T_{\delta}$$
where $C$ is a constant depending   on $J$, $\beta$  and the dimension.
 \end{lm}

\noindent{\bf Proof:} 
\begin{equation*}  
 \begin {split} &  
\frac{ d }{ d t}\| (-\Delta)^{(k+1)/2} v\|_2^2 
=-2\int (-\Delta)^{k+1}v{\partial v\over \partial t}{\rm d}x  = 
-2\| (-\Delta)^{(k+2)/2}v\|_2^2 \cr &
 +  2\beta\int [(-\Delta)^{(k+2)/2}v](-\Delta)^{k/2}{\rm div}\left((1-\bar m^2)\nabla J* v\right){\rm d}x \cr &
 +2 \int [(-\Delta)^{(k+2)/2}v] \left[
(-\Delta)^{k/2}{\rm div}\left(v^2\Psi\right)  +  (-\Delta)^{k/2}{\rm div}\left(v\Phi\right) \right]{\rm d}x . 
 \end {split}
   \end {equation*}  
 As above, we have
$$\|(-\Delta)^{k/2}{\rm div}\left((1-\bar m^2)\nabla J* v\right)\|_2 \le C\delta
$$
and
$$
\|(-\Delta)^{k/2}{\rm div}\left(v\Phi\right)\|_2 \le C\delta\ .$$
Also,
$$
\|(-\Delta)^{k/2}{\rm div}\left(v^2\Psi\right)\|_2 \le C\sum_{j+\ell \le  k+1\ , j,k\ \geq 0} \| |(-\Delta)^{j/2}v|  |(-\Delta)^{\ell/2}v| \|_2\ .$$
Since we have already proved the result for $k\le 3$, we may suppose that $k+1\geq 4$.  
For $k+1\geq 4$, whenever two non-negative integers $j$ and $\ell$ satisfy $j+\ell \le k+1$, at least one of the integers is no 
greater than $k-1$.  Let us suppose that $j \le k-1$. Then we have the sobolev embedding inequality
$$\| (-\Delta)^{j/2}v\|_\infty \le C\left(\| (-\Delta)^{(k+1)/2}v\|_2 +  \| (-\Delta)^{j/2}v\|_2\right)\ ,$$
valid in dimensions $2$ and $3$,  Thus, using the fact that $\delta \le 1$, we have
$$\|(-\Delta)^{k/2}{\rm div}\left(v^2\Psi\right)\|_2 \le C\left(\delta +   \| (-\Delta)^{(k+1)/2}v\|_2\right)\ .$$
Next, by Schwarz
$$\|\| (-\Delta)^{(k+1)/2}v\|_2^2  \le    \|\| (-\Delta)^{(k+2)/2}v\|_2 \|\| (-\Delta)^{k/2}v\|_2^2\ ,$$
and so
\begin{equation*}    \begin {split} &  
\frac{d }{d t}\| (-\Delta)^{(k+1)/2} v\|_2^2  \le  
-2\| (-\Delta)^{(k+2)/2}v\|_2^2 + C\left(\delta +  \| (-\Delta)^{(k+2)/2}v\|_2\right) \cr &
\le  -\| (-\Delta)^{(k+2)/2}v\|_2^2 + C \delta  
 \le  -\frac{\| (-\Delta)^{(k+1)/2}v\|_2^4}
{\| (-\Delta)^{k/2}v\|_2^2} 
+ C\delta \cr &
 \le  -\frac{\| (-\Delta)^{(k+1)/2}v\|_2^4}
{ \delta} 
+ C\delta. 
 \end {split}  
 \end {equation*}  
Thus, $x(t) :=   \| (-\Delta)^{(k+1)/2} v\|_2^2$ satisfies the differential inequality
\eqref {difineq}, and the result now follows as in the proof of  Lemma \ref  {lm1}.\qed 

\noindent  We are now ready to prove the Theorem \ref {smoo}. 
 
 \medskip
 \noindent{\bf Proof of Theorem  \ref {smoo}:}  We proceed by induction on $k$. We shall first show that with $b$ kept fixed,
for any $\epsilon>0$, if $\delta$ is sufficiently small, then
 for any $t_0 > 0$,  and any $k\in \N$, there exists $T_0$  so that
$$\|(-\Delta)^{k/2)}v(t)\|_2^2 \le \epsilon \quad{\rm for\ all}\quad
t_0 \le t \le T_0\ .$$
For $k=1$, this result follows from 
 Lemma \ref  {lm1}.   Suppose that $k\geq2$, and the result has been proved for $k-1$. Then by Lemma \ref  {lm4}
 we have this result for $k$ as well.
 Now, recall that $a(t)$ is defined by
 $$\|m(t) - \bar m_{a(t)}\|_2 \le \|m(t) - \bar m_b\|_2\quad{\rm for\ all}\quad b\in \R\ .$$
 As we have shown, for $\delta$ small enough, this uniquely determines $a(t)$.
 Moreover, as long as $\|m(t) - \bar m_{a(0)}\|_2$ is small, so is $a(t) - a(0)$.
 Then, in the notation of this section,
 $$m(t) - \bar m_{a(t)} =  [m(t)  - \bar m_{a(0)}] + [\bar m_{a(0)} - \bar m_{a(t)}] = v(t) + 
 [\bar m_{a(0)} - \bar m_{a(t)}]\ .$$
 Thus
 $$\|(-\Delta)^{k/2)}(m(t) - \bar m_{a(t)})\|_2 \le   \|(-\Delta)^{k/2)}v(t)\|_2^2 +  
 \|(-\Delta)^{k/2)} [\bar m_{a(t)} - \bar m_{a(0)}]\|_2^2\ .$$
 Next note that
 $$ \|(-\Delta)^{k/2)} [\bar m_{a(t)} - \bar m_{a(0)}]\|_2^2 \le C_k |a(t) - a(0)|\ .$$
 (Note: The constant $C_k$ contains a multiple of $L^d$, so the constant also depend on $L$, which is fixed. This is the
 first place $L$ enters.)
 By  Theorem \ref {5}, $t\mapsto a(t)$ is Lipschitz, and so for any $\epsilon>0$, there is an
 $s_\epsilon$ so that $C_k |a(t) - a(0)| \le \epsilon/2$ provided $t \le s_\epsilon$. 
 Then, by what has been proved above, for any $t_0 < s_\epsilon/2$, there is a $\delta>0$   
 so that if $\|v(t)\|_2 \le \delta$ for all $0\le t \le T_0$, then 
 $ \|(-\Delta)^{k/2)}v(t)\|_2^2 \le \epsilon/2$. Combining results, we have that
 $$ \|(-\Delta)^{k/2)} [m(t) - \bar m_{a(t)}]\|_2^2 \le \epsilon$$
 for all $t_0 \le t \le s_\epsilon$. 
 The same analysis shows that
 $$ \|(-\Delta)^{k/2} [m(t) - \bar m_{a(t)}]\|_2^2 \le \epsilon$$
 for all $s_{\epsilon/2} \le t \le \min \{3s_{\epsilon/2},T_0   \} $,
and inductively,
$$ \|(-\Delta)^{k/2} [m(t) - \bar m_{a(t)}]\|_2^2 \le \epsilon$$
 for all $js_{\epsilon/2} \le t \le \min\{(j+2)s_{\epsilon/2},T_0\}$.
 That is, in steps of fixed length $s_\epsilon$ we cover the interval $(t_0,T_0)$.

 Finally, we prove the assertion in the  Theorem concerning moments.  Namely, in the notation being used here, it suffices to show that
 that for some constant $C$, which may be made small by choosing $\|v_0\|_2$ small,
 $$\int_{\R\times \Lambda} x_1^2v^2(t)\dd x \leq e^{Ct}  \int_{\R\times \Lambda} x_1^2v^2(0)\dd x\ .$$    
 
 Differentiating the left side, we find
\begin{eqnarray}
\frac{{\rm d}}{{\rm d}t} 
\int_{\R\times \Lambda} x_1^2v^2(t)\dd x  &=&
2\int_{\R\times \Lambda} x_1^2v(t)\frac{\partial}{\partial t} v(t)\dd x\nonumber\\
&=&2\int_{\R\times \Lambda} x_1^2v\nabla \cdot 
\left[   
\left(\nabla v - \beta(1-\bar m^2)\nabla J*v\right)  +
v^2\Psi+  v\Phi  
\right]\dd x\nonumber\\
&=& - 2\int_{\R\times \Lambda} x_1^2 v^2 \dd x 
\label{pent1}\\
&+& 2\int_{\R\times \Lambda} x_1^2
\left[   
 \nabla v\cdot \beta\left(1-\bar m^2)\nabla J*v\right)  -
v^2\nabla v\cdot \Psi -  v\nabla v\cdot \Phi  
\right]\dd x\label{pent2}\\
&-& 4 \int_{\R\times \Lambda} x_1v 
\left[   
\left(\nabla_1 v - \beta(1-\bar m^2)\nabla_1 J*v\right)  +
v^2\Psi_1+  v\Phi_1  
\right]\dd x\ .\label{pent3}
\end{eqnarray} 
The term in (\ref{pent1}) has a sign that allows us to ignore it. The next simplest term is the integral in (\ref{pent3}). 
Using the Scwarz inequality, we may bound it in magnitude by
$$\left(\int_{\R\times \Lambda} x_1^2v^2\right)^{1/2} 
\left(\int_{\R\times \Lambda} \left[   
\left(\nabla_1 v - \beta(1-\bar m^2)\nabla_1 J*v\right)  +
v^2\Psi_1+  v\Phi_1  
\right]^2\dd x \right)^{1/2}\ .$$
By what we have proved above, the second square root on the right is bounded uniformly (and small) on the interval under consideration.  

After one more integration by parts in the variable $x_1$, the contribution in (\ref{pent2}) is handled the same way.  
 \qed

 \vskip1.cm
  \section {Appendix  }

  \begin{lm} \label {A2}      Let $m \in \cal M$,  see \eqref {bigmdef},  
 $m =\bar m+v$, where $\bar m$ is the closest instanton  to  $m$  in $L^2 (\R \times \L)$.
 There exists $\k>0 $,   $  \delta >0$,  $c=  c(\k) >0$ so that  for
  and  
$\|v\|_ {W^{s,2}}  \le  \k  $, where     $s >   \frac {D} 2$,     
we have    
\begin{equation}  \label {2.41} \frac 14 \g(L)  ||m -\bar m  ||_2^2 \le {\cal F}\big( m  \big)-{\cal
F}\big(\bar m  
\big) \le   c  ||m -\bar m    ||_2^2.
 \end {equation}  
 Moreover for any $\e>0$  there is a $ \tilde \k (\e,
L,
\b)
$  so that
\begin{equation}  \label {A.1} \frac {1-\e}2 \langle v, \BB v\rangle_{L^2 (\R \times \L)}  \le  {\cal F}\big( m  \big)-{\cal
F}\big(\bar m 
\big) \le \frac {1+\e}2 \langle v, \BB v\rangle_{L^2 (\R \times \L)}  
 \end {equation}  
provided  
 $\|v\|_ {W^{s,2}}  \le  \tilde  \k $     
for    $s >   \frac {D} 2$.  
\end {lm}
 
\noindent{\bf Proof:}   Denote by ${\cal F} '\big (m \big)(v) $ and  ${\cal F}''\big(m \big)  (v, w)$  respectively  the first and the second Frechet derivative of   ${\cal F}(\cdot)$ computed at $m\in \cal M$ in the directions $v$ and $w$ in  $L^2(\R \times \L)$.
 It is easy to see that for $ \|m \|_\infty \le c<1$ repeated Frechet derivatives exist with 
$$   {\cal F} '\big (m \big)(v) = \int_{\R \times \L}   \left [ {1\over\beta}{\rm arctanh} m (x)  -
J\star  m (x) \right ] v(x)  {\rm d} x,$$
and
$$   {\cal F} ''\big (m \big) ( v,w) = \int_{\R \times \L}   \left [ {1\over\beta}\frac {v(x) } { (1 - m^2(x)) } -
J\star  v (x) \right ]  w(x)  {\rm d} x.
$$
When $m = \bar m$ then 
  $$   {\cal F} ''\big (\bar m \big) ( v,w) = \int_{\R \times \L}   \left [ {1\over\beta}\frac {v(x) } { (1 - \bar m^2(x)) } 
-J\star  v (x) \right ]  w(x)  {\rm d} x=    \langle \BB v, w\rangle,$$
where $ \BB$ is the operator defined in \eqref {G.2}.

 Writing $m = \bar m +v $,  we can represent  
$$  {\cal F}\big(\bar m +v \big)-{\cal F}\big(\bar m \big) =\int_0^1 d \tau
{\cal F} '\big (\md + \tau  v  \big)(v)    
 =\int_0^1 d \tau\int_0^\tau  ds  {\cal F}''\big(\bar  m +
s v \big)  ( v, v ) . 
$$
   In order to get a lower bound for the last term  above we expand 
${\cal F}''\big(\bar  m + s v \big) ( v,v ) $  around $s =0$
obtaining 
$${\cal F}''\big(\bar  m + s(m-
\md)\big)  (v, v) 
 ={\cal F}''\big(\bar  m) (v, v)+  {\cal F}'''\big(\tilde
m  \big)  (v,v,v)$$ 
where $\tilde m = \bar  m + s_0v  $ for some $s_0$ between $0$ and
$1$ by the mean value theorem. Therefore 
\begin{equation}  \label {A.2}  {\cal F}\big(m \big)-{\cal F}\big(\bar m \big) = \frac 12 \langle \BB
v, v \rangle + \int_0^1 d \tau\int_0^\tau  ds  {\cal F}'''\big(\tilde m 
\big)(v,v,v).  
 \end {equation}  
Since $\bar m$ is the closest instanton to $m$ in $L^2(\R\times \L)$,  $\int \bar m' (x) v(x) {\rm d} x =0$.  Therefore by  \eqref {G.3}
\begin{equation}  \label {P.1} \langle \BB v, v\rangle
 \ge \g(L) \|v\|_2^2.
  \end {equation}  
We  then need
a lower bound on the term involving the third derivative of the free
energy.  By direct computation,
\begin{equation}  \label {3.6}
   \big |{\cal F}'''\big(\tilde m  \big)(v,v,v)\big|
     = {2\over \beta} 
\big|\int_{\R \times \L }   { \tilde m\over (1-\tilde
m^2)^2}\big(v(x)\big)^3{\rm d}x \big|.
\end {equation}  
 Take  $\|v\|_ {W^{s,2}}  \le \d_1$,     so  that $||\tilde m||_\infty
\le 1- \d_0$, with $\d_0>0$, see Lemma \ref {A1}. 
With this choice of $\d_1$ we have
 $$
   \big |{\cal F}'''\big(\tilde m  \big)(v,v,v)\big|
\le  c (\beta, \d_1) \int_{\R \times \L }   |v(x)|^3{\rm d}x   $$ for
some constant $ c (\beta, \d_1)  $ depending   on $\beta$ and $\d_1$. 
     We have that 
$$  \int_{\R \times \L }   |v(x)|^3{\rm d}x   \le \sup|v(x)| \int_{\R \times \L }   |v(x)|^2 $$
and, see Lemma \eqref {A1},  for $s>\frac {D} 2$ 
$$ \sup|v(x)| \le C(d,s)  \|v\|_ {W^{s,2}} .$$ 
 Therefore 
\begin{equation}  \label {A.3}  \big |{\cal F}'''\big(\tilde m  \big)<v,v,v>\big| \le  
 c (\beta, \d_1,d)  \|v\|_ {W^{s,2}}  \bigl(\|v\|_2 \bigr)^2
  \end {equation}  
and, see \eqref {A.2},
 \begin{equation*}  
   {\cal F}\big(m \big)-{\cal F}\big(\bar m \big) \ge   \|v\|_2^2  \left [   \frac 12 \g(L) -  c (\beta, \d_1,d)   \|v\|_ {W^{s,2}} \right] . 
   \end {equation*}  
Taking     $ \d:=  \min \{ \d_1,  \d (\b, L) $   so that
 \begin{equation*}  
  \frac 14 \g(L) -  c (\beta, \d_1,d)   \|v\|_ {W^{s,2}} \ge 0 
  \end {equation*}  
we have 
  \begin{equation*}    
  {\cal F}\big(m \big)-{\cal F}\big(\bar m \big) \ge  \frac 14 \g(L) \|v\|_2^2.  
   \end {equation*}  
Thus,  we have established a lower bound for  \eqref {2.41}. 
   The upper bound follows from the boundedness of $\BB$ and  estimate \eqref {3.6} of  ${\cal F}'''\big(\tilde m  \big)(v,v,v)$.    Note that one needs always a bound on  $\|v\|_ {W^{s,2}}$  to  get  $\|\tilde m\| \le 1- \delta_0$.  In this way
we proved \eqref {2.41}.

In a similar way  the inequalities \eqref {A.1}follows.
Namely,  from \eqref {A.2},  for any positive $\e$   
\begin{equation}  \label {A.20}   {\cal F}\big(m \big)-{\cal F}\big(\bar m \big) =
\frac 12(1-\e) \langle \BB v, v\rangle  + \left [ \frac 12 \e  \langle \BB v, v\rangle   + 
\int_0^1 d
\tau\int_0^\tau  ds {\cal F}'''\big(\tilde m 
\big)(v,v,v) \right ].  
 \end {equation}  
{ }  From \eqref {G.3} and    \eqref {A.3},  denoting   $\tilde \d$ the  $\d$ appearing in the formula,    the last term in \eqref {A.20} is
bigger or equal to
\begin{equation}  \label {A.21} 
\frac 12 \e \g (L)  ||v||_2^2 -  c (\beta, \tilde \d)  C(d,s) \|v\|_ {W^{s,2}}  \bigl(\|v\|_2 \bigr)^2
= ||v||_2^2  \left [ \frac 12 \e  \g (L)   -
 c (\beta, \tilde \d)  C(d,s)  \|v\|_ {W^{s,2}} \right ].  
  \end {equation}
Choosing   $\tilde \d  $  so that 
  the  term in   \eqref {A.21} is strictly positive 
 we get the lower bound  \eqref {A.1}.
 The upper bound  \eqref {A.1} follows immediately. 
  \qed

 \bigskip
 \begin{lm} \label {TT2}   Let $\rho(x)$ be a probability density with
$$\int |x|\rho(x){\rm d}x < \infty\ .$$
Then for   $v \in W^{1,2} (\R \times \L) $   
\begin{equation*}   \|v - \rho\star v\|_2 \le  \|\nabla v\|_2\
\int |x|\rho(x){\rm d}x  .
 \end {equation*}
\end {lm}
\bigskip
\noindent{\bf Proof:}   We have 
\begin{equation*}     \begin {split}  &\|v - \rho\star v\|_2^2= \cr
& \int_{\R \time \L} dx \left (  \int_{\R \time \L}  \rho (x-y) [v(y)- v(x) ] dy \right )^2  \cr & \le \left (\int |x|\rho(x){\rm d}x \right )^2\|\nabla v\|_2^2.  
\end {split}
\end {equation*}
 \qed

The next lemma shows that  for any function 
$v$ that is orthogonal to $\bar m'$, whenever  $ \| \nabla  v  \|_{L^2 (\R \times \L)} $   is small compared to $\|v\|_{L^2 (\R \times \L)}$,
then ${\cal B}v$
is very close to being a constant multiple of $v$, $\tilde\alpha v$ where 
$\tilde\alpha$ is defined by
\begin{equation*} 
 \tilde\alpha = {1\over \beta(1-m_\beta^2)} - 1
\end {equation*}
 and it is strictly positive for $\beta>1$. The lemma also shows that under the same
condition, $\sigma(\bar m)v$ is very close to $\sigma(m_\b)v$.

\bigskip
 \begin{lm} \label {TT1} 
  Let $v \in W^{1,2} (\R \times \L)$, $\langle
v,\bar m'\rangle_{L^2} = 0$.   There is a finite positive constant $K(\b,J, L, d)$  so that
\begin{equation}  \label {subs1}\|\BB v - \tilde\alpha v\|_{L^2 (\R \times \L)}  \le K(\b,J, L, d)  \| \nabla  v  \|_{L^2 (\R \times \L)}  ,
\end {equation}
and 
\begin{equation}  \label {subs2}\|\sigma (\bar m)  v -  \sigma (m_\b) v\|_{L^2 (\R \times \L)}  \le K(\b,J, L, d)  \| \nabla  v  \|_{L^2 (\R \times \L)}.
\end {equation}
\end {lm}
\bigskip
 \noindent{\bf Proof:}  
  Clearly, 
$$\BB v - \tilde\alpha v =
{1\over \b}\biggl({\bar m^2 - m_\beta^2\over (1-\bar m^2)(1-m_\beta^2)}
\biggr)v +  (v - J\star v).$$
We will estimate these two terms separately. For the second term we apply Lemma \ref {TT2}.  For the first term  split, as done in \eqref {P1}, 
$$ v (x_1,x^\perp)= v_1(x_1)+ w(x_1,x^\perp).$$
Then
$${1\over \b}\biggl({\bar m^2 - m_\beta^2\over (1-\bar m^2)(1-m_\beta^2)}
\biggr)v= {1\over \b}\biggl({\bar m^2 - m_\beta^2\over (1-\bar m^2)(1-m_\beta^2)}
\biggr)[v_1+w] $$
We estimate 
$${1\over \b}\biggl({\bar m^2 - m_\beta^2\over (1-\bar m^2)(1-m_\beta^2)}
\biggr) v_1  $$
 as in     \cite {CCO2}, being   one dimensional. 
First, for any $y_1$ and $x_1$ we have
$$ v_1(x_1) =  v_1(y_1) + \int_{y_1}^{x_1}  v_1'(z){\rm d}z\ .$$
Now multiply both sides by $\bar m'(y_1)$ and integrate in $y_1$. By the orthogonality
of $\bar m'$ and $v_1$, we have
$$ 2m_\beta v_1(x_1)  = 
\int_{-\infty}^\infty\bar m'(y)\biggl(\int_{y_1}^{x_1} v_1'(z) {\rm d}z\biggr){\rm d}y\ .$$
But $|\int_{y_1}^{x_1} v'(z){\rm d}z| \le |x_1-y_1|^{1/2}\|v'\|_2$ so that
$$|v_1(x_1)| \le 
{1\over 2m_\beta}\biggl(\int\bar m'(y)|x_1-y_1|^{1/2}{\rm d}y\biggr)\|v_1'\|_{L^2(\R)}\ ,$$
and clearly there is a finite constant $K(\b,J)$ depending only on $\beta$
and $J$ so that
$${1\over 2m_\beta}\int\bar m'(y)|x_1-y_1|^{1/2}{\rm d}y \le K(\b,J)(1+|x_1|) ,$$
and hence
\begin{equation}  \label {pntwsbnd} | v_1(x_1)| \le K(\b,J)(1+|x_1|)\|v_1'\|_{L^2(\R)}.
\end {equation}
Next, using the pointwise bounds \eqref {pntwsbnd} established above,
\begin{equation*}  
\begin {split} 
&\|(\bar m^2 - m_\beta^2)\bigl((1-\bar m^2)(1-m_\beta^2)\bigr)^{-1}v_1\|_{L^2(\R\times \L)}^2  \cr \le
&\|v_1'\|_{L^2(\R \times \L)}^2K(\b,J)\int (1+|x_1|)^2(\bar m^2 - m_\beta^2)^2
\bigl((1-\bar m^2)(1-m_\beta^2)\bigr)^{-2}
{\rm d}x_1  \cr \le
&\tilde K(\b,J)\|v_1'\|_{L^2(\R \times \L)}^2,
\end {split} 
\end {equation*}
where $\tilde K(\b,J)$ is finite by the rapid decay of $(\bar m^2 - m_\beta^2)^2$.
Further 
$$ \|{1\over \b}\biggl({\bar m^2 - m_\beta^2\over (1-\bar m^2)(1-m_\beta^2)}
\biggr) w\|^2_{L^2(\R \times \L)}     \le  K(\b,J) \|w \|^2_{L^2(\R \times \L)} . $$
Applying the Poincar\'e  inequality as in \eqref {Pa.10} we have 
$$ \|w\|^2_{L^2(\R \times \L)} \le L^2 c(d) \| \nabla^{\perp} w\|^2_{L^2(\R \times \L)}. $$ 
Then  
\begin{equation*} 
 \begin {split}
  & \| {1\over \b}\biggl({\bar m^2 - m_\beta^2\over (1-\bar m^2)(1-m_\beta^2)}
\biggr)v\|^2_{L^2(\R \times \L)} \le    K(\b,J)  \left [ \|v_1'\|_{L^2(\R \times \L)}^2 + L^2 c(d) \| \nabla^{\perp} w\|^2_{L^2(\R \times \L)}\right ] \cr & \le
K(\b,J, d,L) \| \nabla v\|^2_{L^2(\R \times \L)}.  
\end {split}
 \end {equation*}  
  The proof of \eqref {subs2} is very similar to the proof of \eqref {subs1}.
  \qed

\medskip
For function 
 $  v \in  W^{s,2}(\R \times \L)$ we have the following result which can be proven by Fourier analysis, see  \cite {GT}.  

  \vskip0.5cm    
\begin{lm} \label{A1} For   $v \in     W^{s,2}(\R \times \L)$,  if $s > \frac { D} 2$,  we have 
  $$\|v \|_\infty \le  C(d,s) \|v\|_ {W^{s,2}} .$$
\end {lm}
\vskip0.5cm 
 
 \begin{lm} \label {A10}  
For any real number  $a$ and $b$  and for  any  $\lambda$,  $0<\lambda
<1$, 
$$ (a+b)^2  \ge   \l a^2 - \biggl({1\over
1-\l }-1\biggr)b^2.  $$
\end {lm}
\noindent{\bf Proof:}   The proof is immediate: 
\begin{equation}  \label {3.B1}  \begin {split}
(a+b)^2  &\ge a^2 + b^2 - 2ab\cr &=
\l a^2 + \biggl((1-\l )a^2 +b^2 -
2ab\biggr)\cr &\ge \l a^2 - \biggl({1\over
1-\l }-1\biggr)b^2.
 \end {split}  \end {equation}  
The last inequality is obtained adding and subtracting $\frac 1 {1+\lambda} b^2$.


\begin{thebibliography}{99}

\bibitem{AB}   G. Alberti,   G. Bellettini {\it A nonlocal anisotropic model for phase transitions. I. The optimal profile problem.} Math. Ann 310,no. 3,  527-560 (1998)  


\bibitem{BDDP} G. Bellettini, A. De Masi,, N. Dirr,  E. Presutti {\it  
Tunnelling in two dimensions.}
Commun. Math. Phys. 269, 715-763 (2007)


\bibitem{CCO1} E. A. Carlen, M. C. Carvalho,  E. Orlandi, {\it 
Algebraic rate of decay for the excess free energy and stability of
fronts for a non--local phase kinetics equation with a conservation
law I 
}
  J. Stat. Phy.   {\bf 95}  N 5/6,  1069-1117  (1999) 
 
 \bibitem{CCO2}  E. A. Carlen, M. C. Carvalho, E. Orlandi 
{\it  Algebraic rate of decay for the excess free energy and stability of fronts 
for a no--local phase kinetics equation with a conservation law, II}
  Comm. Par. Diff. Eq.  {\bf 25}    N 5/6,  847-886  (2000) 
 
 \bibitem{CCO3} E. A. Carlen, M. C. Carvalho, E. Orlandi 
{\it A simple proof  of stability of fronts for the Cahn-Hilliard
equation}
 Comm. Math. Phy  224-no 1,   323-340  (2001) 
  


\bibitem{DGP}A. De Masi, T.  Gobron, E. Presutti, 
   {\it Travelling fronts in non local evolution equations}
  Arch.  Ration. Mech. Anal. {\bf 132}, 143-205 (1995)
 
 \bibitem{DOPT3}A. De Masi, E. Orlandi, E. Presutti,  L.
 Triolo, {\it  Stability of the interface in a model of phase
 separation}  Proc.Royal Soc. Edinburgh  {\bf 124A} 1013-1022 (1994).

  
\bibitem{DOPT4}A. De Masi, E. Orlandi, E. Presutti,  L.
Triolo, {\it Uniqueness of the instanton profile and global
stability in non local evolution equations}
Rendiconti di Matematica. Serie Vll {\bf 14}, (1994).
 
\bibitem{GL1} G. Giacomin,  J. Lebowitz, 
{\it  Phase segregation dynamics in particle systems with long range
interactions I: macroscopic limits} J. Stat. Phys.  {\bf
87}, 37--61 (1997).

\bibitem{GT}  
D. Gilbarg, N.S. Trudinger, 
{\it Elliptic partial Differential Equations of Second order}
Grundlehren der mathematischen Wissenschaften, 2nd ed., Berlin 1983.

   
\bibitem{KKT}T. Korvola, A. Kupiainen, J. Taskinen, {\it  
Anomalous scaling for 3d Cahn-Hilliard Fronts}
Comm. Pure Appl. Math. 58, no 8, 1077-1115  (2005)

\bibitem{LOP} J.L. Lebowitz,
 E. Orlandi, E. Presutti,  {\it A Particle model for spinodal
decomposition} J. Stat. Phys.  {\bf 63}, 933-974 (1991) 
 
 
\bibitem{Pr}  E. Presutti,    { \it Scaling Limits In Statistical Mechanics and Microstructures in Continuuum Mechanics}  Springer (2009)







\end{thebibliography}
\end{document}